\documentclass[12pt]{iopart}

\usepackage{iopams}  

\usepackage{amssymb}

\usepackage{comment}
\usepackage{graphicx}
\usepackage[table]{xcolor}
\usepackage{todonotes}
\usepackage{placeins}
\usepackage{enumitem}
\usepackage[utf8]{inputenc}
\usepackage{hhline}
\usepackage[english]{babel}

\usepackage[notcite,notref,color]{showkeys}
\usepackage{hyperref}

\DeclareRobustCommand{\DE}[3]{#2} 
\DeclareRobustCommand{\VANDER}[3]{#2} 

\usepackage[displaymath, mathlines]{lineno}

\newcommand{\bR}{\mathbb{R}}

\newcommand{\cA}{\mathcal{A}}

\newcommand{\cC}{\mathcal{C}}
\newcommand{\cD}{\mathcal{D}}
\newcommand{\cE}{\mathcal{E}}
\newcommand{\cF}{\mathcal{F}}
\newcommand{\cG}{\mathcal{G}}

\newcommand{\cM}{\mathcal{M}}
\newcommand{\cN}{\mathcal{N}}
\newcommand{\cP}{\mathcal{P}}
\newcommand{\cU}{\mathcal{U}}

\newtheorem{conj}{Conjecture}
\newtheorem{thm}{Theorem}

\newtheorem{prop}{Proposition}
\newtheorem{lemma}{Lemma}
\newtheorem{definition}{Definition}

\newcommand{\proof}{\noindent{\bf Proof:}~}
\newcommand{\qed}{\hfill $\Box$\newline~}
\newcommand{\sign}{\mathop{{\rm sign}}}

\begin{document}

\title[Bifurcations of critical sets and relaxation oscillations]{Bifurcation of critical sets and relaxation oscillations in singular fast-slow systems}

\author{Karl Nyman$^{1}$, Peter Ashwin$^{2}$ and Peter Ditlevsen$^{1}$}
\address{1. Niels Bohr Institute, University of Copenhagen, 2100 Copenhagen, Denmark\\
2. Department of Mathematics, University of Exeter, Exeter EX4 4QF, UK}

\ead{karl.nyman@nbi.ku.dk}
\vspace{10pt}
\begin{indented}
\item[]\date{January 2019}
\end{indented}

\begin{abstract}
Fast-slow dynamical systems have subsystems that evolve on vastly different time\-scales, and bifurcations in such systems can arise due to changes in any or all subsystems. We classify bifurcations of the critical set (the equilibria of the fast subsystem) and associated fast dynamics, parametrized by the slow variables. Using a distinguished parameter approach we are able to classify bifurcations for one fast and one slow variable. Some of these bifurcations are associated with the critical set losing manifold structure. We also conjecture a list of generic bifurcations of the critical set for one fast and two slow variables. We further consider how the bifurcations of the critical set can be associated with generic bifurcations of attracting relaxation oscillations under an appropriate singular notion of equivalence.
\end{abstract}

%
\vspace{2pc}

%
%
%


\section{Introduction}

Many natural systems are characterized by interactions between dynamical processes that run at very different timescales. These can often be modelled as fast-slow systems, where system dynamics can be separated into interacting fast and slowly changing variables. This has been applied to a wide range of natural phenomena, including electrical circuits \cite{vanderpol}, plasma oscillations\cite{mikikian:2008}, 
surface chemistry \cite{krischer:1992} and cell physiology \cite{harvey:2011} to ecology\cite{rinaldi:2000} and climate\cite{ashwinditlevsen15}. The dynamical behavior of such systems can often be understood in a common mathematical framework.
See \cite{kuehn15} for a recent monograph that summarizes both techniques and applications, and \cite{arnold99,benoit81,guckenheimer03b,ermentrout86,guckenheimer96,guckenheimer03a,guckenheimer12,wechselberger12} for examples of related work. 

The analysis of fast-slow systems is built around a {\em geometric singular perturbation theory} (GSPT) approach \cite{fenichel79}, perturbing from a singular limit where timescales decouple: see also \cite{fenichel71,HirschPughShub1977}. In the singular limit, on the slow timescale there are instantaneous jumps (determined by the fast dynamics) between periods of slow evolution. The slow evolution typically takes place on stable sheets of a {\em critical set} where the fast dynamics is in stable balance, interspersed by fast transitions. If the fast dynamics is one dimensional, then it is typically confined to a manifold (and hence the set is often called a {\em critical manifold}), though at bifurcation it may lose its manifold structure at isolated points. The fast transitions are determined by what we call the {\em umbral map} defined as the map from a fold point on the critical set to another part of the critical set.
In the case of stable periodic behaviour in the singular limit, the resulting limit cycles are referred to as {\em singular relaxation oscillations}. Many examples of bifurcations of relaxation oscillations have been considered \cite{arnold99}, including some associated with bifurcations of the critical set \cite{ashwinditlevsen15, franci12, franci14} although it seems that no exhaustive list of scenarios has been proposed.

Although singular perturbation theory has been developed to explain many aspects of behaviour near the singular limit, there is still no full understanding of generic bifurcations of limit cycles even in the singular limit. Guckenheimer \cite{guckenheimer96,guckenheimer02,guckenheimer04} suggests an approach and several results along these lines, but, as far as we are aware, these conjectures are yet to be framed, let along proved, in rigorous terms. The main aim of this paper is to present an approach to doing precisely this, using singularity theory with distinguished parameters and appropriate notions of equivalence. We classify local and global transitions in the critical set by codimension and consider the consequences for the {\em umbral map}. In doing so, we find a variety of scenarios that give bifurcation of attractors in such singular systems.

We show that it is possible to split the problem of bifurcations of relaxation oscillations into two aspects: bifurcations of the critical set, and bifurcations caused by singularities of the slow flow (possibly interacting with the critical set). In the simplest case of one fast and one slow variable, bifurcations of the critical set can be directly tackled using a global version of the singularity theory with distinguished parameter in \cite{golubitskyscheffer85}. We highlight that this theory needs extension to make it suitable for systems with multiple fast and slow variables. We are able to identify a large number of scenarios that can lead to bifurcation of relaxation oscillations. Note that we only consider fast-slow systems where the fast dynamics is constrained to a subset of the variables; following \cite{fenichel79} it is possible to extend the theory developed here to more general fast-slow systems that are not in standard form: see for example \cite{kuehn2014normal,kuehn2018duck,zagaris2004fast}.

We structure the paper as follows: in Section~\ref{sec:singularlimitintro} we briefly introduce the singular limit of fast-slow systems, critical sets, singular trajectories and global equivalence of critical sets. In Section~\ref{sec:criticalmanifolds} we explore  persistence and bifurcation of critical sets by examining versal unfoldings of the fast dynamics parametrized by the slow variables, using a notion of global equivalence of the fast dynamics. In the case of one fast and one slow variable we classify persistence (Proposition~\ref{prop:CM110}) and bifurcations of the critical set up to codimension one (Proposition~\ref{prop:CM111}) using the theory of \cite{golubitskyscheffer85}. For one fast and several slow variables we highlight the need for an improved theory of bifurcations with multiple distinguished parameters.  We present conjectured statements of persistence and of codimension one bifurcations of the critical set (Conjectures~\ref{conj:CM120} and \ref{conj:CM121} respectively) for one fast and two slow variables. We find a rich variety of distinct mechanisms for typical codimension one bifurcations of the critical set which includes local and global bifurcations in the fast variable.

Section~\ref{sec:bifurcation} turns to the question of bifurcation of attractors in fast-slow systems and in particular bifurcation of stable relaxation oscillations. We introduce a global singular equivalence for the singular trajectories and use this to classify persistence (Proposition~\ref{prop:fastslowpersimpliessropers}) and codimension one bifurcations (Proposition~\ref{prop:srobifcodim1}) of these simple relaxation oscillations. These codimension one bifurcations naturally split into those caused by bifurcations of the critical set, and those caused by interaction of singularities of the slow flow with the critical set: in Section~\ref{sec:examplebifs} we present some numerical examples of various types. Finally, Section~\ref{sec:discussion} is a discussion of some of the challenges for GSPT to describe the unfolding of such bifurcations, and relation to other singularity theory approaches. We include several Appendices that give more details of the tools used for the classification and the examples.

\section{Singular trajectories of fast-slow systems}
\label{sec:singularlimitintro}

A {\em fast-slow system} is a system of coupled ODEs for $z=(x,y)\in\bR^{m+n}$ of the form
\begin{equation}
\left\{\begin{array}{rl}
\epsilon \dot{x} & = g(x,y,\epsilon) \\
\dot{y} & = h(x,y,\epsilon)
\end{array}\right.
\label{eq:mainsystem}
\end{equation}
where $x\in\bR^m$ and $y\in\bR^n$, $\epsilon>0$ is a small constant and $t$ is a time measured relative to the slow dynamics. The functions $g(x,y,\epsilon)$ and $h(x,y,\epsilon)$ are $C^\infty$ functions of their arguments (they are well approximated by Taylor series to arbitrary order). We refer to $x$ as the {\em fast} and $y$ as the {\em slow} subsystems; the dynamics in these subsystems are referred to as {\em fast} and {\em slow} dynamics respectively. These systems have a singular limit $\epsilon\rightarrow 0$, where a typical trajectory can remain close to an equilibrium of the fast system, except at isolated points where it 
``jumps'' along a trajectory of the fast subsystem. The singularly perturbed system with $\epsilon>0$ will have trajectories that typically remain near a trajectory of the singular limit, although especially near bifurcations, trajectories may also explore unstable parts of the slow dynamics along so-called {\em canard solutions} (see e.g. \cite{levinson49, benoit81,arnold99,wechselberger12} and \cite[Chapter 8]{kuehn15}).

In order to understand such systems it is useful to consider the {\em reduced} or {\em slow equations} in slow time $t$:
\begin{equation}
\left\{\begin{array}{rl}
0 & = g(x,y,0) \\
\dot{y} & = h(x,y,0),
\end{array}\right.
\label{eq:reduced}
\end{equation}
describing the slow dynamics in the singular limit $\epsilon \rightarrow 0$. Solutions of (\ref{eq:reduced}) are constrained to the {\em critical set}
\begin{equation}
\cC[g] = \left\{ (x,y) \in \bR^{m+n}: g(x,y,0) = 0 \right\}.
\label{eq:critman}
\end{equation}
Note that by the regular value theorem \cite{loring10}, this critical set is a manifold at all points where the derivative of $g$ has maximal rank. For an open dense set of $g\in C^\infty$ this is true for an open and dense set $(x,y)\in\cC[g]$. The set is often called a critical manifold, however we do not use this notation as at bifurcation points the set may lose its manifold structure.

The flow $g$ may have singular equilibria of the fast dynamics within $\cC[g]$, where the term singular here means non-regular. The {\em regular points} of the critical set $\cC[g]$ we define as
\begin{equation}
\cC_{reg}[g]=\left\{(x,y)\in\cC[g]~:~\partial_xg(x,y,0)\mbox{ is hyperbolic} \right\}.
\end{equation}
The remaining non-hyperbolic (fold) points, also called singular or limit points, form the {\em fold set} of the critical set
\begin{equation}
\cF[g]=\left\{(x,y)\in\cC[g]~:~\partial_xg(x,y,0)\mbox{ is non-hyperbolic} \right\}.
\end{equation}
At all regular points the reduced equations (\ref{eq:reduced}) can be used to describe the flow. At fold points, however, we need to consider the fast dynamics and there will be {\em jumps} in slow time determined by the fast subsystem only.

Changing variable to a {\em fast time} $\tau = t/\epsilon$ and taking the limit $\epsilon\rightarrow 0$ gives quite a different set of equations: the  {\em layer} or {\em fast equations}:
\begin{equation}
\left\{\begin{array}{rl}
x' & = g(x,y,0) \\
y' & = 0,
\end{array}\right.
\label{eq:layer}
\end{equation}
where we write $x'$ to denote $\frac{d}{d\tau}x$ and note that (a) the constant slow variable $y$ now acts as a parameter for evolution of the fast variable $x$ and (b) the layer equation, when restricted to $\cC[g]$ consists entirely of equilibria for $m=1$ (for $m>1$ there may be other objects, such as limit cycles). We split the regular set into a disjoint union of attracting/repelling/saddle points
\begin{equation*}
\cC_{att}[g],~\cC_{rep}[g],~\cC_{sad}[g]
\end{equation*}
so that $\cC_{reg}[g]=\cC_{att}[g]\cup \cC_{rep}[g]\cup \cC_{sad}[g]$. Note that $\cF[g]$ is the union of the set of boundaries of these sets, and also that the set $\cC_{sad}[g]$ only exists for $m\geq 2$. Note that the regular set is the union of all non-singular points
\begin{equation*}
\cC_{reg}[g] = \cC[g] \setminus \cF[g].
\end{equation*}
From now on, we only concern ourselves with the system in the singular limit, and therefore drop the dependency on $\epsilon$ in our notation, such that e.g. $g(x,y,0):=g(x,y)$ and $h(x,y,0):=h(x,y)$. However, we stress once more that the system (\ref{eq:mainsystem}) may depend crucially on $\epsilon$ near, but away from, the singular limit.

We now make the notion of \emph{jumps} more precise: For any point $p=(x,y)\in\cF[g]$ we define the (possibly set-valued) {\em umbral map} to be

\begin{equation*}
U[g](p)=\{\mbox{$\omega$-limits of a non-trivial trajectories in (\ref{eq:layer}) with $\alpha$-limit $p$}\},
\end{equation*}
that is, a map from every point $q$, which limits to $p$ in backward time ($\alpha$-limit), to the set of forward time limits ($\omega$-limit) of every such $q$, excluding $p$ itself. The $\omega$-limits are always non-empty since we will assume bounded global attractors. However, if all $\omega$ limits equal $p$ then the umbral map is empty. 
The {\em umbra} (meaning shadow \cite{guckenheimer96}) or {\em drop set} \cite[Ch.5, p.109]{kuehn15} is the image of the folds under the umbral map:
\begin{equation*}
\cU[g]=U[g](\cF[g]).
\end{equation*}

We define the projection onto the slow variable as
\begin{equation*}
\pi:\bR^{m+n}\rightarrow \bR^n,~~\pi(x,y)=y
\end{equation*}
and for $p=(x,y)$ we define the set of all {\em co-folds} to $p$ as
\begin{equation}
\Pi(p)=\{q\in \cF[g]~:~\pi(p)=\pi(q)\}=\pi^{-1}\left(\pi(p)\right)\cap \cF[g],
\label{eq:cofolds}
\end{equation}
i.e. all fold points with the same slow coordinate as $p$. Similarly, we define the set of folds sharing a given slow ($y$) coordinate to be:
\begin{equation}
P(y) = \cF[g] \cap \pi^{-1}(y).
\label{eq:Pydef}
\end{equation}

\subsection{Trajectories in the singular limit}

Note that typical points in $\bR^{m+n}$ are not on $\cC[g]$: starting at  $(x,y)\not \in\cC[g]$ there will be fast motion governed by the layer equations (\ref{eq:layer}). If this settles to a limit we will typically have arrived at a point on the critical set that is a stable equilibrium, i.e. on $\cC_{att}[g]$. The slow dynamics then carries the trajectory around $\cC_{att}[g]$ until (possibly) it hits a fold point $p=(x,y)\in\cF[g]$. If $U[g](p)$ is a single point then there is a unique non-trivial trajectory of the layer equations from this point, the fast motion will take the dynamics to $U[g](p)$.

Hence typical trajectories in the singular limit (which can be viewed as trajectories of a constrained differential equation \cite{takens76}) are composed of segments of slow trajectories on $\cC_{att}[g]$ interspersed with fast jumps. Depending on the nature of the slow segments and fast jumps, characterised by the umbral maps, there may be a trajectory of the $\epsilon>0$ system that remains close to the singular trajectory. More precisely, we define a singular trajectory as follows (this idea is widely used and sometimes called a candidate trajectory, for example \cite{tikhonov48, levinson49, benoit83,benoit90,maesschalck11} and \cite[Ch.3, p.64]{kuehn15}):

\begin{definition}[Singular trajectory]
\label{def:singulartraj}
A singular trajectory is a homeomorphic image $\gamma_0(s)$ of a real interval $[a,b]$ with $a<b$, where
\begin{itemize}
    \item The interval is partitioned as $a=s_0<s_1< \ldots < s_{n-1}<s_m=b$ into a finite number of subintervals.
    \item The image of each subinterval $\gamma_0(s_{j-1},s_j)$ is a trajectory of either the fast subsystem or the slow subsystem.
    \item The image $\gamma_0(a,b)$ is oriented consistently with the orientations on each subinterval induced by the fast or slow flows.
\end{itemize}
\end{definition}

Note that $s$ is a parametrisation of the curve rather than fast or slow time. Consequently, the image of a subinterval can be complete homoclinic or heteroclinic orbits of the fast or slow subsystem. In typical cases, the attractor will consist of subintervals that alternate between fast and slow segments, but this may not be the case at bifurcation. In the case that all slow segments are on $\cC_{att}[g]$, this will typically perturb \cite{maesschalck11} to similar trajectories $\epsilon\rightarrow 0$. If the slow segments explore other hyperbolic points on the critical manifold, canard trajectories may appear. Several possible cases of fast/slow trajectories are considered in \cite{maesschalck11}.

\subsection{Global equivalence of critical sets}

In order to define persistence and bifurcation of critical sets, we need a notion of equivalence between critical sets. We define an equivalence of critical sets, called global equivalence, through the parametrized fast vector fields that generate them. We adopt the equivalence in \cite{golubitskyscheffer85} for the case $m=n=1$, and conjecture that that equivalence or a similar one may work also for the case $m=1,n=2$. A local version of this equivalence due to \cite{peters91} is presented at the end of this section. Classification of local bifurcations (classification of bifurcation of germs) has been worked out for the local equivalence \cite{peters91}, but it remains to be shown if the equivalence is also suitable for global classification of bifurcations.

We consider {\em fast} and {\em slow} vector fields $g: \bR^{m+n} \to \bR^m$ and $h: \bR^{m+n}\to \bR^n$, with the generic assumptions of smoothness and being in the singular limit $\epsilon \to 0$. More precisely we consider
\begin{equation}
V'_{f}:= C^{\infty}(\bR^{m+n},\bR^m),~~V'_{s}:=C^{\infty}(\bR^{m+n},\bR^n).
\label{eq:restrictedsystem}
\end{equation}
We furthermore consider only systems (\ref{eq:mainsystem}) with bounded global attractors.

We restrict ourselves to vector fields defined on some fixed compact regions $M\subset \bR^m$ and $N\subset \bR^n$ with open interior and smooth boundaries having outward normals $m(x)$ and $n(y)$ respectively. Implicitly fixing these $M$ and $N$, we define the open sets
\begin{equation}
V_{f}:= \left\{g\in V'_{f}~:~ \begin{array}{c}g(x,y) \cdot n(x)<0 ~\forall (x,y)\in(\partial M, N)~
\mbox{ and }\\ g(x,y)=0\Rightarrow g_x(x,y)\neq 0 ~ \forall (x,y)\in (M,\partial N)\end{array}\right \}
\label{eq:Xf}
\end{equation}
and
\begin{equation}
V_{s}:= \{h\in V'_{s}~:~h(x,y) \cdot m(y)<0 ~ \forall (x,y)\in (M,\partial N)\}.
\label{eq:Vs}
\end{equation}
Note that if $(g,h)\in M \times N$ then the forward dynamics must remain in $M\times N$ and that there are no tangencies of the flow, (or the critical set) with the boundary. These conditions ensure that the properties persist under small perturbations of the vector field.

Recall now that the critical set and the fast dynamics depend only on $g$, and suppose that $g,\tilde{g}\in V_{f}$. We say that the zero sets ($\cC[g]$ and $\cC[\tilde{g}]$) of $g$ and $\tilde{g}$ are {\em globally equivalent on $M\times N$}, denoted $g\sim \tilde{g}$, (c.f. \cite[p144]{golubitskyscheffer85}) if there are functions $Y(y):N \rightarrow \bR^n$, $X(x,y):M \times N \to \bR^m$ and $S:M \times N \rightarrow (0,\infty)$ such that:
\begin{equation}
\tilde{g}(x,y)=S(x,y)g(X(x,y),Y(y))
\label{eq:globalequivalence}
\end{equation}
i.e. we only consider changes in coordinate that map fast dynamics to fast dynamics up to a possible change in timescale. More precisely, we assume that:
\begin{itemize}
\item The map $\Phi(x,y):=(X(x,y),Y(y))$ is a diffeomorphism
\item The map $S(x,y)>0$ is smooth on $M\times N$
\end{itemize}
The requirement that $S(\cdot,\cdot)>0$ ensures that trajectories preserve their time orientation under equivalence. Note that we define the equivalence of critical sets through the functions that generate them. Consequently, equivalence of critical sets $\cC[g]$ and $\cC[\tilde{g}]$ does {\em not} imply equivalence of every level set of $g$ and $\tilde{g}$, such as $g=1$ and $\tilde{g}=1$.

For $m=1$ and $n=2$ we state a local equivalence adapted from \cite[Definition 2.1, p6]{peters91}:
\begin{equation}
\tilde{g}(x,y_1,y_2)=S(x,y_1,y_2)g(X(x,y_1,y_2),Y_1(y_1,y_2),Y_2(y_1,y_2)),
\label{eq:globalequivalencetwoslow}
\end{equation}
where we explicitly write $y=(y_1,y_2)$, where $Y_1(y_1,y_2),Y_2(y_1,y_2):N\rightarrow \bR^2$, $X(x,y_1,y_2): M \times N \to \bR$, and $S(x,y_1,y_2):M\times N \rightarrow (0,\infty)$ are smooth functions, and additionally
\begin{equation*}
\left|
\begin{array}{c c}
\frac{\partial Y_1(x,y_1,y_2)}{y_1} & \frac{\partial Y_1(x,y_1,y_2)}{y_2} \\
\frac{\partial Y_2(x,y_1,y_2)}{y_1} & \frac{\partial Y_2(x,y_1,y_2)}{y_2}
\end{array}
\right| \neq 0,
\end{equation*}
and $\frac{\partial X(x,y_1,y_2)}{\partial x}>0$ for every $(x,y_1,y_2) \in M\times N$. Some further conditions are imposed in \cite{peters91} since \cite{peters91} deals with germs. Note that smoothness combined with the conditions on $X$, $Y_1$ and $Y_2$ makes $(x,y_1,y_2) \rightarrow (X(x,y_1,y_2),Y_1(y_1,y_2),Y_2(y_1,y_2))$ a local diffeomorphism.

Since global equivalence should imply local equivalence we expect the classification of local bifurcations under a global equivalence, perhaps (\ref{eq:globalequivalence}), to coincide with the classification under a local equivalence such as (\ref{eq:globalequivalencetwoslow}). The latter was worked out in \cite{peters91}.

We leave the generalisation of the global equivalence (\ref{eq:globalequivalence}) to the case $m>1,n>1$ open.

\section{Persistence and bifurcation of critical sets}
\label{sec:criticalmanifolds}

Assume we define $V_f$ as in (\ref{eq:Xf}) for some compact regions $M$ and $N$. In order to define persistence of the critical sets we define the unfolding of the slow dynamics following \cite[Section III]{golubitskyscheffer85}. We say a smooth function $G(x,y,\lambda)$ for $\lambda \in \bR^r$ is an $r$-parameter \emph{unfolding} of $g(x,y)$ if
\begin{equation*}
G(x,y,0) = g(x,y)
\end{equation*}
for all $(x,y)\in M\times N$. Reference \cite{golubitskyscheffer85} mostly assumes $G$ and $g$ are germs of vector fields, though in \cite[Theorem III.6.1]{golubitskyscheffer85} the equivalence is global within a compact region, as we consider here.

If $G$ and $H$ are both unfoldings of $g$, we say that $H$ {\em factors through} $G$ if there exist smooth mappings $S,X,Y,L$ and $W \subset \bR^r$, a neighbourhood of $0$, such that
\begin{equation*}
H(x,y,\lambda) = S(x,y,\lambda)G(X(x,y,\lambda), Y(y,\lambda), L(\lambda)), 
\end{equation*}
for all $\lambda \in W$ and $(x,y) \in (M,N)$, where $S(x,y,0) = 1, X(x,y,0) = x, Y(y,0) = y, L(0) = 0$ (see \cite{golubitskyscheffer85}). 
We define $G$ to be a {\em versal} unfolding if every unfolding $H$ of $g$ factors through $G$. We say $g$ is {\em persistent} if it is its own unfolding, i.e. for any unfolding $G \in C^{\infty}(M \times N \times \bR^r)$  such that $G(x,y,0)=g(x,y)$, on $M \times N$  there is a neighbourhood $W$ of $0$ in $\bR^r$ such that
\begin{equation*}
G(x,y,\lambda) \sim g(x,y), ~\forall \lambda \in W,
\end{equation*}
where, as before, $\sim$ denotes global equivalence on $M\times N$.

If the unfolding is versal and contains a minimum number of parameters, we call it a \emph{universal unfolding} \cite{golubitskyscheffer85}. The number of parameters $\lambda$ in such a universal unfolding $G$ is the \emph{codimension} of $g$. In particular, if $g$ is \emph{persistent} then $g$ is its own universal unfolding, and in this case we say it has \emph{codimension zero}. We say that a \emph{bifurcation} of $g$ occurs if $g$ is non-persistent, in which case the \emph{codimension} of the bifurcation is that of the universal unfolding of $g$. We emphasise once more that the equivalence relation (\ref{eq:globalequivalence}) concerns the zero sets of $g$, i.e. the critical set $\cC[g]$. Hence, persistence and bifurcation of $g$ under this equivalence is identified with persistence and bifurcation of the critical set.

Note that the aforementioned meaning of \emph{persistence} does not concern persistence to perturbations involving the scale separation parameter $\epsilon$, which is sometimes the case in the fast-slow literature, e.g. \cite{fenichel79,maesschalck11}. Here, we study exclusively the system (\ref{eq:mainsystem}) in the singular limit $\epsilon\to 0$.

\subsection{Persistence and codimension one bifurcation of critical sets for one fast and one slow variable}

The case $m=n=1$ can be directly treated using the global bifurcation theory with distinguished parameter approach of \cite[Section III]{golubitskyscheffer85}. Here, a distinguished parameter is a parameter which is considered integral to the model and separate from unfolding parameters, which represent model perturbations. In the fast-slow setting we consider the slow variable $y$ to be a distinguished parameter from the point of view of the layer equations (\ref{eq:layer}). Consider some $g\in V'_{f}$ and note that the critical set is
\begin{equation*}
\cC[g]=\{p=(x,y)\in \bR\times \bR~:~g(p)=0\},
\end{equation*}
and that the fold set is
\begin{equation*}
\cF[g]=\{p\in \cC[g]~:~g_x(p)=0\}.
\end{equation*}

Table~\ref{tab:D11nonpersistent} lists the three \emph{degenerate fold} sets $\cD_{i}[g],i=\{1,2,3\}$ for $m=n=1$: fold tangency, hysteresis point, and multiple limit point. The term limit point is a historical term for fold point. The set of all degenerate folds is then
\begin{equation}
\cD[g]=\cD_1[g]\cup\cD_2[g]\cup\cD_3[g],
\label{eq:defDg11}
\end{equation}
and any point in $\cF[g]\setminus\cD[g]$ is a {\em non-degenerate fold point}. Note that \cite{golubitskyscheffer85} refers to the fold tangency as a ``simple bifurcation'' and a multiple limit point as a ``double limit point'' but our notation offers easier generalization to higher $n$. The following theorem characterizes the persistent critical sets, using a result from \cite{golubitskyscheffer85}.

\begin{table}
\caption{Degenerate fold sets for $m=n=1$: Proposition~\ref{prop:CM110} states that if $\cD[g]=\cD_{1}[g]\cup\cD_{2}[g]\cup\cD_{3}[g]$ defined in (\ref{eq:defDg11}) is empty then $g\in V_f$ is persistent on $M\times N$
}
\centering
\begin{tabular}{|p{0.3\linewidth}|l|}
\hline
Fold tangency: & $\cD_{1}[g]~=~\{p\in \cF[g]:g_y(p)=0\}$,\\
\hline
Hysteresis point: & $\cD_{2}[g]~=~\{p\in \cF[g]:g_{xx}(p)=0\}$,\\
\hline
Multiple~limit~point: & $\cD_{3}[g]~=~ \{p\in \cF[g]: ~|\Pi(p)|\geq 2\}$.\\
\hline
\end{tabular}
\label{tab:D11nonpersistent}
\end{table}

\begin{prop}[Codimension zero, $m=n=1$]
\label{prop:CM110}
In the case $m=n=1$, if $g\in V_{f}$ has no degenerate folds (i.e. if $\cD[g]=\emptyset$) then the critical set $\cC[g]$ is persistent to smooth perturbations.
\end{prop}

\proof
We apply \cite[Theorem 6.1]{golubitskyscheffer85}: this states that there is bifurcation equivalence if there are no (a) simple bifurcations (here called fold tangencies), (b) hysteresis points  (c) double limit points (here called multiple limit points) or (d) codimension one interactions of equilibria or folds with the boundaries. The assumptions in (\ref{eq:Xf}) are open conditions that ensure that (d) does not happen and that any unfolding of $g$ will remain within $V_f$ for small enough perturbations. Hence the only obstructions are (a-c) which are avoided if $\cD[g]$ is empty.
\qed

\begin{figure}
\begin{center}
\includegraphics[width=15cm]{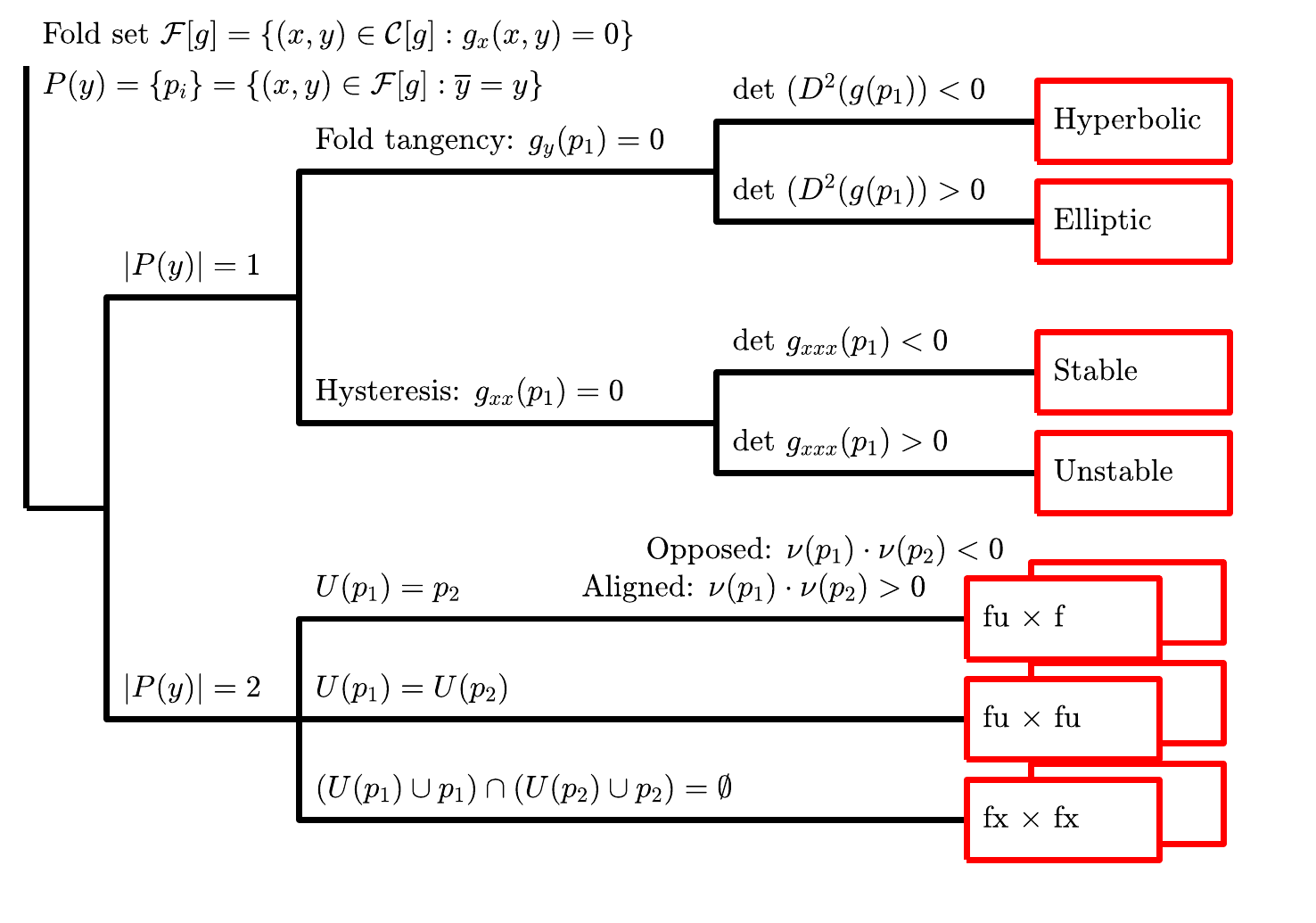}
\end{center}
\caption{
(Color online) Conditions that lead to codimension one degeneracies of the critical set for $m=n=1$ (see also Table~\ref{tab:D11persistentsubcases}). Note that $\nu(p)$ is the direction vector of a fold at a point $p$ and $\det(D^2(g(p)))$ is the Hessian of $g$ at $p$. Similarly, $f$ means fold, $fu$ means fold umbra and $fx$ means non-interacting fold. For ease of notation, we suppress explicit dependence on $g$, such that e.g. $\nu[g](p)=\nu(p)$. For a persistent codimension one bifurcation exactly one branch must be followed for exactly one fast fibre (a single $y$), leading to one of the red boxes. Overlapping red boxes symbolise that the degeneracy can be of either aligned or opposed type. See the text for details}
\label{fig:degeneracytable11}
\end{figure}

To aid the classification of codimension one bifurcations, we define $\cD^1_i[g]$ in Table~\ref{tab:D11persistent}. These are open dense subsets of $\cD_i[g]$ that avoids obvious further degeneracies. We then subdivide these cases further in  Table~\ref{tab:D11persistentsubcases}. A vector field $g$ is degenerate at codimension one if exactly one of these degeneracies $\cD_{i,j}^1[g,y]$ occur for exactly one slow coordinate $y$ (in exactly one fast fibre) which means we only need to compare points in $P(y)$ for some $y$. We avoid higher codimension fold tangency by precluding the cases $\det(D^2g) = g_{xx}g_{yy} - g_{xy}^2=0$ or higher order hysteresis $g_{xxx}=0$. Note that $\cD_{i,j}^1[g,y]$ depends explicitly on $y$. This choice makes it easier notation-wise to preclude the critical set from being degenerate at multiple $y$, which would raise codimension. The next result shows that Table~\ref{tab:D11persistentsubcases} and Figure~\ref{fig:degeneracytable11} give a complete list of codimension one bifurcations for this case.

\begin{table}
\caption{Subsets of $\cD_i[g]$ for $m=n=1$ whose union contains all codimension one bifurcations. Note that $\det(D^2 g)=g_{xx}g_{yy}-g_{xy}^2$ and that the first two are local degeneracies
}
\centering
\begin{tabular}{|p{0.3\linewidth}|l|}
\hline
Quadratic fold tangency: & $\cD_{1}^1[g]~=~\{p\in \cD_1[g]:|\Pi(p)|=1$ and $\det(D^2g(p))\neq 0\}$,\\
\hline
Cubic hysteresis point: & $\cD_{2}^1[g]~=~\{p\in \cD_2[g]:|\Pi(p)|=1$ and  $g_{xxx}(p)\neq 0$\},\\
\hline
Double~limit~point: & $\cD_{3}^1[g]~=~ \{p\in \cD_3[g]: |\Pi(p)|=2\}$.\\
\hline
\end{tabular}
\label{tab:D11persistent}
\end{table}

\begin{table}
\caption{Complete list of degeneracies that lead to codimension 1 bifurcations listed in Proposition~\ref{prop:CM111}. We write $P(y)=\{p_i\}$ as the set of distinct singular points $p_i=(x_i,y)$ of the vector field $g$ with slow coordinate $y$. Note that local degeneracies have $|P(y)|=1$}
\begin{tabular}{|p{0.2\linewidth}|p{0.55\linewidth}|p{0.14\linewidth}|}
\hline
Hyperbolic fold tangency & $\cD_{1,1}[g,y]=\{P(y) \subset D_1^1[g]~:~|P(y)|=1$ and $\det(D^2 g(p))<0\}$ & Fig.~\ref{fig:structureCM111D1D2} a,b,c)\\
\hline
Elliptic fold tangency & $\cD_{1,2}[g,y]= \{P(y)\subset D_1^1[g]~:~|P(y)|=1$ and $\det(D^2 g(p))>0\}$ & Fig.~\ref{fig:structureCM111D1D2} d,e,f)\\
\hline
Stable hysteresis: & $\cD_{2,1}[g,y]=\{P(y)\subset D_2^1[g]~:~|P(y)|=1$ and $g_{xxx}(p)>0\}$ & Fig.~\ref{fig:structureCM111D1D2} g,h,i)\\
\hline
Unstable hysteresis: & $\cD_{2,2}[g,y]=\{P(y)\subset D_2^1[g]~:~|P(y)|=1$ and $g_{xxx}(p)<0\}$ & Fig.~\ref{fig:structureCM111D1D2} j,k,l)\\
\hline
Aligned umbra-fold double limit: & $\cD_{3,1}[g,y]=\{P(y)\subset D_3^1[g]:|P(y)|=2$ and $U[g](p_1)=p_2$ and  $\nu[g](p_1)\cdot\nu[g](p_2) > 0 $, for some $p_1,p_2\in P(y) \}$ & Fig.~\ref{fig:structureCM111D3} a,b,c)\\
\hline
Opposed umbra-fold double limit: & $\cD_{3,2}[g,y]=\{P(y)\subset D_3^1[g]~:~|P(y)|=2$ and $U[g](p_1)=p_2$ and  $\nu[g](p_1)\cdot\nu[g](p_2) < 0$, for some $p_1,p_2\in P(y) \}$ & Fig.~\ref{fig:structureCM111D3} d,e,f)\\
\hline
Aligned umbra-umbra double limit: & $\cD_{3,3}[g,y]=\{P(y)\subset D_3^1[g]~:~|P(y)|=2$ and $U[g](p_1)=U(p_2)$ and $\nu[g](p_1)\cdot\nu[g](p_2) > 0$, for some $p_1,p_2\in P(y) \}$ & Fig.~\ref{fig:structureCM111D3} g,h,i)\\
\hline
Opposed umbra-umbra double limit: & $\cD_{3,4}[g,y]=\{P(y)\subset D_3^1[g]~:~|P(y)|=2$ and $U[g](p_1)=U(p_2)$ and  $\nu[g](p_1)\cdot\nu[g](p_2) < 0$, for some $p_1,p_2\in P(y) \}$ & Fig.~\ref{fig:structureCM111D3} j,k,l)\\
\hline
Aligned non-interacting double limit: & $\cD_{3,5}[g,y]=\{P(y)\subset D_3^1[g]~:~|P(y)|=2$ and $(U[g](p) \cup p) \cap_{p_i\neq p}(U[g](p_i) \cup p_i))=\emptyset$ and  $\nu[g](p_1)\cdot\nu[g](p_2) > 0$, for some $p_1,p_2\in P(y) \}$ & Fig.~\ref{fig:structureCM111D3} m,n,o)\\
\hline
Opposed non-interacting double limit: & $\cD_{3,6}[g,y]=\{P(y)\subset D_3^1[g]~:~|P(y)|=2$ and $ U[g](p) \cap (U(P(y)) \cup P(y) \setminus U[g](p)=\emptyset, ~\forall p\in P(y)$, and  $\nu[g](p_1)\cdot\nu[g](p_2) < 0$, for some $p_1,p_2\in P(y) \}$ & Fig.~\ref{fig:structureCM111D3} p,q,r)f\\
\hline
\end{tabular}
\label{tab:D11persistentsubcases}
\end{table}

\newpage

\begin{figure}
\begin{center}
\includegraphics[width=8cm]{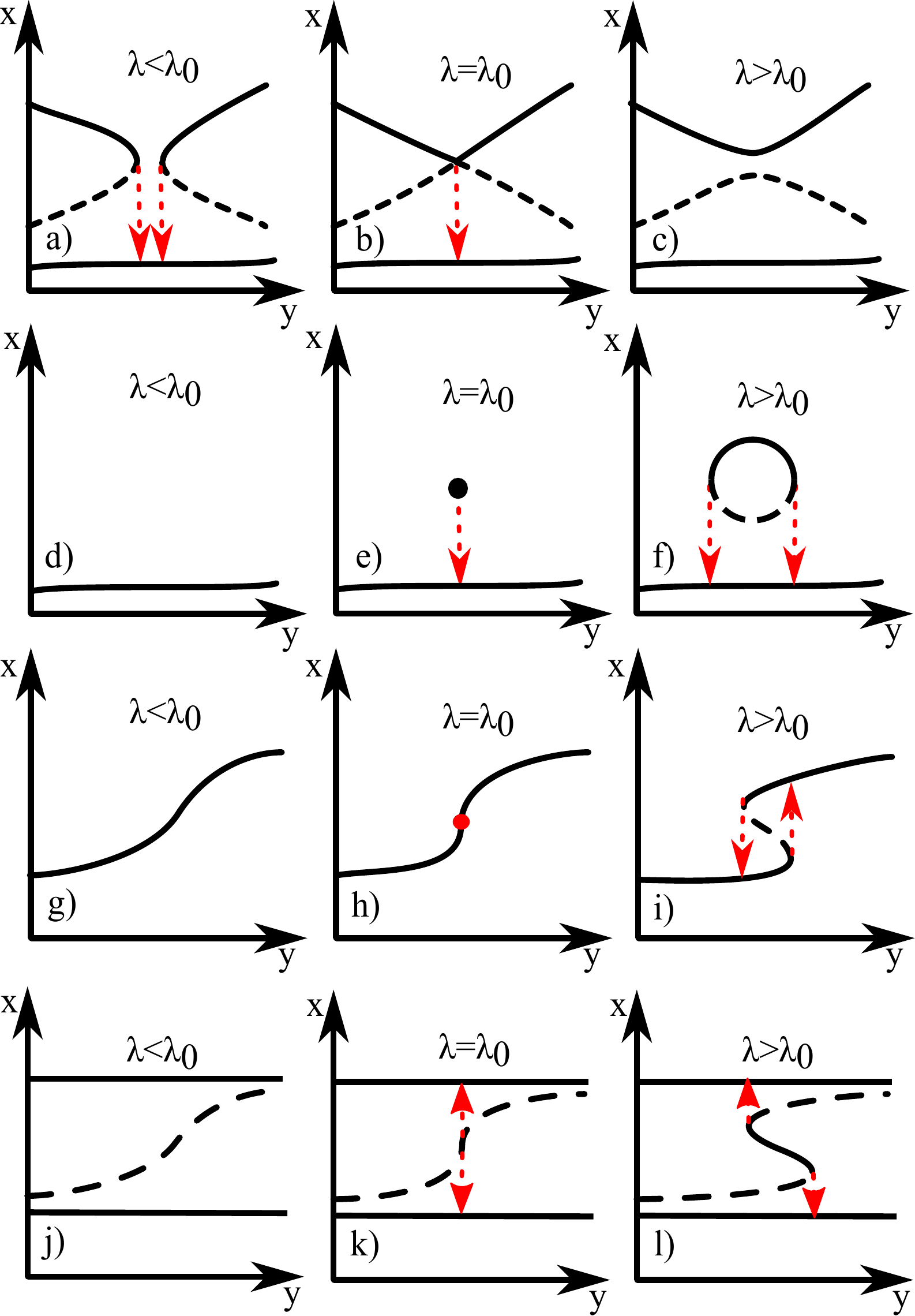}
\end{center}
\caption{
(Color online) Unfoldings of local codimension one bifurcations of the critical set for $m=n=1$. Solid black lines show $\cC_{att}[g]$, dashed black lines show $\cC_{rep}[g]$ while red arrows show the umbral map from fold points. Note that the fast variable $x$ is plotted on the vertical axis, and that the slow variable $y$ is plotted on the horizontal axis.}
\label{fig:structureCM111D1D2}
\end{figure}

\begin{figure}
\begin{center}
\includegraphics[width=15cm]{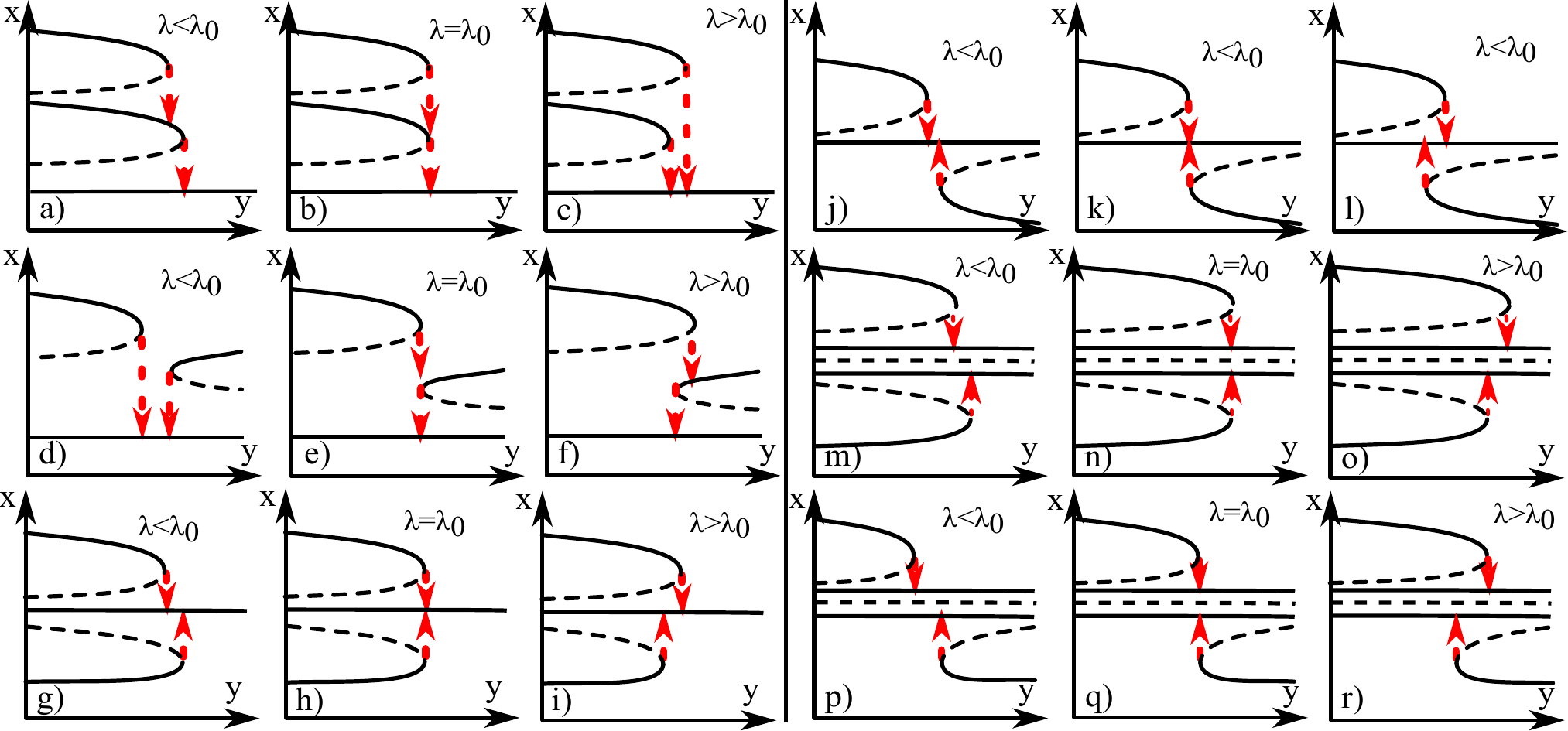}
\end{center}
\caption{
(Color online) Unfoldings of subcases of double limit point degeneracy, the global codimension one bifurcation of the critical set for $m=n=1$ (Table \ref{tab:D11persistentsubcases}). Bifurcation occurs when the unfolding (bifurcation) parameter $\lambda$ equals the critical value $\lambda_0$. Solid black lines show $\cC_{att}[g]$, dashed black lines show $\cC_{rep}[g]$, and red arrows show the umbral map from fold points. As in Figure \ref{fig:structureCM111D1D2}, note that the fast variable $x$ is plotted on the vertical axis, and that the slow variable $y$ is plotted on the horizontal axis}
\label{fig:structureCM111D3}
\end{figure}

\begin{prop}[Codimension one, $m=n=1$]
\label{prop:CM111}
For $n=1$ and $m=1$ the codimension one bifurcations of critical sets $\cC[g]$ for $g\in V_{f}$ are characterised in Fig.~\ref{fig:degeneracytable11}, such that one of the sets $\cD_{j,k}[g,y]$ in Table~\ref{tab:D11persistentsubcases} is non-empty for precisely one $y$. At such a bifurcation, precisely one of the following occurs:
\begin{enumerate}
\item Two folds merge at a fold tangency (e.g. Fig.~\ref{fig:structureCM111D1D2} a,b,c) or d,e,f)).
\item Two folds merge at a hysteresis bifurcation (e.g. Fig.~\ref{fig:structureCM111D1D2} g,h,i) or j,k,l)).
\item Two fold points share the same slow coordinate: there are six distinct ways this can occur (e.g. Fig.~\ref{fig:structureCM111D3})
\end{enumerate}
\end{prop}

\proof
To avoid persistence, at least one of the the degeneracies $\cD_{i}[g]$ listed in Table~\ref{tab:D11nonpersistent} must occur for some $i\in\{1,2,3\}$: as these are independently defined we can assume that only one will occur for an open dense set of unfoldings. Without loss of generality we can assume that the open conditions in Table~\ref{tab:D11persistent} apply.
\qed

The subcases of $\cD^1_1[g]$ follow from examining the sign of $\det(D^2g)$: the hyperbolic fold tangency $\cD_{1,1}[g,y]$ is the simple bifurcation of \cite{golubitskyscheffer85} while the elliptic fold tangency $\cD_{1,2}[g,y]$ is also called the isola. Similarly, the cubic hysteresis $\cD^1_2[g]$ can be either stable or unstable, depending on the sign of the leading order term. These cases can be transformed into the normal forms of Table~\ref{tab:bifnormalforms} (these are given in \cite{golubitskyscheffer85}). The cases $\cD^1_1[g]$ unfold on varying a typical parameter $\lambda$ as shown in Fig.~\ref{fig:structureCM111D1D2} a,b,c) and d,e,f) respectively, while $\cD^1_2[g]$ unfold as shown in Fig.~\ref{fig:structureCM111D1D2} g,h,i) and j,k,l).

The double limit point degeneracy $\cD^1_3[g]$ can be split into several subsets according to the direction of the folds given by the signs of
\begin{equation*}
\nu[g](p)=g_{xx}(p)g_y(p)
\end{equation*}
at the two limit points, and $k$, the number of regular sheets that separate them. The number $k$ determines whether the umbrae and folds intersect. If $k=0$, then the umbra of one fold intersects the other fold, if $k=1$ then the umbrae of the folds intersect, and if $k\geq2$ then the umbrae and folds do not intersect. The six distinct subcases of $\cD^1_3[g]$ are shown in Figure~\ref{fig:structureCM111D3}. 

We conclude this section with a few comments. First, it may appear as if the non-interacting double limit point degeneracy is not a degeneracy. However, the critical sets in e.g. Figure~\ref{fig:structureCM111D3} m,o) cannot be equivalent to that in Figure~\ref{fig:structureCM111D3} n) since the equivalence (\ref{eq:globalequivalence}) preserves the number of zeros in fast fibres: at bifurcation there is a $y$ for which $\cC[g]$ has five zeros, while the perturbed diagrams have at most four zeros for any given $y$. However, non-interacting double limit point degeneracy does not cause bifurcation of singular relaxation oscillations, as is shown in Section \ref{sec:persistbif11}.

\begin{figure}
\begin{center}
\includegraphics[width=10cm]{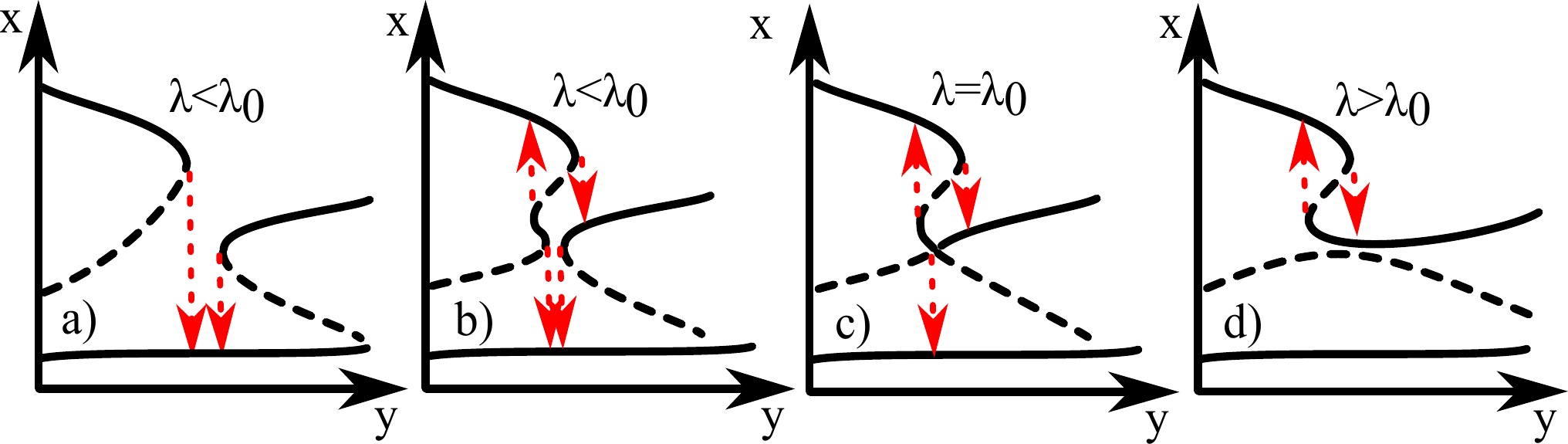}
\end{center}
\caption{
(Color online) Illustration that separate branches of the critical set always meet at a hyperbolic fold tangency with vertical tangent at codimension one. As the bifurcation parameter $\lambda$ increases, it brings the branches of the critical set closer. First the upper branch undergoes hysteresis bifurcation which produces a new fold which eventually merges with the lower branch in a hyperbolic fold tangency}
\label{fig:foldtangencyclarification}
\end{figure}

It may seem from Figure \ref{fig:structureCM111D1D2} that the hyperbolic fold tangency is of codimension higher than one. That this is not the case is shown algebraically in \cite{golubitskyscheffer85} using a considerable technical machinery. However, in Figure \ref{fig:foldtangencyclarification} we aim to provide some intuition that the bifurcation does not require fine tuned perturbations, but rather arises generically as two branches of the critical set approach under one-parameter perturbation.

It might seem puzzling why some combinations of $x$ and $y$ are left out from the normal forms in Table \ref{tab:bifnormalforms}. In general, which terms can be included in a certain normal form is a non-trivial question, treated in e.g. \cite{golubitskyscheffer85}. However, for the particular case of fold tangency the transformation $X(x,y) = x-ay$, $Y(y)=y/\sqrt{|1-a^2|}$ turns the polynomial $g(x,y)=x^2+y^2+2axy +\lambda $ into one of the normal forms in Table \ref{tab:bifnormalforms} as long as $g$ is non-degenerate. This transformation clearly preserves equivalence class under (\ref{eq:globalequivalence}).

In the above example as well as in general, the diffeomorphism $\Phi(x,y) = (X(x,y),Y(y))$ preserves the fast-slow structure of (\ref{eq:mainsystem}). This is seen by letting $(\hat{x},\hat{y}) = \Phi(x,y)$, such that $(x,y) = \Phi^{-1}(\hat{x},\hat{y}):=(\hat{X}(\hat{x},\hat{y}),\hat{Y}(\hat{y}))$, and changing variables in (\ref{eq:mainsystem})
\begin{equation*}
\left\{\begin{array}{rl}
\epsilon \left(\frac{\partial \hat{X}}{\partial \hat{x}}\dot{\hat{x}} + \frac{\partial \hat{X}}{\partial \hat{y}}\dot{\hat{y}} \right) & = g(\hat{x},\hat{y}) \\
\frac{\partial \hat{Y}}{\partial \hat{y}}\dot{\hat{y}} & = h(\hat{x},\hat{y}).
\end{array}\right.
\end{equation*}
After rearranging the above equation, we get the system
\begin{equation*}
\left\{\begin{array}{rl}
\epsilon \dot{\hat{x}} & = \left(g(\hat{x},\hat{y}) - \epsilon \frac{\partial \hat{X}}{\partial \hat{y}} h(\hat{x},\hat{y})/\frac{\partial \hat{Y}}{\partial \hat{y}}\right)/\frac{\partial \hat{X}}{\partial \hat{x}} \\
\dot{\hat{y}} & = h(\hat{x},\hat{y})/\frac{\partial \hat{Y}}{\partial \hat{y}},
\end{array}\right.
\end{equation*}
which is on the form (\ref{eq:mainsystem}). The above expression is well defined since $\Phi(x,y)$ is assumed to be a diffeomorphism.

\begin{table}
\caption{Normal forms (for $m=n=1$) and hypothesised normal forms (for $m=1,n=2$) of local codimension one bifurcations of the critical manifold. Different signs of $\delta_1,\delta_2\neq0$ give different subcases of bifurcation
}
\centering
\begin{tabular}{|p{0.25\linewidth}|p{0.65\linewidth}|}
\hline
$m=n=1$ & \\
\hhline{|=|=|}
Fold tangency: & $g(x,y)=x^2 + \delta_1y^2 + \lambda$,\\
\hline
Hysteresis: & $g(x,y)=\delta_1x^3 + \lambda x + y$,\\
\hline
$m=1$,~$n=2$ & \\
\hhline{|=|=|}
Fold tangency: & $g(x,y_1,y_2)=x^2 + \delta_1y_1^2 + \delta_2y_2^2 + \lambda$,\\
\hline
Cusp tangency: & $g(x,y_1,y_2)=\delta_1x^3 + \delta_2 xy_2^2 +\lambda x + y_1$,\\
\hline
Swallowtail: & $g(x,y_1,y_2)=x^4 + \lambda x^2 + y_1x + y_2$,\\
\hline
\end{tabular}
\label{tab:bifnormalforms}
\end{table}
\FloatBarrier

\subsection{Persistence of critical sets for one fast and two slow variables}

In analogy with the $m=n=1$ case we give a conjectured list of all degeneracies that can cause nonpersistency of vector fields with one fast and two slow variables, up to codimension one. First, we introduce some notation. For any $g\in V'_{f}$ we write 
\begin{equation*}
g_x = \frac{\partial g}{\partial x},~\nabla_y g = (g_{y_1},g_{y_2})~\mbox{ and }\nabla^{\perp}_y g = (-g_{y_2},g_{y_1}).
\end{equation*}
By $u||v$ we mean that vectors $u$ and $v$ are parallel. The non-zero vector $u$ rescaled to unit length is denoted $\overline{u} = u/|u|$. $D^2(g)$ is the Hessian of $g$ with respect to all components of $p=(x,y_1,y_2)$: 
\begin{equation*}
[D^2(g)]_{ij} = \frac{\partial^2g}{\partial p_i p_j},~i,j\in \{1,2,3\}.
\end{equation*}
The slow Hessian $D_y^2(g)$ is defined analogously, but with $i,j\in \{2,3\}$.

Recall that the critical set is
\begin{equation*}
\cC[g]=\{p=(x,y)\in \bR\times \bR^2~:~g(p)=0\}
\end{equation*}
and the fold set is
\begin{equation*}
\cF[g]=\{p\in \cC[g]~:~g_x(p)=0\}.
\end{equation*}
As folds are not typically isolated in this case, we also need to distinguish between quadratic folds, cubic cusps and higher order cusps (Table~\ref{tab:genericsingularities12}).

\begin{table}
\caption{Singularities of the critical set for one fast and two slow variables}
\begin{equation*}
\begin{array}{|l|l|}
\hline
\mbox{Quadratic fold:} &\cF_0[g]=\{p\in \cF[g]~:~g_{xx}(p)\neq 0\}\\
\hline
\mbox{Cubic cusp:} & \cF_1[g]= \{ p\in \cF[g]~:~g_{xx}(p)= 0,~g_{xxx}(p)\neq 0 \}\\
\hline
\mbox{Higher order cusp:} & \cF_2[g]= \{ p\in \cF[g]~:~g_{xx}(p)= 0,~g_{xxx}(p)= 0 \}\\
\hline
\end{array}
\end{equation*}
\label{tab:genericsingularities12}
\end{table}

We believe that the list of degenerate sets $\cD_i[g]$ of $\cF[g]$ given in Table~\ref{tab:nonpersistentdegeneracies12} is an exhaustive list of degeneracies under a suitable equivalence. These degeneracies are natural extensions from the degeneracies for $m=n=1$; generic objects (quadratic fold lines and cubic cusps) can intersect ($\cD_1[g]$, $\cD_2[g]$, $\cD_3[g]$), their projections onto the slow variables can become tangent $\cD_4[g]$ or intersect $\cD_5[g]$ and $\cD_6[g]$. More precisely, we define the set of degenerate points
\begin{equation*}
\cD[g]=\cD_1[g]\cup\cD_2[g]\cup\cD_3[g]\cup\cD_4[g]\cup\cD_5[g] \cup \cD_6[g].
\end{equation*}
Note that the umbral map is single valued for $p\in \cF[g]\setminus (\cD_2[g] \cup \cD_3[g])$.
If $p\in\cD_2 \cup \cD_3[g]$ then it can be zero, one or two-valued (see Figure~\ref{fig:structureCM121D2}). We now conjecture a persistence criterion for $m=1$, $n=2$ that is analogous to the $m=n=1$ case in Proposition~\ref{prop:CM110}.

\begin{conj}[Codimension zero, $m=1$, $n=2$.]
\label{conj:CM120}
For one fast and two slow variables, the critical set $\cC[g]$ is persistent to perturbations for $g\in V_{f}$ if all folds are non-degenerate, i.e. if $\cD[g]=\emptyset$.
\end{conj}

\begin{table}
\caption{Possible degeneracies of the critical set for one fast and two slow variables, $m=1$ and $n=2$. As before, the co-fold set $\Pi(p)$ is the subset in $\cF[g]$ sharing slow coordinate with the point $p$}
\begin{tabular}{|p{0.25\linewidth}|p{0.65\linewidth}|}
\hline
Fold tangency & $\cD_{1}[g] = \{p \in \cF[g]:~\nabla_y g(p)=0\}$\\
\hline
Cusp tangency & $\cD_{2}[g] = \{p \in \cF[g] \setminus \cF_0[g]:~\nabla^\perp_y g(p) \cdot \nabla_y g_x(p)=0\}$\\
\hline
High order cusp & $\cD_{3}[g]=\{p \in \cF[g]: ~g_{xxx}(p)=0\}$\\
\hline
Quadratic-fold projection tangency & 
$\cD_{4}[g]=\{p_1 \in \cF_0[g]:$ $\nabla_yg(p_1)||\nabla_yg(p_2)$ 
for some $p_2 \in (\Pi(p_1) \setminus p_1) \cap \cF_0[g]\}$\\
\hline
Cusp projection intersection & 
$\cD_{5}[g]=\{p\in \cF[g]: |\Pi(p)|\geq2$ and $\Pi(p) 
\cap \cF_1[g] \neq \emptyset \}$ \\
\hline
Triple fold projection intersection & $\cD_{6}[g]=\{p\in \cF[g]:|\Pi(p)|\geq 3 \}$ \\
\hline
\end{tabular}
\label{tab:nonpersistentdegeneracies12}
\end{table}

Unfortunately, the method of proof in \cite[Thm III.6.1]{golubitskyscheffer85} used for Proposition~\ref{prop:CM110} does not easily generalize to this case of ``multiple distinguished parameters'' (two slow variables in our case). This is because degenerate cases appear with codimension infinity \cite{montaldi1994path}, at least if we consider the restricted global equivalence where we require that $q(y)$ is the identity in (\ref{eq:globalequivalence}). Hence proof of this result will require a less stringent (but still natural) form of global equivalence.

\subsection{Codimension one bifurcations of critical sets for one fast and two slow variables}
\label{sec:codimonecritsed12}

A variety of degeneracies can persistently occur for one parameter families, i.e. at codimension one. Table~\ref{tab:persistentdegeneracies12} lists degeneracies that we believe contain all persistent codimension one bifurcations for a suitable notion of global equivalence. We divide these into local or global degeneracies. The \emph{local} degeneracies (Table~\ref{tab:localdegeneracies12}) are denoted $\cD_{i,j}^1[g,y]$ for $i\in \{1,2,3\}$, the others involve interaction of two or more points in the same fast fibre of $y$: these \emph{global} degeneracies are denoted $\cD_{i,j}^1[g,y]$ for $i\in \{4,5,6\}$ and are listed in detail in \ref{app:global12tables}.

We note that our conjectured local codimension one degeneracies coincide with those in the classification of \cite{peters91}, which is found for the local equivalence (\ref{eq:globalequivalencetwoslow}); it might be possible to extend these results to a related global equivalence.

We first state necessary conditions for the degeneracies $\cD_i[g]$ in Table \ref{tab:nonpersistentdegeneracies12} to be codimension one. In this process we introduce some geometric notions, listed in Table \ref{tab:usefuldefinitions}.

\begin{table}
\caption{Quantities used in the classification of degeneracies for one fast and two slow variables, $m=1$ and $n=2$. The last column refers to appendices where all definitions except $\nu[g](p)$ are motivated (we do not motivate the natural definition of $\nu[g](p)$). See Figure \ref{fig:geometricdefinitions} for illustrations of all quantities except $W[g](p)$}
{
\begin{tabular}{|p{0.20\linewidth}|p{0.10\linewidth}|p{0.42\linewidth}|p{0.15\linewidth}|}
\hline
Quadratic fold direction vector & $\nu[g](p)$ & $=g_{xx}(p)\nabla_y g(p)$ & \\
\hline
Scalar quadratic fold curvature & $K[g](p)$ & $=\sign{(g_{xx}(p))}\frac{\overline{{\nabla_y^\perp g(p)}}^T D^2_y(g(p)){\overline{\nabla_y^\perp g(p)}}}{2|\nabla_y^\perp g(p)|}  - \frac{\left(\nabla_y g_x(p) \cdot \overline{\nabla_y^\perp g(p)}\right)^2}{8|g_{xx}(p)||\nabla_y^\perp g(p)|}$ & \ref{app:curvature} \\
\hline
Quadratic fold curvature vector & $\kappa[g](p)$  & $=\Bigg( \frac{\overline{{\nabla_y^\perp g(p)}}^T D^2_y(g(p)){\overline{\nabla_y^\perp g(p)}}}{2|\nabla_y^\perp g(p)|} - \frac{\left(\nabla_y g_x(p) \cdot \overline{\nabla_y^\perp g(p)}\right)^2}{8 g_{xx}(p) |\nabla_y^\perp g(p)|}\Bigg) \overline{\nabla_y g(p)}$ & \ref{app:curvature} \\
\hline
Cubic cusp direction vector & $\mu[g](p)$ & $=\frac{g_{xxx}(p)}{\nabla_y g_x(p) \cdot \overline{\nabla_y^\perp g(p)}}\nabla_y^\perp g(p)$ & \ref{app:cuspdirection} \\
\hline
Cusp quantity & $W[g](p)$ & $ = \frac{g_{xxx}(p)}{2}\overline{\nabla_y^\perp g^T(p)}D^2_y(g_x(p))\overline{\nabla_y^\perp g(p)} - \nabla_y g_{xx}(p) \cdot \overline{\nabla_y^\perp g(p)}$ & \ref{app:W} \\
\hline
\end{tabular}
}
\label{tab:usefuldefinitions}
\end{table}

\begin{figure}
\begin{center}
\includegraphics[width=8cm]{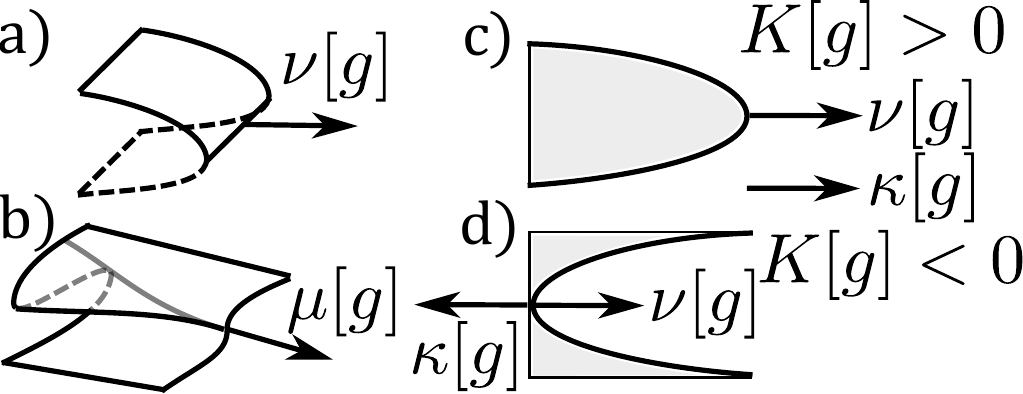}
\end{center}
\caption{
Illustration of fold and cusp direction vectors and fold curvature in the slow plane. a) Quadratic fold direction vector $\nu[g]$, b) cubic cusp direction vector $\mu[g]$, c) convex fold (viewed as a projection onto the slow subsystem) scalar quadratic fold curvature $K[g]$ and quadratic fold curvature vector $\kappa[g]$, d) same as c) but concave. See Table~\ref{tab:usefuldefinitions} for definitions of $\nu[g]$, $\mu[g]$, $\kappa[g]$ and $K[g]$
\label{fig:geometricdefinitions}
}
\end{figure}

Starting with quadratic fold degeneracies, we note that for typical fold tangency we require that the Hessian $D^2(g)$ has no zero eigenvalues. For typical fold projection tangency, we require that the two folds sharing slow coordinate $(y_1,y_2)$ are not aligned and with the same quadratic curvature. This is guaranteed if the quadratic fold curvature vectors $\kappa[g](p)$ (Table \ref{tab:usefuldefinitions} and Figure~\ref{fig:geometricdefinitions} c,d))
at fold points $p_1$ and $p_2\in \Pi(p_1)$ are distinct:
\begin{equation*}
\kappa[g](p_1) \neq \kappa[g](p_2).
\end{equation*}


For typical triple quadratic folds we require that exactly three fold points occur for every slow coordinate $(y_1,y_2)$ and that all folds are quadratic.

Turning to cusps, there is degeneracy if the cubic cusp direction vector $\mu[g](p)$ at a cubic cusp point $p$ (Table~\ref{tab:usefuldefinitions} and Figure~\ref{fig:geometricdefinitions} b)) is either zero or undefined. The case $\mu[g](p) = 0$ corresponds to $g_{xxx}(p)=0$, resulting in a swallowtail bifurcation. The case that $\mu[g](p)$ is undefined occurs if $\nabla^\perp_y g(p) \cdot \nabla_y g_x(p)=0$; this results in a cusp tangency. However, for \emph{typical} cusp tangency we require that the quadratic quantity $W[g](p)$ (Table~\ref{tab:usefuldefinitions}) is non-zero. Furthermore, the sign of $W[g](p)$ separates subcases of typical cusp tangency.

Finally, for typical cusp-fold projection intersection, we require that exactly one cubic cusp and one quadratic fold share slow coordinate $(y_1,y_2)$.

The quantities $\nu[g](p)$, $\mu[g](p)$, $W[g](p)$ and $\kappa[g](p)$ are discussed more in \ref{app:curvature}, \ref{app:cuspdirection} and \ref{app:W}. Hypothesised normal forms for the local codimension one bifurcations are listed in Table \ref{tab:bifnormalforms}.

\begin{table}
\caption{Subsets of $\cD[g]$, for one fast and two slow variables, whose union we conjecture contains all codimension one bifurcation sets. Note that $\cD_i^1[g]$ are local for $i=1,2,3$ and global for $i=4,5,6$. $\Pi(p)$ is the set of all singular points sharing slow coordinate with $p$. The scalars $a_1$, $a_2$, $a_3$ are coefficients in equation (\ref{eq:linearcombination}) and $\kappa[g](p)$, $\nu[g](p)$ and $\mu[g](p)$ are curvature and direction vectors at a point $p$. See the text for details}
\begin{tabular}{|p{0.25\linewidth}|p{0.65\linewidth}|}
\hline
Typical fold tangency & $\cD_{1}^1[g] = \{p \in \cD_1[g]:|\Pi(p)|=1$ and $\det(D^2 g(p)) \neq 0 \}$\\
\hline
Typical cusp tangency & $\cD_{2}^1[g] = \{p \in \cD_2[g]:|\Pi(p)|=1$ and $W[g](p) \neq 0\}$\\
\hline
Swallowtail & $\cD_{3}^1[g]=\{p \in \cD_3[g]: |\Pi(p)|=1$ and $g_{xxxx} \neq 0\}$\\
\hline
Typical double quadratic-fold projection tangency & $\cD_{4}^1[g]=\{p\in \cD_4[g]: |\Pi(p)|=|\Pi(p) \cap \cF_0[g]| = 2$ and $\kappa[g](p) \neq \kappa[g](q), ~\forall q \in \Pi(p) \setminus p \}$\\
\hline
Typical cubic-cusp - quadratic-fold projection intersection & $\cD_{5}^1[g]=\{p\subset \cD_5[g]:|\Pi(p)|=2$ and $|\Pi(p) \cap \cF_0[g]| = 1$ and $|\Pi(p) \cap \cF_1[g]| = 1$ and $\nu[g](p_1) \cdot \mu[g](p_2) \neq 0$ for some $p_1 \in \Pi(p) \cap \cF_0[g]$ and $p_2\in \Pi(q) \cap \cF_1[g] \}$ \\
\hline
Typical triple quadratic-fold projection intersection &  $\cD_{6}^1[g]=\{p\in \cD_6[g]:|\Pi(p)|=3$ and $|\Pi(p) \cap \cF_0[g]|=3$ and $a_1 \cdot a_2 \cdot a_3 \neq 0\}$.\\
\hline
\end{tabular}
\label{tab:persistentdegeneracies12}
\end{table}

\begin{table}
\caption{Subsets of \emph{local} degeneracies, for one fast and two slow variables parametrized by the slow coordinate, conjectured to include all codimension one bifurcations. Note that $P(y)$ is the set of all singular points of the vector field $g$ with slow coordinate $y$. The number of positive eigenvalues of the Hessian $\sign{(g_{xx})}D^2[g]$ at $p$ is written as $|\Sigma_+|$. Note that no eigenvalues are zero, since $\det(D^2 g(p)) \neq 0$ by assumption. Figure D1 is in \ref{app:global12tables}. See the text for details}
\begin{tabular}{|p{0.22\linewidth}|p{0.54\linewidth}|p{0.15\linewidth}|}
\hline
Wormhole fold tangency & $\cD_{1,1}[g,y] = \{P(y) \subset \cD_1^1[g]: |\Sigma_+|=1\}$ & Fig.~\ref{fig:structureCM121DSample} a,b,c)\\
\hline
Tube fold tangency & $\cD_{1,2}[g,y] = \{P(y)\subset \cD_1^1[g]:|\Sigma_+|=2\}$ & Fig.~\ref{fig:structureCM121DSample} d,e,f)\\
\hline
Isola fold tangency & $\cD_{1,3}[g,y] = \{P(y)\subset \cD_1^1[g]: |\Sigma_+|=3\}$ & Fig.~\ref{fig:structureCM121DSample} g,h,i)\\
\hline
Stable lips cusp tangency &$\cD_{2,1}[g,y] = \{P(y)\subset \cD_2^1[g]: W[g](p) > 0 \mbox{ and } g_{xxx}<0\}$ & Fig.~\ref{fig:structureCM121DSample} m,n,o)\\
\hline
Unstable lips cusp tangency &$\cD_{2,2}[g,y] = \{P(y)\subset \cD_2^1[g]:W[g](p) > 0 \mbox{ and } g_{xxx}>0\}$ & Fig.~\ref{fig:structureCM121D2} j,k,l) \\
\hline
Stable beaks cusp tangency &$\cD_{2,3}[g,y] = \{P(y)\subset \cD_2^1[g]:W[g](p) < 0 \mbox{ and } g_{xxx}<0\}$ & Fig.~\ref{fig:structureCM121DSample} j,k,l)\\
\hline
Unstable beaks cusp tangency &$\cD_{2,4}[g,y] = \{P(y)\subset \cD_2^1[g]:W[g](p) < 0 \mbox{ and } g_{xxx}>0\}$ & Fig.~\ref{fig:structureCM121D2} d,e,f)  \\
\hline
\mbox{Swallowtail:} & $\cD_{3,1}[g,y] = \{P(y) \subset \cD_3^1[g]\}$ & Fig.~\ref{fig:structureCM121DSample} p,q,r)\\
\hline
\end{tabular}
\label{tab:localdegeneracies12}
\end{table}
We now go trough the persistent subcases of codimension one degeneracies listed in Table~\ref{tab:persistentdegeneracies12}.

\subsubsection{Fold tangency}
\begin{figure}
\begin{center}
\includegraphics[width=15cm]{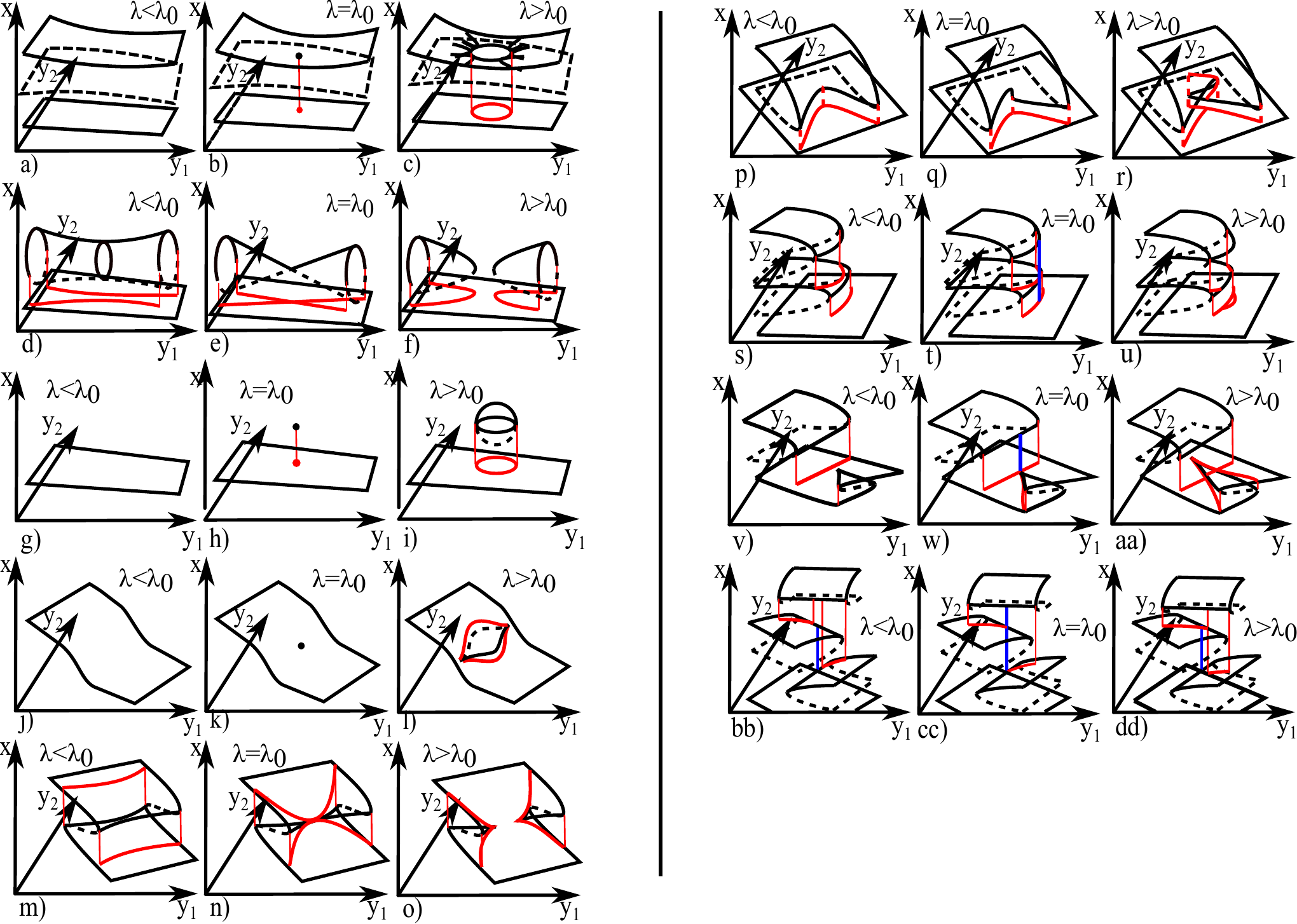}
\end{center}
\caption{
(Color online) Unfolding of examples of codimension one bifurcation of the critical set for $m=1$, $n=2$ (see Table \ref{tab:localdegeneracies12} for a) to r), and Table \ref{tab:globaldegeneracies12P2Folds}, Table \ref{tab:globaldegeneracies12P2Cusps} and Table \ref{tab:globaldegeneracies12P3} in the Appendix for s) to dd) ). Bifurcation occurs when the bifurcation parameter $\lambda$ equals the critical value $\lambda=\lambda_0$. Solid/dashed black lines show the stable/unstable sheets of the critical set while red lines show the image of the fold under the umbral map. Blue lines indicate special points of intersection}
\label{fig:structureCM121DSample}
\end{figure}

Fold tangency occurs when a pair or continuum of folds intersect. Typical fold tangency $\cD^1_1[g]$ is classified by $|\Sigma_+|$, the number of positive eigenvalues in $\Sigma$, the spectrum of $\sign{(g_{xx}(p))}D^2(g)(p)$. The cases of {\em wormhole} ($|\Sigma_+|=1$), {\em tube} $(|\Sigma_+|=2)$ and {\em isola} $(|\Sigma_+|=3)$ are shown in Figure~\ref{fig:structureCM121DSample} a) to i). Note that $g_{xx}(p)\det(D^2(g)(p)\neq 0$ implies that all non-positive eigenvalues are negative and that at least one eigenvalue is positive. 

\subsubsection{Cusp tangency}

At a cusp tangency, two cusps meet locally along a line. The subsets \emph{beaks} and \emph{lips} are distinguished by whether cusps are directed away from or towards each other before bifurcation (see Figure~\ref{fig:structureCM121DSample} and Figure~\ref{fig:structureCM121D2}). These degeneracies are named after the appearance of their projections onto the slow variables (Figure~\ref{fig:projectionsketch}), and the type is determined by the cusp quantity $W[g](p)$. The case $W[g]>0$ gives ``lips'' and $W[g]<0$ gives ``beaks''. We further subdivide these cases depending on their stability, determined by the sign of $g_{xxx}$.

\subsubsection{Swallowtail}

The swallowtail (Figure~\ref{fig:structureCM121DSample} p,q,r)) is well known from catastrophe theory and occurs when a fold ``folds over itself'' to create a degenerate fold that splits up into a pair of cusps.

\subsubsection{Fold projection tangency}

\begin{figure}
\begin{center}
\includegraphics[width=15cm]{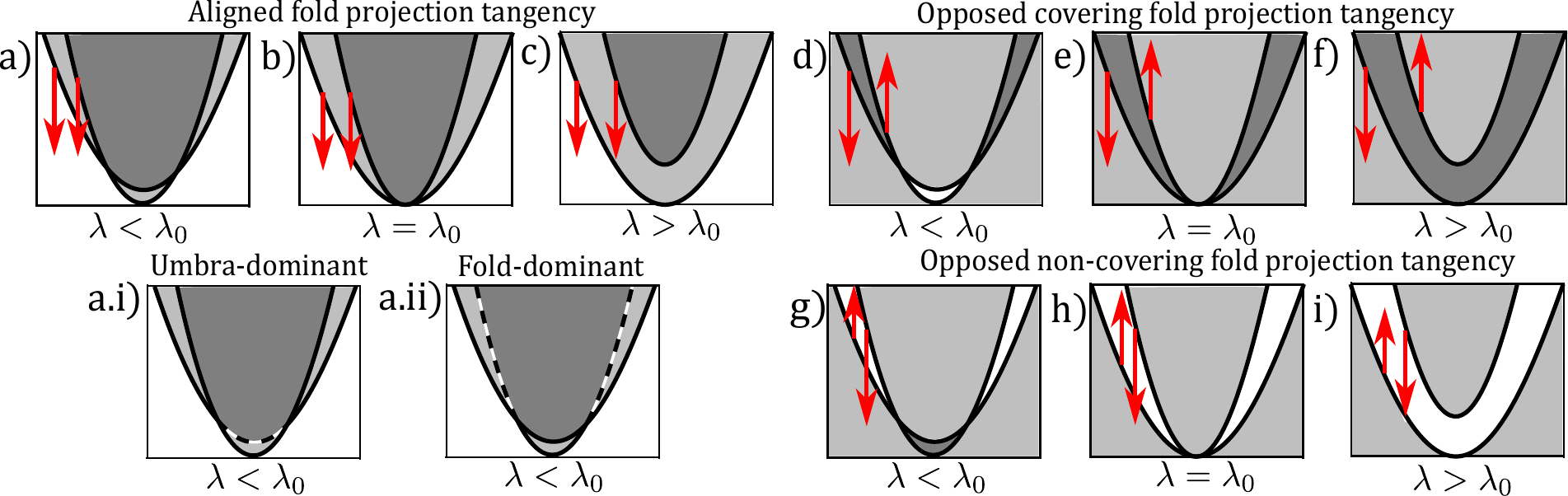}
\end{center}
\caption{
(Color online) Unfolding of a fold projection tangency, viewed in projection onto the slow plane. Three principal cases are shown in a,b,c), d,e,f) and g,h,i). Darker colours mean that more sheets of the critical manifold overlap. Red arrows show quadratic fold direction vectors at tangency points. The aligned fold-fold umbra subcase has umbra-dominant and fold-dominant subcases a.i) and a.ii) respectively. Dotted lines show parts of the destination fold covered by the umbral fold, seen from the stable side of the umbral sheet
\label{fig:foldtangencysketch}
}
\end{figure}

At a fold tangency, the projection of two curves of folds onto the slow variables are tangent. We divide fold projection tangency $\cD^1_{4}[g]$ into subcases depending on whether the folds at points $p_1$ and $p_2$ approach each other from the same direction (aligned) or opposite directions (opposed), captured by the sign of the inner product of the quadratic fold direction vectors $\nu[g](p_1)\cdot \nu[g](p_2)$. We further subdivide the opposed cases depending on the sign of the sum
\begin{equation*}
K[g](p_1) + K[g](p_2),
\end{equation*}
where $K[g](p)$ is the \emph{scalar quadratic fold curvature} at a fold point $p$ (see Table~\ref{tab:usefuldefinitions} and Figure~\ref{fig:geometricdefinitions} c,d)). $K[g](p) > 0$ corresponds to a quadratically convex fold with respect to the fold direction and $K[g](p) < 0$ corresponds to a quadratically concave fold. Hence
\begin{equation*}
K[g](p_1) + K[g](p_2)<0
\end{equation*}
means that the concave curvature dominates, and the degeneracy is called a \emph{covering fold projection tangency} since the folds locally cover the slow plane (see Figure~\ref{fig:foldtangencysketch}). Similarly, if
\begin{equation*}
K[g](p_1) + K[g](p_2)>0,
\end{equation*}
then the degeneracy is called a \emph{non-covering fold projection tangency}. Accounting for whether the fold umbrae interact with each other, or one fold umbra interacts with a fold, or neither, we get six subcases of opposed fold projection tangency (Table \ref{tab:globaldegeneracies12P2Folds}). If the two fold projections are aligned the total curvature does not matter as long as $K[g](p_1)\neq K[g](p_2)$, with one exception. This exceptional case occurs if a fold umbra hits a fold, in which case it matters if the curvature of the umbral fold dominates the destination fold or not (Figure~\ref{fig:foldtangencysketch} a.i) and a.ii)). More details are listed in \ref{app:global12tables}.

\subsubsection{Cusp-fold projection intersection}

At a cusp-fold projection intersection, the projections of a cusp and a fold line coincide in their projection onto the slow variables. We classify the intersection $\cD_{5}[g]$ of a cusp and a fold projection into ten cases (Table \ref{tab:globaldegeneracies12P2Cusps}), depending on the stability of the cusp (determined by $g_{xxx}$), the direction from which the cusp approaches the fold (determined by the sign of $\nu[g](p_2)\cdot \mu[g](p_2)$), and $k$, the number of regular sheets of equilibria separating the fold and cusp. If $k=0$, then one umbra intersect directly with a fold or cusp point (e.g. Figure~\ref{fig:structureCM121DSample}). If $k=1$, then two umbrae intersect, and if $k\geq 2$ then none of the umbrae or folds intersect. The middle columns of Figures \ref{fig:structureCM121D5Aligned} and \ref{fig:structureCM121D5Opposed} show typical cases of these degeneracies. Note that no degeneracies involving the umbrae of a stable cusp exist, since stable cusps have no umbrae.

\subsubsection{Triple fold projection intersection}

\begin{figure}
\begin{center}
\includegraphics[width=6cm]{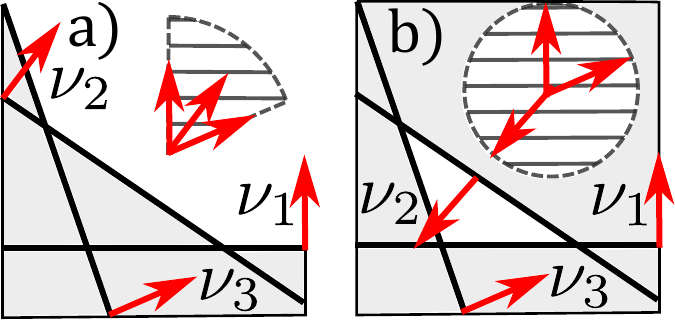}
\end{center}
\caption{(Color online) Sheets of the critical manifold near a) Non-covering and b) covering triple limit points bifurcations, projected onto the slow variables. Solid black lines show folds. Red arrows indicate direction vectors of folds $\nu_i, i\in\{1,2,3\}$, and grey areas indicate overlapping folds. The convex cones spanned by the direction vectors are shown as striped regions}
\label{fig:triplelimitpointcoveringnoncovering}
\end{figure}

The projections of three fold lines $\cD_{6}[g]$ onto the slow variables can intersect transversally in two ways: as a \emph{covering triple limit} or as a \emph{non-covering triple limit} (see Figure~\ref{fig:triplelimitpointcoveringnoncovering} b)). For brevity we write $\nu_i:=\nu[g](p_i)$. In the covering case, all folds are opposed in the sense that their direction vectors span a convex cone covering all of $\bR^2$. Therefore, the zero vector can be written as a linear combination of the direction vectors using only non-negative coefficients $a_i\geq 0$, not all zero:
\begin{equation}
\nu_1 a_1 + \nu_2 a_2 + \nu_3 a_3 = 0.
\label{eq:linearcombination}
\end{equation}
In the non-covering case the convex cone of the direction vectors does not cover $\bR^2$, meaning that at least one coefficient has to be negative in order for the vector sum to be zero (see Figure~\ref{fig:triplelimitpointcoveringnoncovering} a)). Therefore, the two subcases are defined by the signs of the coefficients in (\ref{eq:linearcombination})
\begin{equation}
    \left\{\begin{array}{rl}
        \mbox{Non-covering triple limit} & \mbox{ if } \pm\sign{(a_1,a_2,a_3)} = (+,+,-)\\
        \mbox{Covering triple limit} & \mbox{ if } \pm\sign{(a_1,a_2,a_3)} = (+,+,+)\\
    \end{array}\right.,
\end{equation}
for some choice of prefactor sign. Note that a higher codimension degeneracy will occur if $a_i=0$ for at least one $i$. Interactions of umbrae of the folds with other folds or umbrae give additional subclasses of triple limit points: these cases are detailed in \ref{app:global12tables} in Table~\ref{tab:globaldegeneracies12P3}. Note that it is not possible for all three fold umbrae to intersect.

We summarise the discussion in this section with the following classification of codimension one bifurcations for the case of one fast and two slow variables, analogous to Proposition~\ref{prop:CM111}.

\begin{figure}
\begin{center}
\includegraphics[width=15cm]{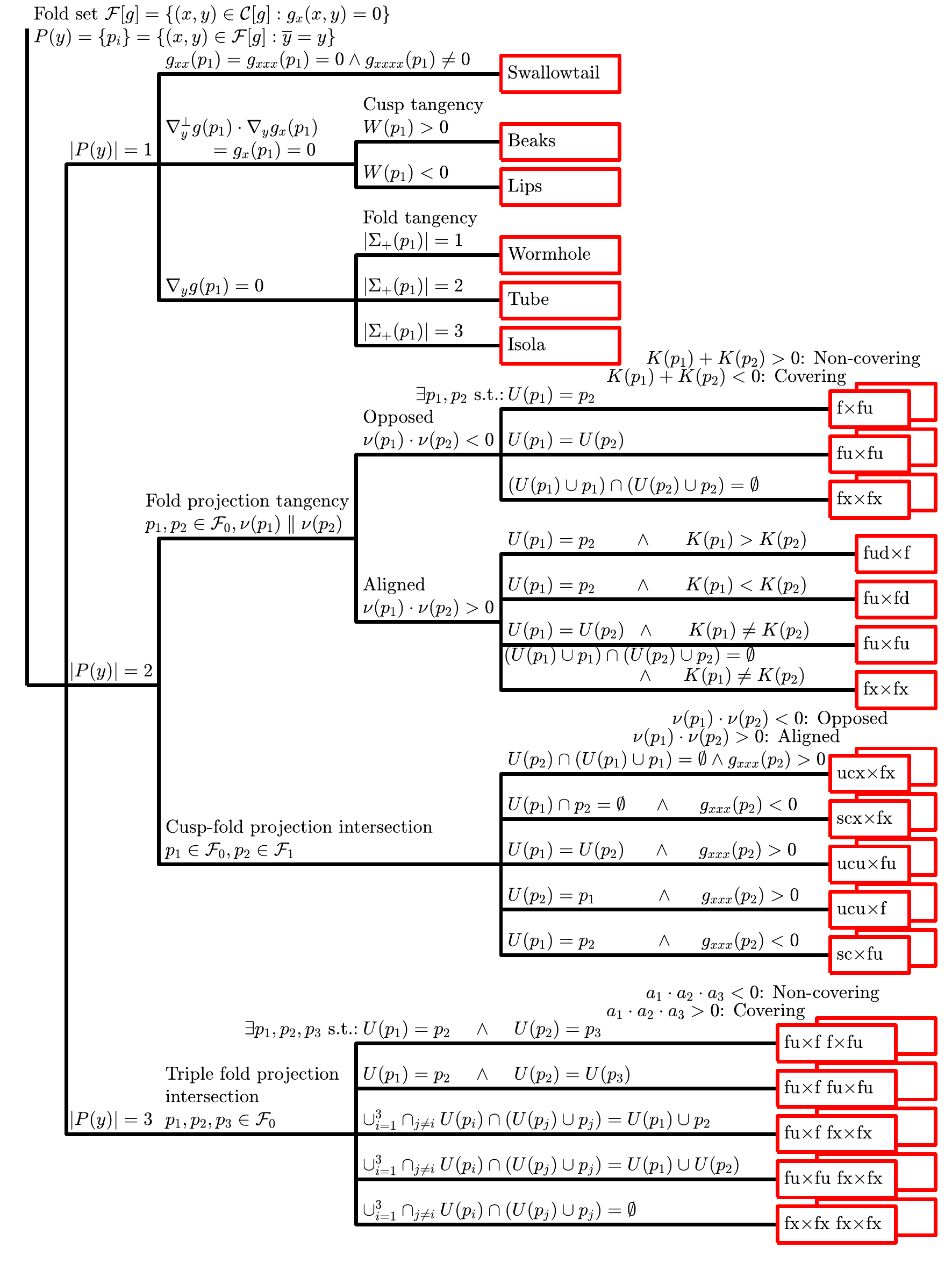}
\end{center}
\caption{
(Color online) Classification of codimension one degeneracies of the critical set for $m=1$ and $n=2$. $\nu$ and $\mu$ are direction vectors of folds and cusps, and $K$, $a_i$, $W$ and $\Sigma_+$ are described in Section \ref{sec:codimonecritsed12} and the appendix. Additionally, f means fold, sc/uc means stable/unstable cusp, fu/cu means fold/cusp umbra and fx/cx means non-interacting fold/cusp. Similarly, fd means fold dominant and fud means fold umbra dominant. We suppress dependence on $g$ e.g. $\nu[g](p)=\nu(p)$. Each red box corresponds to one persistent codimension one bifurcation, if degeneracy occurs for one $y$ only. Overlapping red boxes means there are two subcases separated by a condition. See the text for details}
\label{fig:degeneracytable12}
\end{figure}

\begin{conj}[Codimenson one, $m=1,n=2$.]
\label{conj:CM121}
For $m=1$ and $n=2$ the codimension one bifurcations of critical sets $\cC[g]$ for $g\in V_{f}$ are characterised in Figure~\ref{fig:degeneracytable12}, such that precisely one of the sets $\cD_{j,k}[g,y]$ in Table~\ref{tab:D11persistentsubcases} is non-empty, for precisely one $y\in\bR^2$. At such a bifurcation, precisely one of the following occurs:
\begin{enumerate}
\item A loop or pair of hyperbolae appears in the fold projections at a fold tangency $\cD_{1,k}$. [e.g. Fig.~\ref{fig:structureCM121DSample} a) to i)]f
\item Two cusps annihilate at a cusp tangency $\cD_{2,k}$. [e.g. Fig.~\ref{fig:structureCM121DSample} j) to o)]
\item A quadratic fold line folds over to form two cusps in a swallowtail $\cD_{3,k}$. [e.g. Fig.~\ref{fig:structureCM121DSample} p,q,r)]
\item The projections of two quadratic fold curves onto the slow variables become tangent $\cD_{4,k}$. [e.g. Fig.~\ref{fig:structureCM121DSample} s,t,u)]
\item The projections of a quadratic fold curve and a cubic cusp intersect $\cD_{5,k}$. [e.g. Fig.~\ref{fig:structureCM121DSample} v,w,aa)]
\item The projections of three fold lines intersect $\cD_{6,k}$. [e.g. Fig.~\ref{fig:structureCM121DSample} bb,cc,dd)]
\end{enumerate}
\end{conj}

\begin{figure}
\begin{center}
\includegraphics[width=15cm]{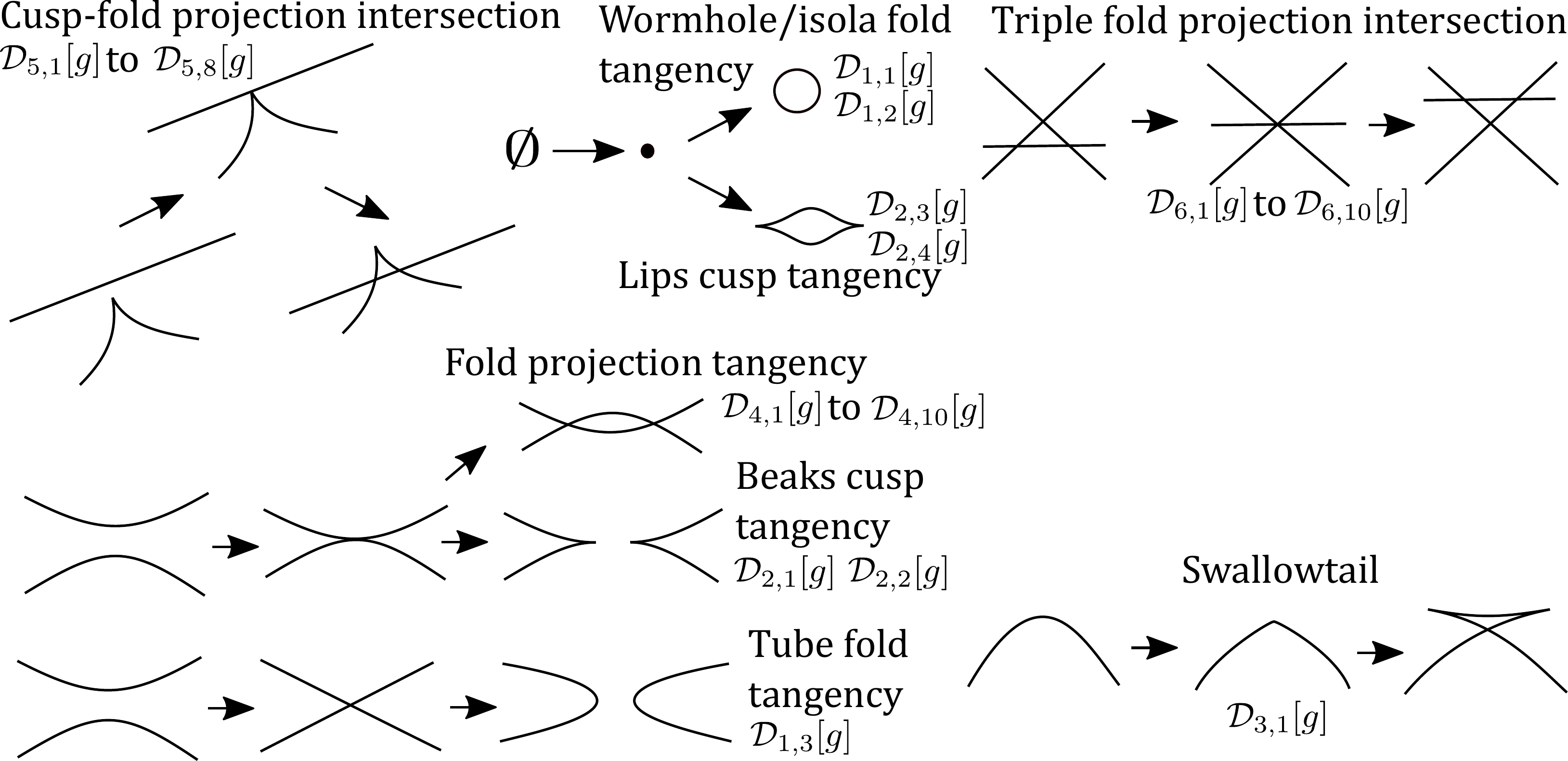}
\end{center}
\caption{(Color online) Types of codimension one bifurcations of the critical set for one fast and two slow variables, shown in terms of changes to the fold set projected onto slow variables. The cases are enumerated more precisely in Figure~\ref{fig:degeneracytable12}
}
\label{fig:projectionsketch}
\end{figure}

Figure~\ref{fig:projectionsketch} shows shows the projections of fold lines and cusps that correspond to the possible codimension one degeneracies of the critical set.  \ref{app:global12tables} gives a detailed listing of inequivalent subcases of codimension one bifurcations associated with projection intersection: we do not attempt to suggest global normal forms for these cases.

\FloatBarrier

\section{Global singular equivalence, persistence and bifurcation} 

\subsection{Global singular equivalence of systems}

To define a useful notion of global equivalence of system (\ref{eq:mainsystem}) in the singular limit, we fix compact regions $M$ and $N$ as above and suppose that $\{g,h\}$ and $\{\tilde{g},\tilde{h}\}$ are both in $V_f\times V_s$ where these are defined as in the previous section. 

We say $\{g,h\}$ is {\em globally singularly equivalent} to $\{\tilde{g},\tilde{h}\}$ (on $M\times N$) if one can write
\begin{equation}
\left\{\begin{array}{rl}
\tilde{g}(x,y) = S(x,y) g(X(x,y),Y(y))& \mbox{ for all }(x,y)\in M\times N\\
\tilde{h}(x,y) = T(x,y) h(X(x,y),Y(y)) & \mbox{ for all }(x,y)\in \cC[\tilde{g}]
\end{array}\right.
\label{eq:globalsingequiv}
\end{equation}
where:
\begin{itemize}
\item The map $\Phi(x,y)=(X(x,y),Y(y))$ is a diffeomorphism on $M\times N$.
\item The function $S(x,y)>0$ is smooth and positive on $M\times N$.
\item The function $T(x,y)>0$ is smooth and positive on $M\times N$.
\end{itemize}
Note that because we are only interested in equivalence of the singular systems, we allow independent re-parametrization of the fast and slow timescales. Note that $T(x,y)$ is globally defined but only evaluated on $\cC[\tilde{g}]$. Clearly, if $\{g,h\}$ is globally singularly equivalent to $\{\tilde{g},\tilde{h}\}$ then $g$ is globally equivalent to $\tilde{g}$ in the sense of (\ref{eq:globalequivalence}). One can check that this is an equivalence relation - it is transitive and reflexive, and one can check it is symmetric by noting that if (\ref{eq:globalsingequiv}) holds then
\begin{equation}
\left\{\begin{array}{rl}
g(x,y) = \tilde{S}(x,y) \tilde{g}(\tilde{X}(x,y),\tilde{Y}(y))& \mbox{ for all }(x,y)\in M\times N\\
h(x,y) = \tilde{T}(x,y) \tilde{h}(\tilde{X}(x,y),\tilde{Y}(y)) & \mbox{ for all }(x,y)\in \cC[g]
\end{array}\right.
\end{equation}
because $(x,y)\in\cC[g]$ if and only if $\Phi(x,y)\in \cC[\tilde{g}]$, and one can verify that:
\begin{itemize}
\item The map $\tilde{\Phi}(x,y)=(\tilde{X}(x,y),\tilde{Y}(y))$ is a diffeomorphism that is the inverse of $\Phi$ on $M\times N$.
\item The function $\tilde{S}(x,y)= 1/S(\tilde{\Phi}(x,y))$ is smooth and positive on $M\times N$.
\item The function $\tilde{T}(x,y)=1/T(\tilde{\Phi}(x,y))$ is smooth and positive on $M\times N$.
\end{itemize}
Note that singular trajectories are mapped onto each other by global singular equivalence as expressed in the following result.

\begin{lemma}
Suppose that $\{g,h\}$ is globally singularly equivalent to $\{\tilde{g},\tilde{h}\}$ on $M\times N$. Then the singular trajectories of these systems are equivalent via a diffeomorphism.
\end{lemma}

\proof
To see this, suppose that $\Phi,S,T$ are found that satisfy (\ref{eq:globalsingequiv}) and suppose that $\gamma_0:[a,b]\rightarrow M\times N$ is a singular trajectory for $\{g,h\}$  as in Definition~\ref{def:singulartraj} for $a=s_1<\cdots<s_m=b$. If $J_j=(\tilde{s}_j,\tilde{s}_{j+1})$ is any fast trajectory segment then $\tilde{\Phi}(J_j)$ is a fast trajectory segment for $\{\tilde{g},\tilde{h}\}$ with the same orientation (and time scaled by $\tilde{S}$). If $J_j$ is a slow segment then it lies within $\cC[g]$ and so $\tilde{\Phi}(J_j)$ is a slow trajectory segment for $\{\tilde{g},\tilde{h}\}$ that lies within $\cC[\tilde{g}]$ with the same orientation (and time scaled by $\tilde{T}$).
\qed

\subsection{Persistence under global singular equivalence}

If system (\ref{eq:mainsystem}) is its own universal unfolding under global singular equivalence then we say the system is {\em persistent}. Clearly, in such a case the fast vector fields indexed by the slow variables must be persistent, but also we cannot have degeneracies of the slow system on the critical set. We expect (\ref{eq:mainsystem}) to be persistent under global singular equivalence if the fast subsystem is persistent and, in addition, the slow system has persistent behaviour on the critical set.

For the case $m=n=1$ we can make this statement more precise. We define the slow nullcline
\begin{equation*}
\cN[h]=\{(x,y)~:~h(x,y)=0\}.
\end{equation*}
For any regular point $p_0=(x_0,y_0)\in\cC_{reg}[g]$ there will be a curve $(X_{p_0}(y),y)\in \cC[g]$ with $X_{p_0}(x_0)=y_0$ such that $g(X_{p_0}(y),y)=0$. Implicitly differentiating this gives
\begin{equation*}
\frac{dX_{p_0}}{dy}(y)=-\frac{g_y(X_{p_0}(y),y)}{g_x(X_{p_0}(y),y)}
\end{equation*}
for $y$ close to $y_0$. Then we can locally reduce (\ref{eq:reduced}) to an equation on the critical set of the form
\begin{equation*}
\dot{y}=H_{p_0}(y):=h(X_{p_0}(y),y).
\end{equation*}
If $p_0=(x_0,y_0)\in \cN[h]\cap\cC_{reg}[g]$ then
$H_{p_0}(y_0)=0$ is an equilibrium and its linear stability is determined via
\begin{equation*}
H_{p_0}'=\frac{dH_{p_0}}{dy}=-h_x\frac{g_y}{g_x}+h_y
\end{equation*}
evaluated at $p_0$. This highlights that the slow dynamics are essentially one-dimensional when restricted to $\cC_{reg}[g]$. Before stating a result on persistence of fast-slow systems, we give some definitions. We define the restriction of the slow nullcline onto the critical manifold as
\begin{equation*}
\cN_r[g,h]=\cN[h] \cap \cC[g],
\end{equation*}
and define the slow degenerate set $\cE[g,h]$ as the union of two subsets, defined shortly: $\cE[g,h] = \cE_1[g,h] \cup \cE_2[g,h]$. The slow locally degenerate set is
\begin{equation*}
\cE_1[g,h]=\{p=(x,y)\in\cN_r[g,h]~:~g_x(p)h_y(p) - h_x(p)g_y(p)=0\},
\end{equation*}
which occurs when the fast and the slow nullclines intersect tangentially. Note that this condition is equivalent to the determinant of the Jacobian of the full system being zero, and for $p_0\in\cC_{reg}[g]$ this implies that $H_{p_0}'(y)=0$.

Recalling that $\pi(p)$ is the projection onto the slow variables, we define the set of {\em slow co-equilibria}
\begin{equation*}
\Xi[g,h](p)=\{\pi^{-1}(\pi(p))\cap \cN_r[g,h]\},
\end{equation*}
which we use to define the {\em multiple slow equilibrium} set
\begin{equation*}
\cE_2[g,h]=\{p\in \cN_r[g,h]~:~|\Xi[g,h](p)|\geq2\}.
\end{equation*}

We define the (mixed) fold-equilibrium multiple projection set as
\begin{equation*}
\cM[g,h]=\{p\in \cN_r[g,h]~:~\pi(p)\cap \pi(\cF[g]) \neq \emptyset\},
\end{equation*}
that is, the set of equilibria that share slow coordinate with at least one fold of the critical manifold. This finally allows us to define the fast-slow degenerate set as
\begin{equation*}
\cG[g,h]=\cD[g]\cup \cE[g,h]\cup \cM[g,h].
\end{equation*}
Theorem~\ref{thm:globalsingpersist} below establishes that this set contains all degeneracies under global singular equivalence.

\begin{thm}
\label{thm:globalsingpersist}
In the case $m=n=1$, if $g\in V_f$ and $h\in V_s$ then (\ref{eq:mainsystem}) is persistent under global singular equivalence for $\{g,h\}$ if and only if the all of the following hold (i.e. $\cG[g,h]=\emptyset$):
\begin{enumerate}
\item The critical set $\cC[g]$ has no degenerate folds (i.e. $\cD[g]=\emptyset$).
\item There is at most one equilibrium per slow coordinate $y$ (i.e. $\cE_2[g,h]=\emptyset$)
\item There is no intersection of the slow nullcline and folds or co-folds, i.e. ($\cM[g,h]=\emptyset$).
\item There are no degenerate slow equilibria on the critical manifold ($\cE_1[g,h]=\emptyset$)
\end{enumerate}
\end{thm}

\proof
We begin with the ``if'' part. If the critical set has degenerate folds $\cD[g]\neq\emptyset$ then $g$ is non-persistent under global equivalence. Hence, we need that $\cD[g]=\emptyset$. Assume by contradiction that $\cE_2[g,h]\neq\emptyset$. Then at least two equilibria share slow coordinate $y_1$. Under a generic perturbation of $\{g,h\}$ these will have different slow coordinates, and since the base is preserved under global singular equivalence, $\{g,h\}$ cannot be deformed to make the equilibria share $y$ coordinate. Hence, for persistence we need $\cE_2[g,h]=\emptyset$. A similar argument implies that $\cM[g,h]=\emptyset$ is required for persistence. Given that $\cM[g,h]=\emptyset$, $\cE_1[g,h]\neq\emptyset$ implies that there is a $p_0\in \cC_{reg}[g]$ and a $y$ such that $H_{p_0}(y)=H'_{p_0}(y)=0$. But this is a non-hyperbolic equilibrium, and thus  is not persistent to perturbation. Hence, for persistence we need $\cE_1[g,h]=\emptyset$.

For the ``only if'' part, we need to argue that there is no other way for (\ref{eq:mainsystem}) to be non-persistent than if $\cG[g,h]\neq\emptyset$. Assume that $\cG[g,h]=\emptyset$. Then there is a neighbourhood of every singularity of $g$ and equilibrium of $\{g,h\}$ that has only one fold or equilibrium in the fast fibre. These are either quadratic folds of $g$ or hyperbolic equilibria of $\{g,h\}$, both which are persistent under perturbation. Hence, $\{g,h\}$ must be persistent.
\qed

Note that the global singular equivalence (\ref{eq:globalsingequiv}) does not depend on the nullcline $\cN(x,y)$ away from the critical manifold. Bifurcation occurs if one of the assumptions in Theorem~\ref{thm:globalsingpersist} is broken. Note that the assumptions that $g\in V_f$ implies the critical set does not intersect $\partial M\times N$, and that $h\in V_s$ implies that the nullcline does not intersect $M\times \partial N$; more generally there will be additional persistence conditions that require persistent intersection with these boundaries. 

\subsection{Generic bifurcations in singular fast-slow systems}

We can understand generic codimension one bifurcations of fast-slow systems (\ref{eq:mainsystem}) by examining the ways that the persistence conditions of Proposition~\ref{prop:CM110} are violated. For $m=1$ and $n=1,2$ this means that the codimension one bifurcations of the critical manifold are possible bifurcations under global singular equivalence. In addition, there are many ways that a change in the slow subsystem can lead to a bifurcation.

For $m=n=1$ we define, as for the critical manifold, subsets of $\cE_1[g,h]$, $\cE_2[g,h]$ and $\cM[g,h]$ containing all codimension one degeneracies of $\cG[g,h]$, which are not only due to degeneracy of the fast subsystem, in Table \ref{tab:EM11persistent}. We further define subsets of these, which give codimension one bifurcation if all except for one of the subsets are empty, and if the nonempty subset is nonempty for only one slow coordinate $y$ (Table \ref{tab:D11persistentsubcases}). Equipped with the subsets in Table \ref{tab:slowcodim1subcases11}, Proposition~\ref{prop:fastslowbifcodim1} lists the codimension one degeneracies of (\ref{eq:mainsystem}) for $m=n=1$.

\begin{table}
\caption{Subsets of $\cE_1[g,h]$, $\cE_2[g,h]$ and $\cM[g,h]$ for $m=n=1$ whose union contains all codimension one bifurcations, not only due to bifurcation in the fast subsystem. Note that $\det(D \{g,h\})=g_{x}h_{y}-g_{y}h_{x}$ is the Jacobian of the full system (\ref{eq:mainsystem}) and that the first subset is a local degeneracy
}
\centering
\begin{tabular}{|p{0.25\linewidth}|p{0.65\linewidth}|}
\hline
Saddle-node: & $\cE_{1}^1[g,h]~=~\{p\in \cE_1[g,h]:|\Xi[g,h](p)|=1$ and $\det(D\{g,h\}(p))\neq 0\}$\\
\hline
Double slow equilibrium: & $\cE_{2}^1[g,h]~=~\{p\in \cE_2[g,h]:|\Xi[g,h](p)|=2$\}\\
\hline
Fold-equilibrium~double~projection~set: & $\cM^1[g,h]~=~ \{p\in \cM[g,h]: |\Xi[g,h](p)|=1\}$\\
\hline
\end{tabular}
\label{tab:EM11persistent}
\end{table}

\begin{table}
\caption{Degeneracies that lead to codimension one bifurcation of the singular fast-slow system (\ref{eq:mainsystem}) up to global singular equivalence. The sets $\cF[g]$ and $\cU[g]$ are the fold and the umbral sets and $\pi(p)$ is the projection map. $R(y)$ is the set of all equilibria sharing slow coordinate $y$. In saddle-node non-degeneracy condition, $r$ and $q$ are eigenvectors of the Jacobian of the full system and its adjoint respectively, and $B=\sum_{j,k}^2q_j q_k\frac{\partial^2}{\partial \xi_1 \partial \xi_2}(g,h)$, where $(\xi_1,\xi_2)=(x,y)$ \cite[p.175]{kuznetsov04}. The last column associates the degeneracy to a possible bifurcation of relaxation oscillations. Degeneracies which do not lead to bifurcation of singular relaxation oscillations have no associated figures}
\begin{tabular}{|p{0.2\linewidth}|p{0.5\linewidth}|p{0.2\linewidth}|}
\hline
Non-degenerate saddle-node & $\cE_{1,1}[g,h,y]=\{R(y)=\pi^{-1}(y)\cap \cN_r[g,h] \subset \cE^1_{1}[g,y]~:~|R(y)|=1$ and $b = {\frac{1}{2}(r \cdot B(q,q))} \neq 0\}$ & Fig.~\ref{fig:singularbifCM111slowsub} a,b,c), Saddle-node on invariant circle (SNIC) \cite{ermentrout86} \\
\hline
Double slow equilibrium: & $\cE_{2,1}[g,h,y]=\{R(y)\subset \cE_2^1[g]:|R(y)|=1\}$ & \\
\hline
Sink-fold intersection: & $\cM_{1,1}[g,h,y]=\{R(y)\subset \cM^1[g,h]~:~|R(y)\cap \cF[g]|=1$ and $H'_{p_0}(y)<0$ for $p_0\in R(y)\}$ & Fig.~\ref{fig:singularbifCM111slowsub} d,e,f), Singular Hopf\\
\hline
Source-fold intersection: & $\cM_{1,2}[g,h,y]=\{R(y)\subset \cM^1[g,h]~:~|R(y)\cap \cF[g]|=1$ and $H'_{p_0}(y)>0$ for $p_0\in R(y)\}$ & \\
\hline
Sink-fold umbra intersection: & $\cM_{1,3}[g,y]=\{R(y)\subset \cM^1[g,h]~:~|R(y)\cap \cU[g]|=1$ and $H'_{p_0}(y)<0$ for $p_0\in R(y)\}$ & \\
\hline
Source-fold umbra intersection: & $\cM_{1,4}[g,h,y]=\{R(y)\subset \cM^1[g,h]~:~|R(y)\cap \cU[g]|=1$ and $H'_{p_0}(y)>0$ for $p_0\in R(y)\}$ & Fig.~\ref{fig:singularbifCM111slowsub} g,h,i), Singular homoclinic\\
\hline
Non-interacting source-fold umbra: & $\cM_{1,5}[g,h,y]=\{R(y)\subset \cM^1[g,h]~:~R(y)\cap (\cU[g]\cup\cF[g])= \emptyset\}$ & \\
\hline
\end{tabular}
\label{tab:slowcodim1subcases11}
\end{table}

\begin{prop}
In the case $m=n=1$, codimension one bifurcation of the fast-slow system (\ref{eq:mainsystem}) for $\{g,h\}$ occurs due to exactly one of the following reasons, for exactly one slow coordinate $y\in N$. 
\begin{enumerate}
    \item Two folds of $\cC[g]$ merge at a quadratic fold tangency of the critical manifold at some $(x,y)$, and $\cG[g,h]\setminus \cD^1_{1}[g,y]=\emptyset$.  \label{itm:codimonero1}
    \item There is a cubic hysteresis of $\cC[g]$ at some $(x,y)$, and $\cG[g,h]\setminus \cD^1_{2}[g,y]=\emptyset$.  \label{itm:codimonero2}
    \item There is a double limit point degeneracy of $\cC[g]$ for some $y$ and $\cG[g,h]\setminus \cD^1_3[g]=\emptyset$  \label{itm:codimonero3}
    \item There is a nondegenerate slow saddle-node equilibrium on the regular part of the critical manifold, and $\cG[g,h]\setminus \cE_{1,1}[g,h,y]=\emptyset$  \label{itm:codimonero4}
    \item There are exactly two hyperbolic equilibria that share the same slow coordinate $y$ (i.e. there are exactly two points $p_1,p_2\in\cE_2[g,h]$ for which $\pi(p_1)=\pi(p_2)=y$, and $\cG[g,h]\setminus \cE_{2,1}[g,h,y]=\emptyset$.
     \label{itm:codimonero5}
    \item The slow nullcline intersects the critical set transversally at exactly one point $(x,y)$ that shares slow coordinate with a quadratic fold, and $\cG[g,h]\setminus \cM^1[g,h]=\emptyset$
     \label{itm:codimonero6}
\end{enumerate}
\label{prop:fastslowbifcodim1}
\end{prop}

\proof
Note that all degeneracies are contained in the set $\cG[x,y]=\cD_1[g]\cup\cD_2[g]\cup\cD_3[g]\cup \cE_1[g,h]\cup\cE_2[g,h]\cup\cM[g,h]$. Because of this, and because the defining conditions are independent, codimension one degeneracy will occur at a point that is in exactly one of those sets. Furthermore, bifurcation must occur for exactly one $y$ since otherwise more than one equality constraint is imposed, raising the codimension.

Case \ref{itm:codimonero1}) describes the only subset $\cD_1^1[g]$ of $\cD_1[g]$ containing codimension one degeneracies exclusively in $\cD_1[g]$, and therefore it produces codimension one degeneracy of $\{g,h\}$. The same is true for cases \ref{itm:codimonero2}) and \ref{itm:codimonero3}). Case \ref{itm:codimonero4}) is codimension one since we impose just one equality condition and exclude higher codimension degeneracy with the non-degeneracy condition in Table \ref{tab:slowcodim1subcases11}. Case \ref{itm:codimonero5}) is codimension one since hyperbolic equilibria are persistent, and more than two hyperbolic equilibria sharing slow coordinate would impose more than one equality constraint. Case \ref{itm:codimonero6}) is codimension one for the same reason.
\qed

Note that codimension two bifurcations may combine degeneracies in more than one of these sets.

\section{Generic bifurcations of relaxation oscillations in singular fast-slow systems}
\label{sec:bifurcation}

Not all bifurcations of the singular fast-slow system listed in Proposition \ref{prop:fastslowbifcodim1} will lead to bifurcation of singular relaxation oscillations, as the degeneracy in the fast-slow system must interact with a limit cycle. We focus on bifurcation of singular relaxation oscillations and  \emph{simple} relaxation oscillations, a generic subclass due to \cite{guckenheimer96}. Several cases of these bifurcations have been considered in the literature, see for example \cite{maesschalck11}. 

\subsection{Singular relaxation oscillations}

Consider a fast-slow system with $m=1$ fast variables (\ref{eq:mainsystem}) in the singular limit $\epsilon=0$. A relaxation oscillation is a singular periodic trajectory $\gamma:[a,b]\rightarrow M\times N$ (i.e. such that $\gamma(b)=\gamma(a)$) where the slow segments are in $\overline{\cC_{att}}$.  If the oscillation consists of alternating stable slow segments on $\cC_{att}[g]$ up to non-degenerate folds, fast segments from these folds to their umbra, and satisfies certain other non-degeneracy conditions then we say it is a simple relaxation oscillation. These are called {\em strongly common slow-fast cycles} in \cite{maesschalck11} where it is shown that these singular trajectories will be shadowed by a stable periodic orbit for small enough $\epsilon$. Guckenheimer stated a similar persistence theorem in \cite{guckenheimer96}; in Section \ref{sec:persistbif11} we state and prove a version of it.

We say a continuous curve $s_k:[0,1]\rightarrow M\times N$ is a {\em slow segment} of a singular trajectory if there is a continuous and monotonic increasing $\theta:[0,1]\rightarrow \bR$ such that $s_k(\theta(s))$ is a trajectory of (\ref{eq:reduced}). We say a slow segment $s_k$ has {\em slow time duration} $T_k>0$ if it can be parameterised by $\theta(s)=s/T_k$. If not, and if $\theta(0)=\theta(1)$ then $T_k=0$, otherwise $T_k=\infty$.

Up to equivalence of the fast segments joining the slow segment end-points, we define a {\em relaxation oscillation} in terms of its slow segments as
\begin{equation}
\cA=\{s_k(\theta)~:~\theta\in[0,1]\}_{k=0}^{d-1}
\label{eq:RA}
\end{equation}
a sequence of continuously parametrized slow segments $s_k:[0,1]\rightarrow M\times N$. We will assume that {\em either} $\cA$ is a loop entirely within $\overline{\cC_{att}}[g]$ {\em or}
\begin{itemize}
    \item $s_k(\theta)\subset \overline{\cC_{att}}[g]$ for all $\theta\in(0,1)$
    \item There is a trajectory $\phi(t)$ of the fast system such that $\alpha(\phi(0))=s_k(1)$ and $\omega(\phi(0))=s_{k+1}(0)$ for $k$ modulo $d$.
\end{itemize}
for each $k$, where $\omega(p)$ and $\alpha(p)$ are the omega and alpha limits of a point $p$ respectively. This equivalence class of relaxation oscillations has more than one member if there is more than one fast segment joining two consecutive slow segments.

We define the {\em slow period} $\cP(\cA)$ of a relaxation oscillation $\cA$ to be the total slow time duration of its slow segments. This is
\begin{equation}
\cP(\cA)=\sum_{k=0}^{d-1} T_k
\end{equation}
which may be infinite, where orbits in the equivalence class of $\cA$ clearly all have the same slow period. We allow the possibility that $\cA$ is a loop on $\overline{\cC_{att}}[g]$ without jumps, in which case $d=1$ and $s_0(1)=s_0(0)$, or that the jumps are trivial and on $\cC_{att}[g]$. Infinite slow period relaxation oscillation ($\cP(\cA)=\infty$) of a variety of types are covered by this definition. 
We define a \emph{simple} relaxation oscillation (cf Guckenheimer \cite{guckenheimer02,guckenheimer95}) as follows:

\begin{definition}[Simple relaxation oscillation ($m=1,n\geq1$)]
A relaxation oscillation $\cA$ in (\ref{eq:RA}) is \emph{simple} if all of the following hold:
\begin{enumerate}[label=\roman*)]
\item The slow period $\cP(\cA)$ is finite.
\item The slow segments are on $\cC_{att}[g]$, except possibly the last point.
\item Either $s_k(1)\in \cF[g]\setminus \cD[g]$, or $d=1$ and $s_0(1)=s_0(0)$. \label{itm:sro4}
\item The slow segments are not tangent to either fold set or umbral set.
\item The singular return map local to $s_k(1)$ is well-defined with a hyperbolic equilibrium at $s_k(1)$. 
\label{itm:sro3}
\end{enumerate}
\label{def:simplero}
\end{definition}

Note that the assumption $P(\cA)<\infty$ implies that the slow segments do not limit to any equilibria of the slow flow. For one fast and one slow variable, our definition of a simple relaxation oscillation can be expressed in a simpler way:

\begin{definition}[Simple relaxation oscillation $(m=n=1)$]
A relaxation oscillation (\ref{eq:RA}) with one fast and one slow variable ($m=n=1$) is \emph{simple} if
\begin{enumerate}[label=\roman*)]
\item The slow period $\cP(\cA)$ is finite. 
\item We have $s_k(\theta)\in\cC_{att}[g]$ for all $\theta\in[0,1)$.
\item Either $s_k(1)\in \cF[g]\setminus\cD[g]$ or $d=1$ and $s_0(1)=s_0(0)$.
\end{enumerate}
\label{def:simplero11}
\end{definition}

\subsection{Persistence and bifurcation}
\label{sec:persistbif11}

We say that a simple relaxation oscillation undergoes bifurcation if the relaxation oscillation ceases to be simple under perturbation of the singular fast-slow system. If not, we say that the relaxation oscillation is persistent. Note that as we only consider fast-slow systems on absorbing regions in $\bR^2$, singular relaxation oscillations can bifurcate to either equilibrium points or other singular relaxation oscillations. The following proposition links bifurcation of simple relaxation oscillations to degeneracy in the fast-slow system (cf \cite{guckenheimer96,guckenheimer02,guckenheimer04}):

\begin{prop}
A simple relaxation oscillation $\cA$ is persistent for $n=m=1$ if the fast-slow system is persistent under global singular equivalence.
\label{prop:fastslowpersimpliessropers}
\end{prop}
\proof
Assume $\{g,h\}$ is persistent. Because the slow period is finite, no slow equilibrium can intersect a slow segment $s_k(\theta)$ in a degenerate way; intersection for $\theta\in\{0,1\}$ implies that $\cM[g,h]\neq\emptyset$ and intersection for $\theta\in(0,1)$ implies $\cE_1[g,h]\neq\emptyset$, since the flow must have the same direction on both sides of the equilibrium, which is not possible if the equilibrium is hyperbolic. Hence Condition i) is persistent. Condition ii) is persistent since $s_k(0)\not\in \cC_{att}[g]$ implies that both $s_{k-1}(1)$ and $s_k(0)$ are in $\cF[g]$, which in turn would imply that $\cD_3[g]\neq\emptyset$. Condition iii) is persistent since the fold set is non-degenerate everywhere, and since trivial relaxation oscillations $s_0(1)=s_0(0)$ coincide with hyperbolic equilibria in the interior of $\overline{C_{att}[g]}$.
\qed

\begin{figure}
\begin{center}
\includegraphics[width=10cm]{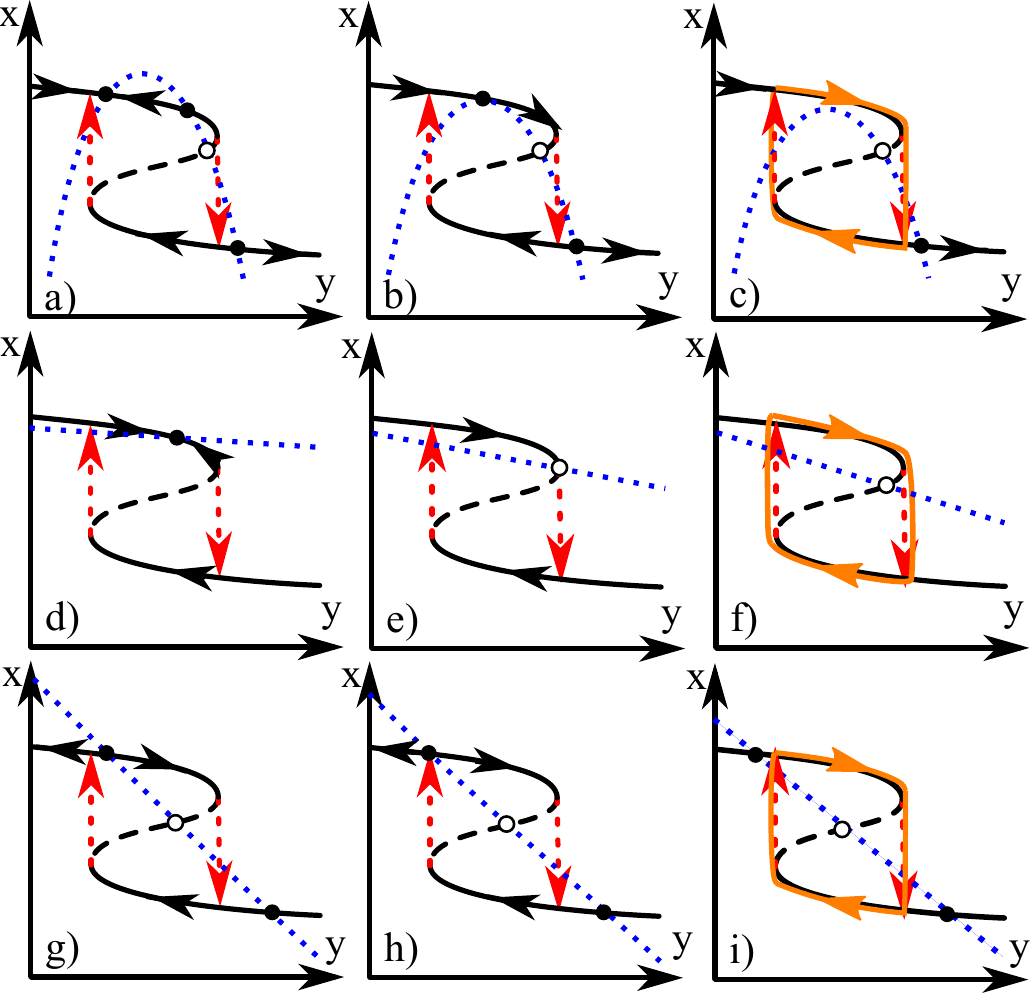}
\end{center}
\caption{
(Color online) The middle column shows typical examples of the classes of codimension one bifurcations of equilibria for $m=n=1$ that are not due to bifurcation of the critical set (Table \ref{tab:bifofRO11}). In all cases solid black lines show the critical set, red lines show the image of the fold under the umbral map, blue lines show nullclines of the slow variable and orange lines show stylised solutions of the non-singular system. Filled/open dots are stable/unstable equilibria of the fast subsystem}
\label{fig:singularbifCM111slowsub}
\end{figure}

The converse is not true, however: only some of the degeneracies in the fast-slow system will give rise to a bifurcation of a simple relaxation oscillation. The reasons are listed in the following Proposition and in Table~\ref{tab:bifofRO11}. Examples are portrayed in Figures~\ref{fig:singularbifCM111slowsub} and Figure~\ref{fig:slowfastbifexamples}.

\begin{table}
\caption{Bifurcations of singular relaxation oscillations for $m=n=1$ due to codimension one bifurcation of the fast-slow system. Definitions of the calligraphic sets are found in Table \ref{tab:D11persistentsubcases} and Table \ref{tab:slowcodim1subcases11}. It is understood that $\cA(\lambda)$ is perturbed due to perturbed $g(\lambda)=g(x,y,\lambda)$ and $h(\lambda)=h(x,y,\lambda)$. The final column indicates whether the bifurcation is due to bifurcation of the critical set
}
\centering
\begin{tabular}{|p{0.20\linewidth}|p{0.6\linewidth}|p{0.08\linewidth}|}
\hline
Type:\newline Example & Conditions for a simple singular relaxation oscillation $\cA(\lambda)$ to bifurcate at $\lambda=\lambda_0$ &  Bifur\-cation of $\cC[g]$? \\
\hhline{|=|=|=|}
Saddle node on invariant circle \cite{ermentrout86}:\newline Figure~\ref{fig:singularbifCM111slowsub} a,b,c) & There exists a unique $y\in N$ such that $\cA_{lim}(\lambda_0)\neq \emptyset$ and $\cE_{1,1}[g(\lambda_0),h(\lambda_0),y]\neq \emptyset$ and $\cG[g(\lambda_0),h(\lambda_0)]\setminus\cE_{1,1}[g(\lambda_0),h(\lambda_0),y]=\emptyset$,& No\\
\hline
Singular Hopf: \newline Figure~\ref{fig:singularbifCM111slowsub} d,e,f) & There exists a unique $y\in N$ such that $\cA_{lim}(\lambda_0)\neq \emptyset$ and $\cM_{1,1}[g(\lambda_0),h(\lambda_0),y]\neq \emptyset$ and $\cG[g(\lambda_0),h(\lambda_0)]\setminus\cM_{1,1}[g(\lambda_0),h(\lambda_0),y]=\emptyset$, & No\\
\hline
Singular homoclinic: \newline Figure~\ref{fig:singularbifCM111slowsub} g,h,i) & There exists a unique $y\in N$ such that $\cA_{lim}(\lambda_0)\neq \emptyset$ and $\cM_{1,4}[g(\lambda_0),h(\lambda_0),y]\neq \emptyset$ and $\cG[g(\lambda_0),h(\lambda_0)]\setminus\cM_{1,4}[g(\lambda_0),h(\lambda_0),y]=\emptyset$, & No\\
\hline
Hyperbolic fold tangency:\newline Figure~\ref{fig:slowfastbifexamples} a,b,c) & There exists a unique $y\in N$ such that $\cA_{lim}(\lambda_0)\neq \emptyset$ and $\cD_{1,1}[g(\lambda_0),y]\neq \emptyset$ and $\cG[g(\lambda_0),h(\lambda_0)]\setminus\cD_{1,1}[g(\lambda_0),y]=\emptyset$, & Yes\\
\hline
Hysteresis: \newline Figure~\ref{fig:slowfastbifexamples} d,e,f) & There exists a unique $y\in N$ such that $\cA_{lim}(\lambda_0)\neq \emptyset$ and $\cD_{2,1}[g(\lambda_0),y]\neq \emptyset$ and $\cG[g(\lambda_0),h(\lambda_0)]\setminus\cD_{2,1}[g(\lambda_0),y]=\emptyset$, & Yes\\
\hline
Aligned fold-fold umbra double limit:\newline Figure~\ref{fig:slowfastbifexamples} g,h,i) & There exists a unique $y\in N$ such that $\cA_{lim}(\lambda_0)\neq \emptyset$ and $\cD_{3,1}[g(\lambda_0),y]\neq \emptyset$ and $\cG[g(\lambda_0),h(\lambda_0)]\setminus\cD_{3,1}[g(\lambda_0),y]=\emptyset$, & Yes\\
\hline
Opposed fold-fold umbra double limit:\newline Figure~\ref{fig:slowfastbifexamples} j,k,l) & There exists a unique $y\in N$ such that $\cA_{lim}(\lambda_0)\neq \emptyset$ and $\cD_{3,2}[g(\lambda_0),y]\neq \emptyset$ and $\cG[g(\lambda_0),h(\lambda_0)]\setminus\cD_{3,2}[g(\lambda_0),y]=\emptyset$, & Yes\\
\hline
\end{tabular}
\label{tab:bifofRO11}
\end{table}

\begin{prop}
Bifurcation of a singular relaxation oscillation $\cA$ for $m=n=1$ due to codimension one bifurcation of the singular fast-slow system $\{f,g\}$ occurs for exactly one of the following reasons, at exactly one point $(x,y)\in \cA$. 
\begin{enumerate}
    \item There is saddle-node bifurcation of the slow subsystem in the interior of $\cC_{att}[g]$: $(x,y)\in \cE_{1,1}[g,h,y]$.
    \item There is a hyperbolic fold tangency of the critical manifold: $(x,y)\in \cD_{1,1}[g,y]$.
    \item There is a stable hysteresis of the critical manifold: $(x,y)\in \cD_{2,1}[g,y]$.
    \item There is an aligned double fold-umbra or umbra-umbra limit point of the critical manifold: $(x,y)\in \cD_{3,1}[g,y]$ or $(x,y)\in \cD_{3,3}[g,y]$.
    \item There is an opposed double fold-umbra or umbra-umbra limit point of the critical manifold: $(x,y)\in \cD_{3,2}[g,y]$ or $(x,y)\in \cD_{3,4}[g,y]$.
    \item A sink in the slow subsystem intersects a quadratic fold in the fast subsystem: $(x,y)\in \cM_{1,1}[g,h,y]$.
    \item A source in the slow subsystem intersects the umbra of a quadratic fold in the fast subsystem: $(x,y)\in \cM_{1,4}[g,h,y]$.
\end{enumerate}
\label{prop:srobifcodim1}
\end{prop}

\proof
We determine which of the codimension one degeneracies in Table~\ref{fig:degeneracytable11} and Table~\ref{tab:slowcodim1subcases11} can cause bifurcation of singular relaxation oscillations by ruling out those that cannot, and by providing examples in the following section.

First, two hyperbolic equilibria sharing y-coordinate $\cE_{2,1}[g,h,y]$, a non-interacting double limit point degeneracy $\cD_{3,5}[g,y]$ or $\cD_{3,6}[g,y]$ or an equilibrium sharing slow coordinate with a fold point but not intersecting it or its umbra $\cM_{1,5}[g,h,y]$ are ruled out because these degeneracies cannot break the simple property of limit cycles at codimension one.

Furthermore, codimension one bifurcation of limit cycles requires that a regular stable part of the critical manifold exists in a neighbourhood of the bifurcation point, or else the limit cycle does not generically pass through that point. This excludes elliptic fold tangency $\cD_{1,2}[g,y]$ and unstable hysteresis $\cD_{2,2}[g,y]$.

Moreover, a source equilibrium intersecting a fold $\cM_{1,2}[g,h,y]$ or a sink interacting with a fold umbra $\cM_{1,3}[g,h,y]$ are excluded since no relaxation oscillation can exist either at or in a neighbourhood of the bifurcation parameter at codimension one. Additionally, a ``fold umbra - fold umbra'' double limit point degeneracy cannot cause bifurcation at codimension one, since the umbra is generically on $\cC_{reg}[g]$, and not intersecting any equilibria (making the period finite). Therefore, there is no way for a simple relaxation oscillation to be lost at such a point.
\qed

The remaining codimension one bifurcations can break the simple property by violating one of the defining conditions: these are listed in Table \ref{tab:bifofRO11}. If a vector field is perturbed by a distinguished parameter $\lambda\in \bR$, a simple singular relaxation oscillation may cease to exist for some critical $\lambda_0$ where we assume the limit is from below. In such cases we define a limit relaxation oscillation $\cA_{lim}(\lambda_0)$ as the limit, in the Hausdorff distance, of a sequence of relaxation oscillations $\cA(\lambda)$ parametrized by $\lambda$
\begin{equation}
    \cA_{lim}(\lambda_0) = \lim_{\lambda \to \lambda_0-} \cA(\lambda).
\label{eq:limitro}
\end{equation}
The limit is well defined for simple relaxation oscillations but may be empty. The bifurcation of relaxation oscillations will unfold for the non-singular systems $\epsilon>0$ to give a variety of canards that will appear on a case to case basis: see for example \cite[Chapter 8]{kuehn15} and \cite{maesschalck11}. Outside a small (in $\epsilon$) range of parameters $\lambda(\epsilon)$ near the critical $\lambda_c$, many of the solutions will closely resemble those of the singular system. Hence, the bifurcations in Proposition \ref{prop:srobifcodim1} will give rise to a detectable qualitative change even for non-singular systems.

\subsection{Examples of bifurcations of relaxation oscillations}
\label{sec:examplebifs}

\begin{figure}
\begin{center}
\includegraphics[width=12.8cm]{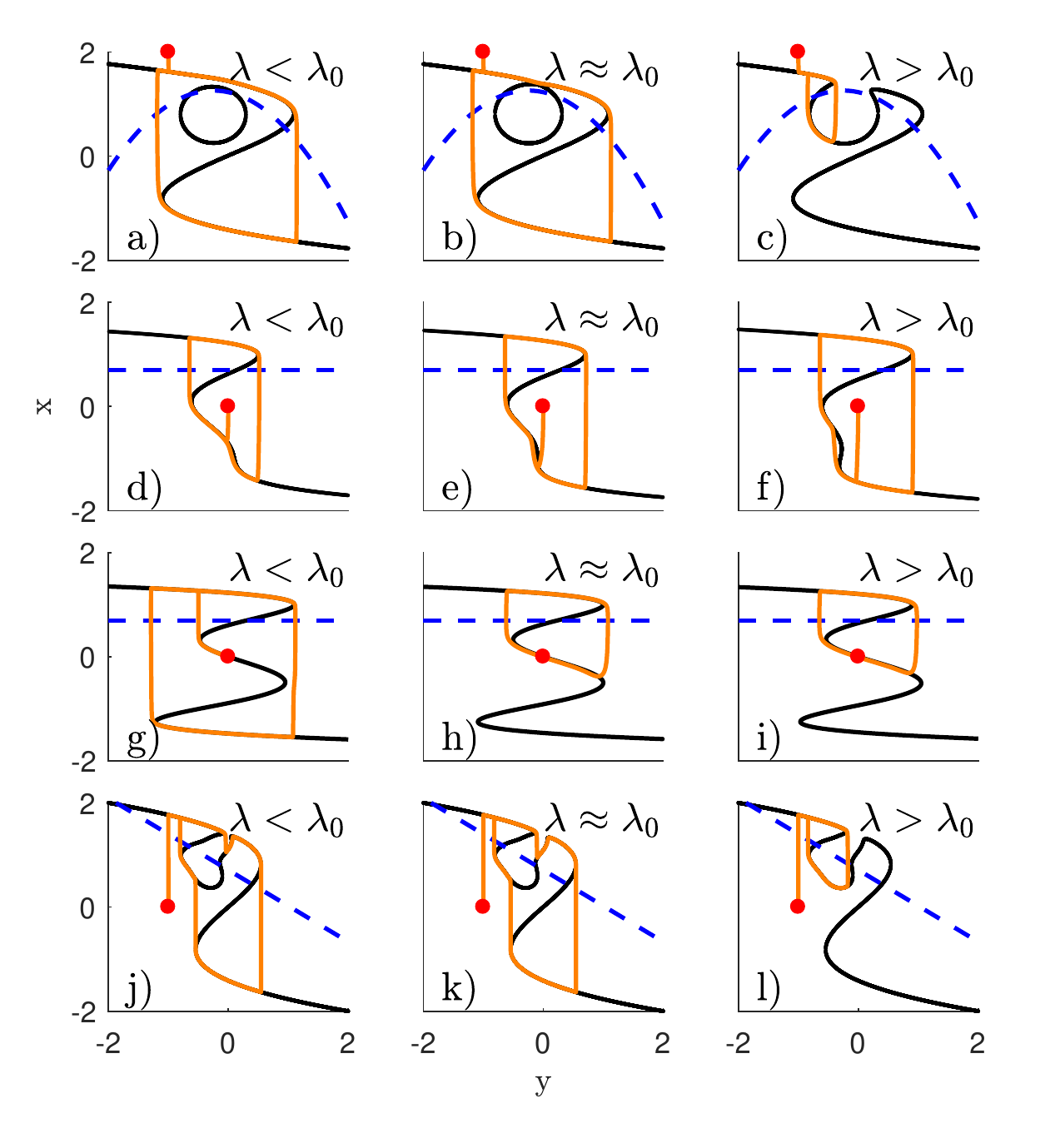}
\end{center}
\caption{(Color online) Examples of bifurcation of relaxation oscillation due to bifurcation of the critical set for $m=n=1$. Bifurcation occurs at the critical value $\lambda_0$ of the bifurcation parameter $\lambda$. Black lines are the critical set, blue dashed lines are nullclines of the slow subsystem, and orange lines show example trajectories of a nearly singular system started at the indicated red point, that evolve towards a relaxation oscillation. In panels j,k,l) the point of interest is the left facing fold near $y=0$ which nearly intersects the slow nullcline; at bifurcation that fold shares $y$-value with another right facing fold at a larger value of $x$. As $\lambda$ increases the left facing fold moves rightward and singular relaxation oscillations change from passing right of the fold to the left of it. Polynomial equations for the critical set are listed in Table \ref{tab:critmanbifpoly11}, bifurcation conditions are listed in Table \ref{tab:bifofRO11}, and a detailed description of the systems is found in \ref{sec:examplesection}}
\label{fig:slowfastbifexamples}
\end{figure}

To illustrate how the bifurcations of limit cycles can be realised, we show in Figure~\ref{fig:slowfastbifexamples} some examples of bifurcations of relaxation oscillations in fast-slow systems (\ref{eq:mainsystem}) with $m=n=1$, near the singular limit. In all cases the critical manifold is expressed as a relatively low-order polynomial (Table \ref{tab:critmanbifpoly11}). We explain how these were derived in  \ref{sec:examplesection}.

\begin{table}
\caption{Examples of fast subsystems $g$ that undergo each of the codimension one bifurcations of the critical set $g(x,y)=0$ for $m=n=1$ at $\lambda=\lambda_c\approx 0$, shown in Figure~\ref{fig:slowfastbifexamples}. Indefinite integrals are taken to have zero constant term. $(\hat{x},\hat{y})$ are scaled, rotated and translated coordinates. Details  how $g(x,y)$ is constructed, and parameters for the opposed double limit degenerate case can be found in  \ref{sec:examplesection}
}
\centering
\begin{tabular}{|p{0.2\linewidth}|p{0.70\linewidth}|}
\hline
Fold tangency \newline (Fig.~\ref{fig:slowfastbifexamples} a,b,c)) & $g(x,y)=-(g_1(x,y)g_2(x,y) + \lambda x +q),\mbox{ with } g_1(x,y)=x^3 - 2x + y \mbox{ and } g_2(x,y)=(x-x_c)^2 + (y-y_c)^2 - R^2 \mbox{ and } (x_c,y_c,R,q)=(81/100,-1/4,11/20,1/100)$,\\
\hline
Hysteresis  \newline (Fig.~\ref{fig:slowfastbifexamples} d,e,f) & $g(x,y)=\int -a(x+x_1)(x+x_2)(x+x_3)^2dx + \lambda x - b - y, (a,b,x_1,x_2,x_3) = (15/4,6/10,-1,1/25,-1)$,\\
\hline
Aligned double limit \newline (Fig.~\ref{fig:slowfastbifexamples} g,h,i) & $g(x,y)=\int -a(x+x_1)(x+x_2)(x+x_3)(x+x_4)dx + \lambda x - y, (a,x_1,x_2,x_3,x_4) = (640/49,-1,-13/40,1/2,5/4)$,\\
\hline
Opposed double limit \newline (Fig.~\ref{fig:slowfastbifexamples} j,k,l) & $g(x,y)=-(g_1(x,y)g_2(\hat{x},\hat{y})+\lambda x + q),\mbox{ with } g_1(x,y) = 0.5x^3-x + y,~\mbox{ and } g_2(\hat{x},\hat{y}) = (\hat{x}^2 + \hat{y}^2)^3 - (\hat{x}^2 + (\hat{x}^2 + \hat{y}^2)^2\hat{y}^2)$,\\
\hline
\end{tabular}
\label{tab:critmanbifpoly11}
\end{table}

\FloatBarrier

\section{Discussion}
\label{sec:discussion}

Almost two decades ago, Guckenheimer \cite{guckenheimer02} called for a classification of bifurcations of relaxation oscillations in fast-slow systems up to two slow and two fast variables.  In this paper we have used bifurcation theory with distinguished parameters and singular equivalence to take some steps towards such a classification. Indeed, in \cite{guckenheimer02} Guckenheimer gives the following list of codimension one degeneracies that we can relate to our classification:
\begin{itemize}
\item[\bf G1:] A fast segment ends at a regular fold point. There are two cases depending on whether the slow flow approaches or leaves the fold near this point.
\item[\bf G2:] A slow segment ends at a folded saddle.
\item[\bf G3:] A fast segment encounters a saddle point.
\item[\bf G4:] There is a point of Hopf bifurcation at a fold.
\item[\bf G5:] A slow segment ends at a cusp.
\item[\bf G6:] The reduced system has a quadratic umbral tangency between projections of fold and umbra.
\end{itemize}
The degeneracy {\bf G1} is a subcase of fold projection intersection for $m=1,n=1,2$ and fold  projection tangency for $m=1,n=2$. The degeneracy {\bf G2} appears when the slow flow is tangent to a fold line: this can occur for $n\geq2$. Degeneracy {\bf G3} can appear at a saddle for $m\geq2$ or at an unstable node for $m=1$. Note that the formulation of {\bf G3} is slightly modified from \cite{guckenheimer02}. Degeneracy {\bf G4} corresponds to a singular Hopf bifurcation, which we discussed in the context of $m=n=1$. Degeneracy {\bf G5} corresponds to a hysteresis bifurcation for $m=n=1$, and to a limit cycle hitting a cusp on the slow manifold for $m=1$, $n=2$. Finally, degeneracy {\bf G6} requires $n\geq2$.

Guckenheimer states in \cite{guckenheimer02} that the list is incomplete, and mentions the case that a slow segment ends at a folded node as an example. Degeneracy due to fold tangency of the critical set is missing from the list, since the slow variables were not regarded as distinguished parameters in \cite{guckenheimer02}. We believe that Proposition \ref{prop:srobifcodim1} completes the list for $m=n=1$.

For $n=2$, the degeneracies involve tangencies of generic one-dimensional objects such as relaxation oscillations, fold lines and fold umbrae, as well as intersections of one-dimensional objects and generic zero-dimensional objects such as cusp points and equilibria in the slow subsystem. Some of these cases are in Guckenheimer's list. Note that  degeneracies of the critical manifold do not cause codimension one bifurcation of relaxation oscillations since they occur at points, which do not generically intersect relaxation oscillations. They will be involved in bifurcation of invariant tori or more complex singular attractors or of relaxation oscillations at higher codimension however.

\subsection{Relation to other singularity theory problems}

There appears to be connections between global equivalence of critical sets with two distinguished parameters (for $m=1$ and $n=2$) and the equivalence of vector fields under projection to the slow plane. In particular, several hypothesised degeneracies under strong equivalence also appear as degeneracies of orthogonal projections of vector fields, see for instance \cite{alharbi14,bruce84,ohmoto06,thomas11,yoshida14} and references. Crucially, however, such equivalences do not produce degeneracies such as the fold tangency where the manifold structure is lost.

Our approach to bifurcation of the critical manifold in fast-slow systems uses a singularity theory approach with distinguished parameters from \cite{golubitskyscheffer85}. In the following we briefly indicate how this approach is related to singularity theory, catastrophe theory, the theory of constrained equations, and projections from manifolds to manifolds.

Much of singularity theory concerns the stability of zero sets of smooth functions $g(x)=0$ under perturbation \cite{whitney55}. In the singularity theory approach to bifurcation theory of \cite{golubitskyscheffer85}, there is a distinguished (bifurcation) parameter $y$ that is not ``mixed up'' with the unfolding parameters. Hence, e.g. the quadratic fold $g=x^2 + y$ is a codimension one degeneracy of $g(x)=x^2$ in singularity theory, but is codimension zero in bifurcation theory with one distinguished parameter $y$. For our interpretation, $y$ is identified with the slow variables.

Catastrophe theory \cite{arnold99,PostonStewart1978} classifies the changes to stationary points of potentials $V(x)$, with $\nabla V(x)=g(x)$, by codimension of  deformation. Although different equivalences and objects are studied in catastrophe theory and singularity theory, for one (fast) variable, the classification of local singularities of the critical set is the same. Constrained systems \cite{takens76,jardon14} correspond to singular fast-slow systems where equivalence of critical manifolds is defined by potential functions, as in catastrophe theory, together with a slow flow local to a point. Unlike our approach, there are no distinguished parameters and the unfolding parameters are identified with the slow variables. This means that some local bifurcations (notably fold tangencies) that are present for the distinguished parameter approach are missed, because ``slow variables" never appear in powers higher than one in the local normal forms.

Singularity theory, catastrophe theory and constrained equations have been framed in terms of germs, which are local notions of functions. This means that global intersections of projections of singularities (such as double limit points which are important for bifurcations of relaxation oscillations) have not been widely studied in these contexts, some exceptions being \cite{arnold99,broer13,guckenheimer02,maesschalck11}.

\subsection{Further perspectives}

Persistence and codimension one bifurcation of the critical set for one fast ($m=1$) and two slow variables ($n=2$) remains to be proved. This requires a suitable equivalence, which should give rise to the degeneracies that we have listed, but possibly more.

A full investigation of bifurcations of singular relaxation oscillations for $m=1,n=2$ is outside the scope of this paper. Some specific examples have been studied by Guckenheimer \cite{guckenheimer03a,guckenheimer12,guckenheimer04} who outlined a scheme for investigation of bifurcation of solutions to singular fast-slow systems in \cite{guckenheimer96}.

The general case of two fast variables $m=2$ is considerably more complicated as the vector fields cannot be written as gradients of potentials, and hence there can be other asymptotic behaviour than fixed points. If $n=2$ then for generic asymptotic fast dynamics, the system will approach a critical set that is a union of all equilibria, periodic orbits and homoclinic/heteroclinic cycles of the fast system. The persistence of bifurcations on the critical set will depend on the number of slow variables. For $n=1$ then we expect persistence precisely when (a) All singularities of equilibria within the critical set are quadratic folds or Hopf points. (b) All singularities of limit cycles within the critical set are one of saddle-nodes of limit cycles, saddle node on a periodic orbit, or homoclinic bifurcation.
(c) The slow flow has generic intersection with umbrae of the singularities. For $n=2$ we will get in addition generic local and global codimension two singularities at isolated points in the slow variables; this will include, for example, cusp points, Bogdanov Takens points and Bautin points at singular equilibria, and a wide variety of possible generic codimension two bifurcations of homoclinic orbits \cite{ChampneysKuznetsov1994}.

We have ignored the phenomena that arise when the scale separation is imperfect, that is for $\epsilon>0$. In that case, the fast and slow subsystems evolve at similar speeds close to singular points; this gives rise to \emph{canards} and mixed mode oscillations \cite{desroches14}. Canard behaviour has been extensively studied, especially near regular values of the critical set see e.g. \cite{benoit81,kuehn15,guckenheimer96,wechselberger12}. Canards for degenerate critical sets are discussed in \cite{arnold99,broer13,maesschalck11}, but we are unaware of any systematic treatment.

\subsection*{Acknowledgements}

We thank James Montaldi, Christian Kuehn, Hildeberto Jardón-Kojakhmetov, Bernd Krauskopf, Christian Bick and Daniele Proverbio for valuable discussions. We also thank two anonymous reviewers and two editors whose comments significantly improved the manuscript. This research has been funded by the European Union's Horizon 2020 innovation and research programme for the ITN CRITICS under the Marie Sk\l{}odow\-ska-Curie grant agreement No. 643073.

\subsection*{References}

\bibliographystyle{plain}

\DeclareRobustCommand{\DE}[3]{#3}
\DeclareRobustCommand{\VANDER}[3]{#3}

\bibliography{singbifs_references_rev}

\appendix

\section{The quadratic curvature}
\label{app:curvature}

We define the scalar quadratic fold curvature (see Table~\ref{tab:nonpersistentdegeneracies12})) at a fold point $p$ in the direction of a fold of the critical manifold as
\begin{equation}
\begin{array}{rl}
    K[g](p) = & \sign{(g_{xx}(p))}\frac{\overline{{\nabla_y^\perp g(p)}}^T D^2_y(g(p)){\overline{\nabla_y^\perp g(p)}}}{2|\nabla_y^\perp g(p)|} - \frac{\left(\nabla_y g_x(p) \cdot \overline{\nabla_y^\perp g(p)}\right)^2}{8|g_{xx}(p)||\nabla_y^\perp g(p)|},
\label{eq:curvaturecriterion}
\end{array}
\end{equation}
where $D_y^2(g)$ is the Hessian of the slow subsystem, superscript $T$ denotes transpose, $\perp$ denotes perpendicular, and $\overline{{\nabla_y^\perp g}}$ denotes ${\nabla_y^\perp g}$ scaled to unit length. For the remainder of this section we drop the dependency on $p$, such that e.g. $g(p)$ is written just $g$. If $\nabla_y^\perp g = 0$, then $K[g]$ is undefined. Note that this notion of scalar quadratic fold curvature only captures the quadratic curvature of the fold curve as a projection onto the slow $(0,y_1,y_2)$ plane.

In this section, we first motivate our definition of (\ref{eq:curvaturecriterion}) and then offer an interpretation.

To motivate (\ref{eq:curvaturecriterion}), we start with a quadratic fold at the origin and completely in the slow plane whose curvature in the slow plane reasonably should be considered to be along the $y_1$-axis
\begin{equation*}
g(x,y_1,y_2) = \xi x^2 + c_1 y_1 + c_2 y_2^2,
\end{equation*}
where $\xi$, $c_1$ and $c_2$ are real constants.

Recall that the quadratic direction vector of a fold is given by $\nu[g] = g_{xx}\nabla_y g=(2\xi c_1,0)$, so if $\xi c_1>0$ then the fold is directed rightward, while if $\xi c_1 < 0$ the direction is directed leftward. It makes sense to define the slow quadratic fold curvature of this fold as
\begin{equation}
K[g] = \sign{(\xi)}\frac{c_2}{2c_1},
\label{eq:normalformcurvature}
\end{equation}
since then the curvature is independent of the magnitude of $\xi$, proportional to $c_2=0$, and inversely proportional to $c_1$. The fold is convex in the direction of the fold if $K[g](p) > 0$ and concave if $K[g](p) < 0$; see Figure \ref{fig:geometricdefinitions} for a graphical representation.

We now consider a general quadratic polynomial function of a quadratic fold at the origin with terms of relevant order
\begin{equation}
g(x,y_1,y_2) = \xi x^2 + ay_1 + by_2 + \alpha y_1^2 + \beta y_2^2 + 2\gamma y_1y_2 + \delta x y_1 + \eta x y_2,
\label{eq:generalfold}
\end{equation}
where $a$, $b$, $\alpha$, $\beta$, $\gamma$, $\delta$ and $\eta$ are real constants unrelated to any quantities with the same names elsewhere in this text. The $x$ term and constant term are missing because of the quadratic fold condition $g=g_x=0$, and higher order terms are not present since they should not matter for the {\em quadratic} curvature.

Next, we find an expression for the projection of the quadratic approximation of the fold curve onto the slow plane $(0,y_1,y_2)$. To this end, we solve $(\ref{eq:generalfold})=0$ for $x$ to get
\begin{equation}
\begin{array}{r l}
x = & \frac{1}{2 \xi}\bigg((\delta + \eta) \\
& \pm \sqrt{(\delta y_1 + \eta y_2)^2 - 4 \xi(\alpha  y_1^2 + \beta y_2^2 + 2\gamma y_1y_2 + a y_1 + by_2)} \bigg).
\label{eq:quadraticslowprojection}
\end{array}
\end{equation}
At the quadratic approximation of the fold curve the discriminant of (\ref{eq:quadraticslowprojection}) is zero
\begin{equation*}
(\delta y_1 + \eta y_2)^2 - 4 \xi(\alpha  y_1^2 + \beta y_2^2 + 2\gamma y_1y_2 + a y_1 + by_2)=0,
\end{equation*}
giving a condition on $y_1$ and $y_2$. Expanding parentheses, collecting terms and dividing by $-4\xi$ gives that
\begin{equation*}
a y_1 + b y_2 + \left(\alpha - \frac{\delta^2}{4\xi}\right)y_1^2 +\left (\beta - \frac{\eta^2}{4\xi}\right) y_2^2 + 2\left(\gamma - \frac{\delta \eta}{4\xi}\right) y_1 y_2=0.
\end{equation*}
We define the new coefficients $\tilde{\alpha}=\alpha - \delta^2/4\xi$, $\tilde{\beta} = \beta - \eta^2/4\xi$ and $\tilde{\gamma} =\gamma - \delta \eta/4\xi$, such that
\begin{equation}
a y_1 + b y_2 + \tilde{\alpha}y_1^2 +\tilde{\beta} y_2^2 + 2\tilde{\gamma} y_1 y_2=0.
\label{eq:quadraticfoldconditionnomixedx}
\end{equation}
Based on (\ref{eq:quadraticfoldconditionnomixedx}) we define a new function 
\begin{equation}
\tilde{g}(x,y_1,y_2)= \xi x^2 + a y_1 + b y_2 + \tilde{\alpha}y_1^2 +\tilde{\beta} y_2^2 + 2\tilde{\gamma} y_1 y_2,
\label{eq:quadraticfoldfucntionnomixedx}
\end{equation}
which defines a quadratic fold curve at the origin having the same slow quadratic curvature as $g$ in (\ref{eq:generalfold}) but lying entirely in the slow plane.

Next, we seek a rotation $R$ of the slow variables which brings the quadratic fold direction vector $\nu[g]=\sign{(g_{xx})}(a,b)^T$ of $g$ (as well as $\tilde{g}$) in the positive $y_1$ direction, that is: $R \nu[g] = (|\nu[g]|,0)^T$. This rotation does not change the slow curvature, but it allows us to identify the relevant coefficients corresponding to $c_1$ and $c_2$ in (\ref{eq:normalformcurvature}).

The sought rotation in matrix form is
\begin{equation*}
R = \sign{(g_{xx})}\frac{1}{\sqrt{a^2 + b^2}}\left(\begin{array}{cc} a & b \\
-b & a
\end{array}\right),\end{equation*}
and consequently
\begin{equation*}R^{-1} = \sign{(g_{xx})}\frac{1}{\sqrt{a^2 + b^2}}\left(\begin{array}{cc} a & -b \\
b & a
\end{array}\right).
\end{equation*}
The old slow coordinates $(y_1,y_2)$ are expressed in the new ones $(\hat{y_1},\hat{y_1})$ as
\begin{equation*}
\left(\begin{array}{cc}
y_1 \\ y_2 
\end{array}\right)
= 
R^{-1}
\left(\begin{array}{cc}
\hat{y_1} \\ \hat{y_2} 
\end{array}\right).
\end{equation*}
In the new coordinates (\ref{eq:quadraticfoldfucntionnomixedx}) becomes
\begin{equation*}
\begin{array}{rl}
g(x,\hat{y_1},\hat{y_2}) = & \xi x^2 + \frac{1}{a^2 + b^2} [ \sign{(g_{xx})}\sqrt{a^2 + b^2}a(a\hat{y_1} - b\hat{y_2}) \\
& + \sign{(g_{xx})}\sqrt{a^2 + b^2}b(b\hat{y_1} + a\hat{y_2}) \\
& + \tilde{\alpha} (a\hat{y_1} - b\hat{y_2})^2 + \tilde{\beta} (b\hat{y_1} + a\hat{y_2})^2 \\
& + 2\tilde{\gamma} (a\hat{y_1} - b\hat{y_2})(b\hat{y_1} + a\hat{y_2}) ] \\
= & \xi x^2 + \frac{1}{a^2 + b^2}[\sign{(g_{xx})}(a^2 + b^2)^{3/2}\hat{y_1} \\
& + (\tilde{\alpha} b^2 + \tilde{\beta} a^2 + 2\tilde{\gamma} ab)\hat{y_1}^2 \\
& + (\tilde{\alpha} b^2 + \tilde{\beta} a^2 - 2\tilde{\gamma} ab)\hat{y_2}^2 \\
& + 2(-\tilde{\alpha} ab + \tilde{\beta} ab + (a^2 - b^2)\tilde{\gamma})\hat{y_1}\hat{y_2}].
\end{array}
\end{equation*}
Hence, reading off the coefficients of $\hat{y}_1$ and $\hat{y}_2^2$ in analogy with (\ref{eq:normalformcurvature}), we get the scalar quadratic fold curvature
\begin{equation}
\begin{array}{rl}
    K[g] = & \sign{(g_{xx})}\frac{(\tilde{\alpha} b^2 + \tilde{\beta} a^2 - 2\tilde{\gamma} ab)}{(a^2+b^2)^{3/2}} \\
    = & \sign{(g_{xx})}\frac{(\alpha b^2 + \beta a^2 - 2\gamma ab)}{(a^2+b^2)^{3/2}} + \frac{(\delta^2 b^2 + \eta^2 a^2 - 2\delta\eta ab)}{4|\xi|(a^2+b^2)^{3/2}} \\
    = & \sign{(g_{xx})}\frac{(g_{y_1y_1}g_{y_2}^2 - g_{y_2y_2}g_{y_1}^2 - 2g_{y_1y_2}g_{y_1}g_{y_2})}{2(g_{y_1}^2+g_{y_2}^2)^{3/2}} \\ 
    & - \frac{(g_{xy_1}^2g_{y_2}^2 + g_{xy_2}^2g_{y_1}^2 - 2g_{xy_1}g_{xy_2}g_{y_1}g_{y_2})}{8|g_{xx}|(g_{y_1}^2+g_{y_2}^2)^{3/2}} \\
    = & \sign{(g_{xx})}\frac{\overline{{\nabla_y^\perp g}}^T D^2_y(g){\overline{\nabla_y^\perp g}}}{2|\nabla_y^\perp g|} - \frac{\left(\nabla_y g_x \cdot \overline{\nabla_y^\perp g}\right)^2}{8|g_{xx}||\nabla_y^\perp g|},
\end{array}
\label{eq:scalarquadraticcurvaturederived}
\end{equation}
where $ \overline{\nabla_y^\perp g} = \nabla_y^\perp g/|\nabla_y^\perp g|$ and $\nabla_y g_x = (g_{x y_1},g_{x y_2})$. The case $K[g]>0$ implies a locally convex fold and $K[g]<0$ implies a locally concave fold. The degenerate case $K=0$ implies that the fold line is locally straight, or not quadratic.

We can now define the quadratic fold curvature vector as
\begin{equation*}
\begin{array}{rl}
\kappa[g](p) = & K[g](p)\overline{\nu[g] (p)}\\
= & \left( \frac{\overline{{\nabla_y^\perp g(p)}}^T D^2_y(g(p)){\overline{\nabla_y^\perp g(p)}}}{2|\nabla_y^\perp g(p)|} - \frac{\left(\nabla_y g_x(p) \cdot \overline{\nabla_y^\perp g(p)}\right)^2}{8 g_{xx}(p) |\nabla_y^\perp g(p)|}\right) \overline{\nabla_y g(p)}.
\end{array}
\end{equation*}
The quadratic fold curvature vector points in the direction of the fold if the fold is convex, and against the direction if it is concave (see Figure \ref{fig:geometricdefinitions}).

Some special cases of Equation (\ref{eq:curvaturecriterion}) are insightful. First, if the fold lies locally in the slow plane (that is, if the coefficients of the $x y_1$ and $x y_2$ terms in (\ref{eq:generalfold}) are zero), then a positive definite matrix $\sign{(g_{xx})}D^2_y(g)$ implies that the fold is convex and a negative definite $\sign{(g_{xx})}D^2_y(g)$ implies that it is concave. However, if $\sign{(g_{xx})}D^2_y(g)$ is indefinite or has a zero eigenvalue, then the sign of $K[g]$ can be either positive, negative or zero depending on the direction of the fold.

On the other hand, if the $x y_1$ and $x y_2$ terms are present in (\ref{eq:generalfold}), then the second term in $(\ref{eq:curvaturecriterion})$ only serves to reduce convexity (or equivalently increase concavity). The extent to which convexity is reduced depends quadratically on the component of the gradient of $\nabla_y g_x$ perpendicular to the gradient, and inversely on the magnitude of the gradient and the magnitude of the curvature in the $x$ direction.

\subsection{Persistent subcases of the fold projection tangency}

We classify the fold projection tangency degeneracy (degeneracy subset $\cD_4[g]$ for $m=1$ fast and $n=2$ slow variables) into a number of qualitatively different subcases. The subcases are separated by the scalar quadratic fold curvatures at the points of degeneracy, whether the umbrae interact with each other or a fold, and in the case of fold-umbra degeneracy, whether the dominant curvature belongs to the fold curve with the largest $x$-coordinate.

Aligned folds generate four distinct subcases. In analogy with the situation for $m=n=1$, we have fold-fold, fold-umbra, and non-interacting fold cases. But the fold-umbra has two subcases, depending on whether the fold with interacting umbra also has the greatest scalar quadratic fold curvature $K[g]$. Hence, there are four subcases.

Opposed folds have six distinct subcases, three for each case that either the sum of curvatures is positive (net convex) or negative (net concave). Opposed fold projection tangency does not have the two  fold-umbra subcases of aligned fold projection tangency, since the dominant $x$-values are reversed by one half rotation of the slow variables. Hence there are six such subcases, and ten subcases in total.

The nonpersistence condition (at codimension one) for fold projection tangency degeneracy at points $p_1\in \cD_4[g]$ and $p_2 \in \Pi(p_1) \cap \cD_4[g]$ is
\begin{equation*}
\kappa[g](p_1) \neq \kappa[g](p_2)
\end{equation*}
that is, the folds must have distinct quadratic curvature vectors (see \ref{app:curvature}). However, to distinguish subcases of degeneracy we will later use the scalar quadratic fold curvature $K[g]$ and information about whether folds are aligned or opposed.

\begin{figure}
\begin{center}
\includegraphics[width=15cm]{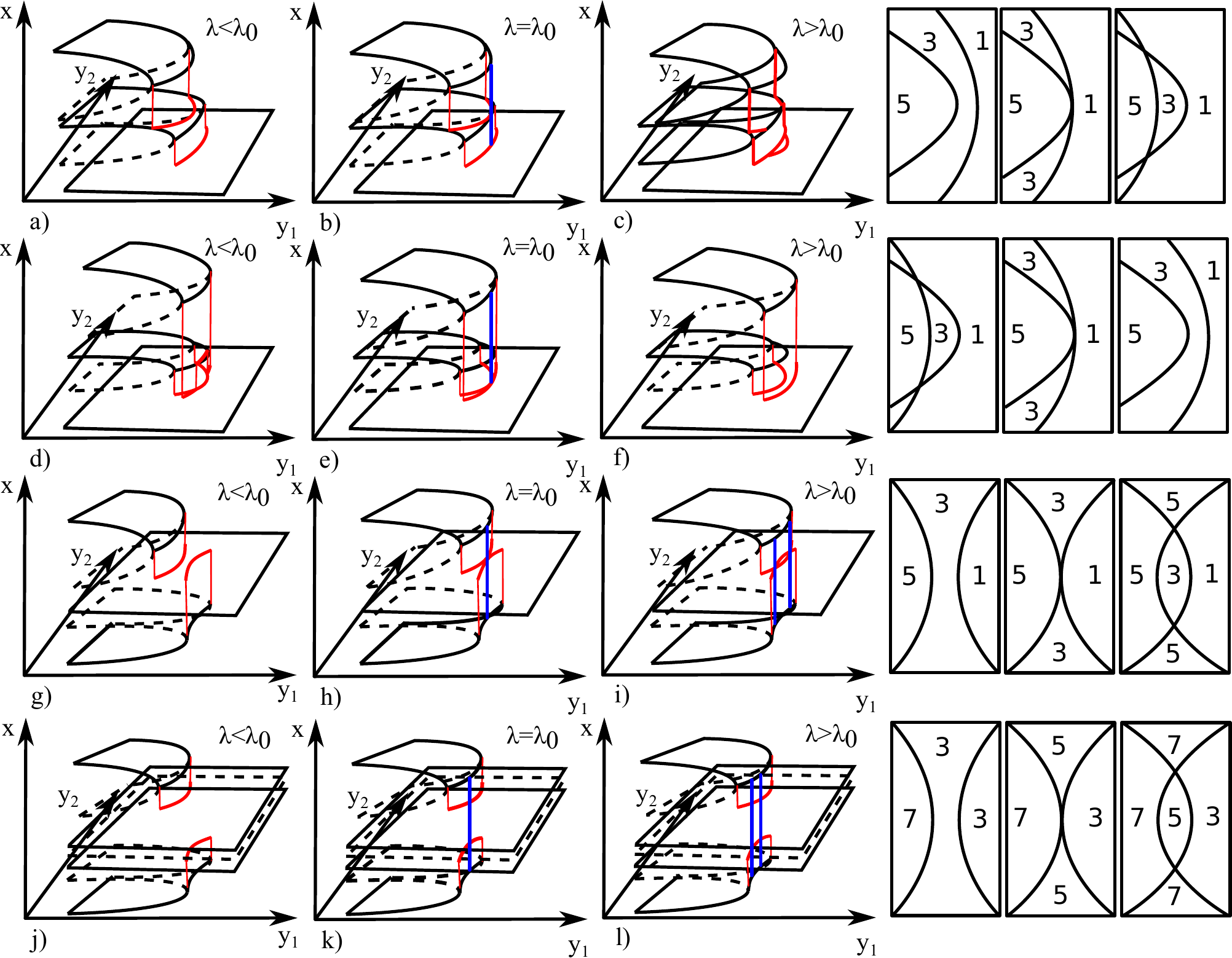}
\end{center}
\caption{
(Color online) Examples of codimension one aligned fold projection tangency bifurcation (see Table \ref{tab:globaldegeneracies12P2Folds}). Each row shows unfolding with a bifucation parameter $\lambda$. Bifurcation occurs as $\lambda=\lambda_0$. Solid/dashed black lines show the stable/unstable sheets of the critical set, red lines show the image of the fold under the umbral map and blue lines indicate tangency of projections of fold sheets. As a visual aid, the number of sheets of the critical set in a neighbourhood of the bifurcation is shown to the right, to be viewed as a projection onto the slow variables
\label{fig:structureCM121D4Aligned}
}
\end{figure}

We now separate subcases depending on whether folds are aligned or opposed. If $\nu[g](p_1) \cdot  \nu[g](p_2) > 0$ then the folds are aligned, and qualitatively indistinguishable unless the umbra of the fold with the larger $x$ component hits the other fold. If the umbra $U[g](p_1)$ of an aligned fold at a point $p_1$ hits another fold at a point $p_2$ and the scalar quadratic fold curvatures satisfy $K[g](p_1)>K[g](p_2)$, then the degeneracy is called \emph{umbra dominant}. Otherwise if $K[g](p_1)<K[g](p_2)$, then the degeneracy is called \emph{fold dominant}. Therefore, we get the four subcases illustrated in Figure~\ref{fig:structureCM121D4Aligned} and tabulated in Table~\ref{tab:globaldegeneracies12P2Folds}. (Figure~\ref{fig:foldtangencysketch} a,b,c) shows a slow projection sketch of aligned fold projection tangency).

If, on the other hand $\nu[g](p_1) \cdot \nu[g](p_2) < 0$, then the folds are opposed and therefore distinguishable at bifurcation. (Figure~\ref{fig:foldtangencysketch} d,e,f) shows a slow projection sketch of opposed fold projection tangency).

Assuming that folds are opposed, the defining condition for covering opposed fold projection tangency is 
\begin{equation*}
K[g](p_1) + K[g](p_2) < 0,
\end{equation*}
that is, if at least one fold is concave, and the magnitude of the curvature of the convex fold is smaller than that of the convex curve (see Figure~\ref{fig:foldtangencysketch} d,e,f)).
If, on the other hand
\begin{equation*}
K[g](p_1) + K[g](p_2) > 0,
\end{equation*}
then we have covering opposed fold projection tangency (see  Figure~\ref{fig:foldtangencysketch} g,h,i) ).

Opposed fold curves do not have the umbra-dominant and fold-dominant subcases that aligned fold curves do, since a rotation of the slow variables by 180 degrees turns one such case into the other. Therefore, there are four aligned cases and six opposed cases, giving in total the ten subcases in Table~\ref{tab:globaldegeneracies12P2Folds}.

\section{Definition of the cusp direction vector}
\label{app:cuspdirection}

We define the direction of a cubic cusp (see Table~\ref{tab:nonpersistentdegeneracies12})) to be
\begin{equation}
\mu[g] = \frac{g_{xxx}(p)}{\nabla_y g_x(p) \cdot \overline{\nabla_y^\perp g(p)}}\nabla_y^\perp g(p).
\end{equation}
For the remainder of this section we do not explicitly write out the dependence on $p$, so that for instance $g(p)$ is written $g$. For the example in this section we assume that $p=(0,0,0)$. As for the definition of scalar quadratic fold curvature, we motivate the definition starting from a normal form of the cubic cusp
\begin{equation}
g(x,y_1,y_1) = \xi x^3 +c_1 xy_2 + c_2 y_1,
\end{equation}
where $\xi$, $c_1$ and $c_2$ are real constants. In this case we naturally let the cusp direction vector be
\begin{equation}
    \mu[g] = \left(\frac{\xi}{c_1}\right)(0,-c_2),
\end{equation}
where $(0,-c_2)$ is a vector perpendicular to the gradient $(c_2,0)$. It makes sense for the direction of the cusp to be along the gradient perpendicular, since it is parallel to the pair of fold curves which emanate from the cusp.

Next, we consider a more general expression
\begin{equation}
g(x,y_1,y_2) = \xi x^3 + \alpha x y_1 + d x y_2 + ay_1 + by_2,
\label{eq:generalcuspnormalform}
\end{equation}
where $\xi, a,b,\alpha$ and $\beta$ are coefficients, and where we do not keep terms of intermediate orders $xy_1^2$ or $x^2y_2$, since visually they do not seem to alter the direction or sharpness of the cusp. We then rotate the slow subsystem as to make the gradient in the old coordinates $(a,b)$ directed along the positive $y_1$ axis. We accomplish this with a rotation $R$
\begin{equation*}
R = \sqrt{a^2 + b^2}
\left(\begin{array}{cc} a & b \\
-b & a
\end{array}\right),
\end{equation*}
mapping the new coordinates $\hat{y_1},\hat{y_2}$ to the old ones
\begin{equation*}
\left(\begin{array}{c} y_1 \\
y_2
\end{array}\right) = R^{-1}
\left(\begin{array}{c} \hat{y_1} \\
\hat{y_2}
\end{array}\right),
\end{equation*}
such that in the new coordinates Equation (\ref{eq:generalcuspnormalform}) becomes
\begin{equation}
\begin{array}{rl}
g(x,\hat{y_1},\hat{y_2}) = & \xi x^3 + \frac{1}{\sqrt{a^2+b^2}}[\alpha(a\hat{y_1} -b\hat{y_2})x + \beta(b\hat{y_1} + a\hat{y_2})x \\
& + a(a\hat{y_1} -b\hat{y_2}) + b(b\hat{y_1} + a\hat{y_2})] \\
= & \xi x^3 + \frac{(a\alpha + \beta b)}{\sqrt{a^2+b^2}}\hat{y_1}x + \frac{(-b\alpha + a \beta)}{\sqrt{a^2+b^2}}\hat{y_2}x + \sqrt{a^2+  b^2}\hat{y_1}.
\end{array}
\end{equation}
By reading off the coefficients of the $x^3$, $y_2x$ and $y_1$ terms, we find that the cubic cusp direction vector in the new coordinates is
\begin{equation}
\mu[g] = \frac{\xi \sqrt{a^2 + b^2}}{-bc + ad}(\sqrt{a^2 + b^2},0).
\end{equation}
Rotating this vector back to the original coordinates we get that
\begin{equation}
\mu[g] = \frac{\xi \sqrt{a^2 + b^2}}{-bc + ad}(-b,a),
\end{equation}
which we recognise can be written (up to constant scaling)
\begin{equation}
\mu[g] = \frac{g_{xxx}}{\nabla_y g_x \cdot \overline{\nabla_y^\perp g}},
\end{equation}
where $\overline{\nabla_y^\perp g}=\nabla_y^\perp g/|\nabla_y^\perp g|$. Hence, the cusp is always perpendicular to the slow gradient, with magnitude inversely proportional to the projection of $\nabla_y g_x$ onto the gradient perpendicular. As a consequence, the magnitude of the cubic cusp direction vector blows up (becomes undefined) if $\nabla_y g_x$ is parallel to the gradient (its perpendicular component vanishes).

\section{The quantity \texorpdfstring{$W[g]$}{W[g]}}
\label{app:W}

To construct the cusp quantity $W[g](p)$ (Table~\ref{tab:nonpersistentdegeneracies12}) at a point $p=(x,y_1,y_2)$, we start from a modified normal form for the cusp tangency bifurcation (see Table \ref{tab:bifnormalforms})
\begin{equation*}
g(x,y_1,y_2) = \delta_1x^3 + \delta_2 x y_1^2 + \delta_3 x^2y_1 + x \lambda + y_2,
\end{equation*}
where $\lambda$ is an unfolding parameter. The sign of the cusp quantity $W[g] = 6\delta_1\delta_2 - 2\delta_3\neq0$ gives the type of bifurcation ($W[g]<0$ gives beaks and $W[g]$ gives lips), see \cite{peters91}. Our aim is to express $W[g]$ for a general function equivalent to $g$. In this section we do not explicitly write out the dependence on $p$, which for the normal form is $p=(0,0,0)$, so that e.g. $W[g](p)$ is written $W[g]$. The parameter $\delta_3$ does not appear in the normal form in Table \ref{tab:bifnormalforms} because there are only four subcases of cusp tangency, which can be represented with just $\delta_1$ and $\delta_2$. However, all of $\delta_1$, $\delta_2$ and $\delta_3$ are needed in the non-degeneracy condition $W[g]\neq 0$ and for separating different subcases.

The derivation of $W[g]$ is as for $K[g]$ in \ref{app:curvature}. Adding coefficients $a,b,\alpha_1,\alpha_2,\beta_1,\beta_2,\gamma$ to the above normal form at bifurcation ($\lambda=0$), with terms of relevant order only, we get that
\begin{equation*}
\begin{array}{r l}
g(x,y,\lambda) = & \delta_1x^3 + \alpha_1 x y_1^2 + \beta_1 x y_2^2 + 2\gamma x y_1 y_2 \\
+ & \alpha_2 x^2y_1 + \beta_2 x^2y_2 + ay_1 + by_2.
\end{array}
\end{equation*}
Rotating the coordinate system to have gradient in the positive $y_1$ direction by the transformation
\begin{equation*}
\left(
\begin{array}{c}
y_1 \\
y_2
\end{array}
\right)
 = \frac{1}{\sqrt{a^2+b^2}}
\left(\begin{array}{cc}
 a & -b \\
 b & a
\end{array}\right)
\left(\begin{array}{c}
\hat{y_1} \\
\hat{y_2}
\end{array}\right)
\end{equation*}
gives that the coefficient of the $x\hat{y_2}^2$ term is
\begin{equation*}
\delta_2 = \frac{1}{a^2+b^2}
\left(\begin{array}{cc}
-b & a
\end{array}\right)
\left(\begin{array}{cc}
\alpha & \gamma \\
\gamma & \beta
\end{array}\right)
\left(\begin{array}{c}
-b \\ a
\end{array}\right)
=
\frac{1}{2|\nabla_y^\perp g|^2}\nabla_y^\perp g^TD^2_y(g_x)\nabla_y^\perp g,
\end{equation*}
and the coefficient of the $x^2\hat{y_2}$ term is
\begin{equation*}
\delta_3 = \frac{1}{\sqrt{a^2+b^2}}(-\alpha_2b + \beta_2 a) = \frac{1}{2\left| \nabla_y^\perp g \right|}\nabla_y g_{xx} \cdot \nabla_y^\perp g.
\end{equation*}
Note that the gradients and the slow Hessian $D_y^2(g_x)$ are expressed in the original slow coordinates $(y_1,y_2)$. Combining the expressions for $\delta_1$, $\delta_2$ and $\delta_3$, we have that
\begin{equation*}
\begin{array}{rl}
W[g] = & 6\delta_1\delta_2 - 2\delta_3 = \\
& 6\frac{g_{xxx}}{6} \frac{1}{2|\nabla_y^\perp g|^2}\nabla_y^\perp g^TD^2_y(g_x)\nabla_y^\perp g - 2\frac{1}{2\left| \nabla_y^\perp g \right|}\nabla_y g_{xx} \cdot \nabla_y^\perp g =  \\
&  \frac{g_{xxx}}{2}\overline{\nabla_y^\perp g^T}D^2_y(g_x)\overline{\nabla_y^\perp g} - \nabla_y g_{xx} \cdot \overline{\nabla_y^\perp g},
\end{array}
\end{equation*}
where $\overline{\nabla_y^\perp g} = \nabla_y^\perp g/|\nabla_y^\perp g|$ is the unit length slow gradient perpendicular. We repeat that if $W[g] <0$ then the cusp tangency is of beaks type, whereas if $W[g] >0$ it is of lips type.

\section{Global codimension one bifurcations for one fast and two slow variables}
\label{app:global12tables}

Tables~\ref{tab:globaldegeneracies12P2Folds},~\ref{tab:globaldegeneracies12P2Cusps} and \ref{tab:globaldegeneracies12P3} list the various inequivalent subclasses of degeneracies $\cD_{4,5,6}$ in Table~\ref{tab:nonpersistentdegeneracies12}. We also include a number of figures which illustrate these degeneracies.

\begin{table}
\caption{Subclasses of special \emph{global} degeneracies for one fast and two slow variables with tangency of fold projection. $P(y)$ is the set of all singular points of the vector field $g$ with slow coordinate $y$. $\nu[g](p)$ is the direction vector of a fold at a point $p$ and $K[g](p)$ is the scalar quadratic curvature of a fold (see the text for details). \emph{f}=fold, \emph{fu}=fold umbra, \emph{fx}=non-interacting fold. For example, fu$\times$fu means that the umbrae of two folds intersect}
{\small
\begin{tabular}{|p{0.18\linewidth}|p{0.54\linewidth}|p{0.16\linewidth}|}
\hline
Aligned umbra-dominant fu$\times$f tangency: &$\cD_{4,1}[g,y]=\{P(y) \subset \cD_4^1[g]: |P(y)|=2 $ and $ U[g](p_1)=p_2 \mbox{ and } \nu[g](p_1) \cdot \nu[g](p_2) > 0,$ and $K[g](p_1)>K[g](p_2)$ for some $ p_1,p_2\in P(y) \}$ &
Fig.~\ref{fig:structureCM121D4Aligned} a,b,c)\\
\hline
Aligned fold-dominant fu$\times$f tangency: &$\cD_{4,2}[g,y]=\{P(y) \subset \cD_4^1[g]: |P(y)|=2 $ and $ U[g](p_1)=p_2 \mbox{ and } \nu[g](p_1) \cdot \nu[g](p_2) > 0$ and $K[g](p_1)<K[g](p_2)$ for some $ p_1,p_2\in P(y) \}$ & Fig.~\ref{fig:structureCM121D4Aligned} d,e,f)\\
\hline
Aligned fu$\times$fu tangency: &$\cD_{4,3}[g,y]=\{P(y) \subset \cD_4^1[g]: |P(y)|=2 $ and $U[g](p_1)=U[g](p_2) \mbox{ and } \nu[g](p_1) \cdot \nu[g](p_2) > 0$ for some $ p_1,p_2\in P(y) \}$ & Fig.~\ref{fig:structureCM121D4Aligned} g,h,i)\\
\hline
Aligned fx$\times$fx tangency: &$\cD_{4,4}[g,y]=\{P(y) \subset \cD_4^1[g]: |P(y)|=2 $ and $\forall p_1,~U[g](p_1)\neq p_2,U[g](p_2) \mbox{ and } \nu[g](p_1) \cdot \nu[g](p_2) > 0$ and $K[g](p_1)<K[g](p_2)$ for some $ p_1,p_2\in P(y) \}$ & Fig.~\ref{fig:structureCM121D4Aligned} j,k,l)\\
\hline
Opposed non-covering fu$\times$f tangency: &$\cD_{4,5}[g,y]=\{P(y) \subset \cD_4^1[g]: |P(y)|=2 $ and $U[g](p_1)=p_2$ and $\nu[g](p_1) \cdot \nu[g](p_2) < 0$ and $K[g](p_1) + K[g](p_2)<0$ and $K[g](p_1)<K[g](p_2)$ for some $ p_1,p_2\in P(y) \}$ & Fig.~\ref{fig:structureCM121D4Opposed} a,b,c)\\
\hline
Opposed non-covering fu$\times$fu tangency: &$\cD_{4,6}[g,y]=\{P(y) \subset \cD_4^1[g]: |P(y)|=2 $ and $U[g](p_1)=U[g](p_2) \mbox{ and } \nu[g](p_1) \cdot \nu[g](p_2) < 0$ and $K[g](p_1)+K[g](p_2)<0,$ for some $ p_1,p_2\in P(y) \}$ & Fig.~\ref{fig:structureCM121D4Opposed} d,e,f)\\
\hline
Opposed non-covering fx$\times$fx tangency: &$\cD_{4,7}[g,y]=\{P(y) \subset \cD_4^1[g]: |P(y)|=2 $ and $ \forall p_1,~U[g](p_1)\neq p_2,U[g](p_2) \mbox{ and } \nu[g](p_1) \cdot \nu[g](p_2) < 0$ and $K[g](p_1)+K[g](p_2)<0,$ for some $ p_1,p_2\in P(y) \}$ & Fig.~\ref{fig:structureCM121D4Opposed} g,h,i)\\
\hline
Opposed covering fu$\times$f tangency: &$\cD_{4,8}[g,y]=\{P(y) \subset \cD_4^1[g]: |P(y)|=2 $ and $U[g](p_1)=p_2 \mbox{ and } \nu[g](p_1) \cdot \nu[g](p_2) < 0$ and $K[g](p_1)+K[g](p_2)>0,$ for some $p_1,p_2\in P(y) \}$ & Fig.~\ref{fig:structureCM121D4Opposed} j,k,l)\\
\hline
Opposed covering fu$\times$fu tangency: &$\cD_{4,9}[g,y]=\{P(y) \subset \cD_4^1[g]: |P(y)|=2 $ and $U[g](p_1)=U[g](p_2) \mbox{ and } \nu[g](p_1) \cdot \nu[g](p_2) < 0$ and $K[g](p_1)+K[g](p_2)>0,$ for some $ p_1,p_2\in P(y) \}$ & Fig.~\ref{fig:structureCM121D4Opposed} m,n,o)\\
\hline
Opposed covering fx$\times$fx tangency: &$\cD_{4,10}[g,y]=\{P(y) \subset \cD_4^1[g]: |P(y)|=2 $ and $ \forall p_1,~U[g](p_1)\neq p_2,U[g](p_2) \mbox{ and } \nu[g](p_1) \cdot \nu[g](p_2) < 0$ and $K[g](p_1)+K[g](p_2)>0,$ for some $ p_1,p_2\in P(y) \}$ & Fig.~\ref{fig:structureCM121D4Opposed} p,q,r)\\
\hline
\end{tabular}
}
\label{tab:globaldegeneracies12P2Folds}
\end{table}

\begin{table}
\caption{Subclasses of special \emph{global} degeneracies for one fast and two slow variables involving the intersection of projections of a fold and cusp. $P(y)$ is the set of all singular points of the vector field $g$ with slow coordinate $y$. $\nu[g](p)$ and $\mu[g](p)$ are direction vectors of folds and cusps respectively (see the text for details). \emph{f}=fold, \emph{fu}=fold umbra, \emph{uc}=unstable cusp, \emph{sc}=stable cusp,\emph{ucu}=unstable cusp umbra, \emph{fx}=non-interacting fold, \emph{cx}=non-interacting cusp. E.g. fu $\times$ sc means that the umbra of a fold intersects a stable cusp}
{\small
\begin{tabular}{|p{0.18\linewidth}|p{0.54\linewidth}|p{0.16\linewidth}|}
\hline
Aligned {fu$\times$sc} intersection: &$\cD_{5,1}[g,y]=\{P(y) \subset \cD_5^1[g]: |P(y)|=2 $ and $ U[g](p_1)=p_2 \mbox{ and } g_{xxx}(p_2) < 0 \mbox{ and } \nu[g](p_1) \cdot \mu[g](p_2) > 0,$ for some $ p_1,p_2\in P(y) \}$ & Fig.~\ref{fig:structureCM121D5Aligned} a,b,c)\\
\hline
Aligned f$\times$ucu intersection: &$\cD_{5,2}[g,y]=\{P(y) \subset \cD_5^1[g]: |P(y)|=2 $ and $ U[g](p_2)\cap p_1 \neq \emptyset \mbox{ and } g_{xxx}(p_2) > 0 \mbox{ and } \nu[g](p_1) \cdot \mu[g](p_2) > 0,$ for some $ p_1,p_2\in P(y) \}$ & Fig.~\ref{fig:structureCM121D5Aligned} d,e,f)\\
\hline
Aligned fu$\times$ucu intersection: &$\cD_{5,3}[g,y]=\{P(y) \subset \cD_5^1[g]: |P(y)|=2 $ and $ U[g](p_2)\cap U[g](p_1) \neq \emptyset \mbox{ and } g_{xxx}(p_2) > 0 \mbox{ and } \nu[g](p_1) \cdot \mu[g](p_2) > 0,$ for some $ p_1,p_2\in P(y) \}$ & Fig.~\ref{fig:structureCM121D5Aligned} g,h,i)\\
\hline
Aligned fx$\times$scx intersection: &$\cD_{5,4}[g,y]=\{P(y) \subset \cD_5^1[g]: |P(y)|=2 $ and $(U[g](p_1) \cup p_1) \cap (U[g](p_2) \cup p_2) = \emptyset \mbox{ and } \nu[g](p_1) \cdot \mu[g](p_2) > 0$ and $g_{xxx}(p_2)<0,$ for some $ p_1,p_2\in P(y) \}$ & Fig.~\ref{fig:structureCM121D5Aligned} j,k,l)\\
\hline
Aligned fx$\times$ucx intersection: &$\cD_{5,5}[g,y]=\{P(y) \subset \cD_5^1[g]: |P(y)|=2 $ and $(U[g](p_1) \cup p_1) \cap (U[g](p_2) \cup p_2) = \emptyset \mbox{ and } \nu[g](p_1) \cdot \mu[g](p_2) > 0$ and $g_{xxx}(p_2)>0,$ for some $ p_1,p_2\in P(y) \}$ & Fig.~\ref{fig:structureCM121D5Aligned} m,n,o) \\
\hline
Opposed fu$\times$sc intersection: &$\cD_{5,6}[g,y]=\{P(y) \subset \cD_5^1[g]: |P(y)|=2 $ and $ U[g](p_1)=p_2 \mbox{ and } g_{xxx}(p_2) < 0 \mbox{ and } \nu[g](p_1) \cdot \mu[g](p_2) < 0,$ for some $ p_1,p_2\in P(y) \}$ & Fig.~\ref{fig:structureCM121D5Opposed} a,b,c)\\
\hline
Opposed f$\times$ucu intersection: &$\cD_{5,7}[g,y]=\{P(y) \subset \cD_5^1[g]: |P(y)|=2 $ and $ U[g](p_2)\cap p_1 \neq \emptyset \mbox{ and } g_{xxx}(p_2) > 0 \mbox{ and } \nu[g](p_1) \cdot \mu[g](p_2) < 0,$ for some $ p_1,p_2\in P(y) \}$ & Fig.~\ref{fig:structureCM121D5Opposed} d,e,f)\\
\hline
Opposed fu$\times$ucu intersection: &$\cD_{5,8}[g,y]=\{P(y) \subset \cD_5^1[g]: |P(y)|=2 $ and $ U[g](p_2)\cap U[g](p_1) \neq \emptyset \mbox{ and } g_{xxx}(p_2) > 0 \mbox{ and } \nu[g](p_1) \cdot \mu[g](p_2) < 0,$ for some $ p_1,p_2\in P(y) \}$ & Fig.~\ref{fig:structureCM121D5Opposed} g,h,i)\\
\hline
Opposed fx$\times$scx intersection: &$\cD_{5,9}[g,y]=\{P(y) \subset \cD_5^1[g]: |P(y)|=2 $ and $(U[g](p_1) \cup p_1) \cap (U[g](p_2) \cup p_2) = \emptyset \mbox{ and } \nu[g](p_1) \cdot \mu[g](p_2) < 0$ and $g_{xxx}(p_2)<0,$ for some $ p_1,p_2\in P(y) \}$ & Fig.~\ref{fig:structureCM121D5Opposed} j,k,l)\\
\hline
Opposed fx$\times$ucx intersection: &$\cD_{5,10}[g,y]=\{P(y) \subset \cD_5^1[g]: |P(y)|=2 $ and $(U[g](p_1) \cup p_1) \cap (U[g](p_2) \cup p_2) = \emptyset \mbox{ and } \nu[g](p_1) \cdot \mu[g](p_2) < 0$ and $g_{xxx}(p_2)>0,$ for some $ p_1,p_2\in P(y) \}$ & Fig.~\ref{fig:structureCM121D5Opposed} m,n,o) \\
\hline
\end{tabular}
}
\label{tab:globaldegeneracies12P2Cusps}
\end{table}

\begin{table}
\caption{Subclasses of special \emph{global} degeneracies for one fast and two slow variables that involve $|P(y)|=3$ singular points. $P(y)$ is the set of all singular points of the vector field $g$ with slow coordinate $y$. The scalars $a_1$,$a_2$ and $a_3$ are coefficients of fold direction vectors. See Section \ref{sec:codimonecritsed12} for details. \emph{f}=fold, \emph{fu}=fold umbra, \emph{fx}=non-interacting fold. E.g. fu$\times$f fu$\times$fu means that one fold umbra intersects another fold, whose umbra interacts with the umbra of another fold}
{\small
\begin{tabular}{|p{0.17\linewidth}|p{0.55\linewidth}|p{0.16\linewidth}|}
\hline
Non-covering fu$\times$f f$\times$fu intersection: &$\cD_{6,1}[g,y]=\{P(y) \subset \cD_6^1[g]: |P(y)|=3 $ and $ U[g](p_1) = p_2 \mbox{ and } U[g](p_2) = p_3 \mbox{ and } a_1 \cdot a_2 \cdot a_3 < 0 $ for some $p_1,p_2,p_3 \in P(y) \}$ & Fig.~\ref{fig:structureCM121D6Sideview} d,e,f)\\
\hline
Covering fu$\times$f f$\times$fu intersection: &$\cD_{6,2}[g,y]=\{P(y) \subset \cD_6^1[g]: |P(y)|=3 $ and $U[g](p_1) = p_2 \mbox{ and } U[g](p_2) = p_3 \mbox{ and } a_1 \cdot a_2 \cdot a_3 > 0 $ for some $p_1,p_2,p_3 \in P(y) \}$ & Fig.~\ref{fig:structureCM121D6Sideview} d,e,f)\\
\hline
Non-covering fu$\times$f fu$\times$fu intersection: &$\cD_{6,3}[g,y]=\{P(y) \subset \cD_6^1[g]: |P(y)|=3 $ and $U[g](p_1) = p_2 \mbox{ and } U[g](p_2) = U(p_3) \mbox{ and } a_1 \cdot a_2 \cdot a_3 < 0 $ for some $p_1,p_2,p_3 \in P(y) \}$ & Fig.~\ref{fig:structureCM121D6Sideview} g,h,i)\\
\hline
Covering fu$\times$f fu$\times$fu intersection: &$\cD_{6,4}[g,y]=\{P(y) \subset \cD_6^1[g]: |P(y)|=3 $ and $U[g](p_1) = p_2 \mbox{ and } U[g](p_2) = U(p_3) \mbox{ and } a_1 \cdot a_2 \cdot a_3 > 0 \}$ & Fig.~\ref{fig:structureCM121D6Sideview} g,h,i)\\
\hline
Non-covering fu$\times$f fx$\times$fx intersection: &$\cD_{6,5}[g,y]=\{P(y) \subset \cD_6^1[g]: |P(y)|=3 $ and $U[g](p_1) = p_2 \mbox{ and } U[g](p_i) \cap (P(y) \cup U[g](P(y))\setminus U[g](p_i)) = \emptyset,~\forall p_i\neq p_1 \mbox{ and } a_1 \cdot a_2 \cdot a_3 < 0 $ for some $p_1,p_2,p_3 \in P(y) \}$ & Fig.~\ref{fig:structureCM121D6Sideview} j,k,l)\\
\hline
Covering fu$\times$f fx$\times$fx intersection: &$\cD_{6,6}[g,y]=\{P(y) \subset \cD_6^1[g]: |P(y)|=3 $ and $U[g](p_1) = p_2 \mbox{ and } U[g](p_i) \cap (P(y) \cup U[g](P(y))\setminus U[g](p_i)) = \emptyset,~\forall p_i\neq p_1 \mbox{ and } a_1 \cdot a_2 \cdot a_3 > 0 $ for some $p_1,p_2,p_3 \in P(y) \}$ & Fig.~\ref{fig:structureCM121D6Sideview} j,k,l)\\
\hline
Non-covering fu$\times$fu fx$\times$fx intersection: &$\cD_{6,7}[g,y]=\{P(y) \subset \cD_6^1[g]: |P(y)|=3 $ and $U[g](p_1) = U[g](p_2) \mbox{ and } U[g](p_i) \cap (P(y) \cup U[g](P(y))\setminus U[g](p_i)) = \emptyset,~\forall p_i\neq p_1 \mbox{ and } a_1 \cdot a_2 \cdot a_3 < 0 $ for some $p_1,p_2,p_3 \in P(y) \}$ & Fig.~\ref{fig:structureCM121D6Sideview} m,n,o)\\
\hline
Covering fu$\times$fu fx$\times$fx intersection: &$\cD_{6,8}[g,y]=\{P(y) \subset \cD_6^1[g]: |P(y)|=3 $ and $U[g](p_1) = U[g](p_2) \mbox{ and } U[g](p_i) \cap (P(y) \cup U[g](P(y))\setminus U[g](p_i)) = \emptyset,~\forall p_i\neq p_1 \mbox{ and } a_1 \cdot a_2 \cdot a_3 > 0 $ for some $p_1,p_2,p_3 \in P(y) \}$ & Fig.~\ref{fig:structureCM121D6Sideview} m,n,o)\\
\hline
Non-covering fx$\times$fx fx$\times$fx intersection: &$\cD_{6,9}[g,y]=\{P(y) \subset \cD_6^1[g]: |P(y)|=3 $ and $U[g](p_i) \cap (P(y) \cup U[g](P(y))\setminus U[g](p_i)) = \emptyset,~\forall p_i \mbox{ and } a_1 \cdot a_2 \cdot a_3 < 0 $ for some $p_1,p_2,p_3 \in P(y) \}$ & Fig.~\ref{fig:structureCM121D6Sideview} p,q,r)\\
\hline
Covering fx$\times$fx fx$\times$fx intersection: &$\cD_{6,10}[g,y]=\{P(y) \subset \cD_6^1[g]: |P(y)|=3 $ and $U[g](p) \cap (P(y) \cup U[g](P(y))\setminus U[g](p)) = \emptyset,~\forall p\in P(y) \mbox{ and } a_1 \cdot a_2 \cdot a_3 > 0 $ for some $p_1,p_2,p_3 \in P(y) \}$ & Fig.~\ref{fig:structureCM121D6Sideview} p,q,r)\\
\hline
\end{tabular}
}
\label{tab:globaldegeneracies12P3}
\end{table}

\begin{figure}
\begin{center}
\includegraphics[width=10cm]{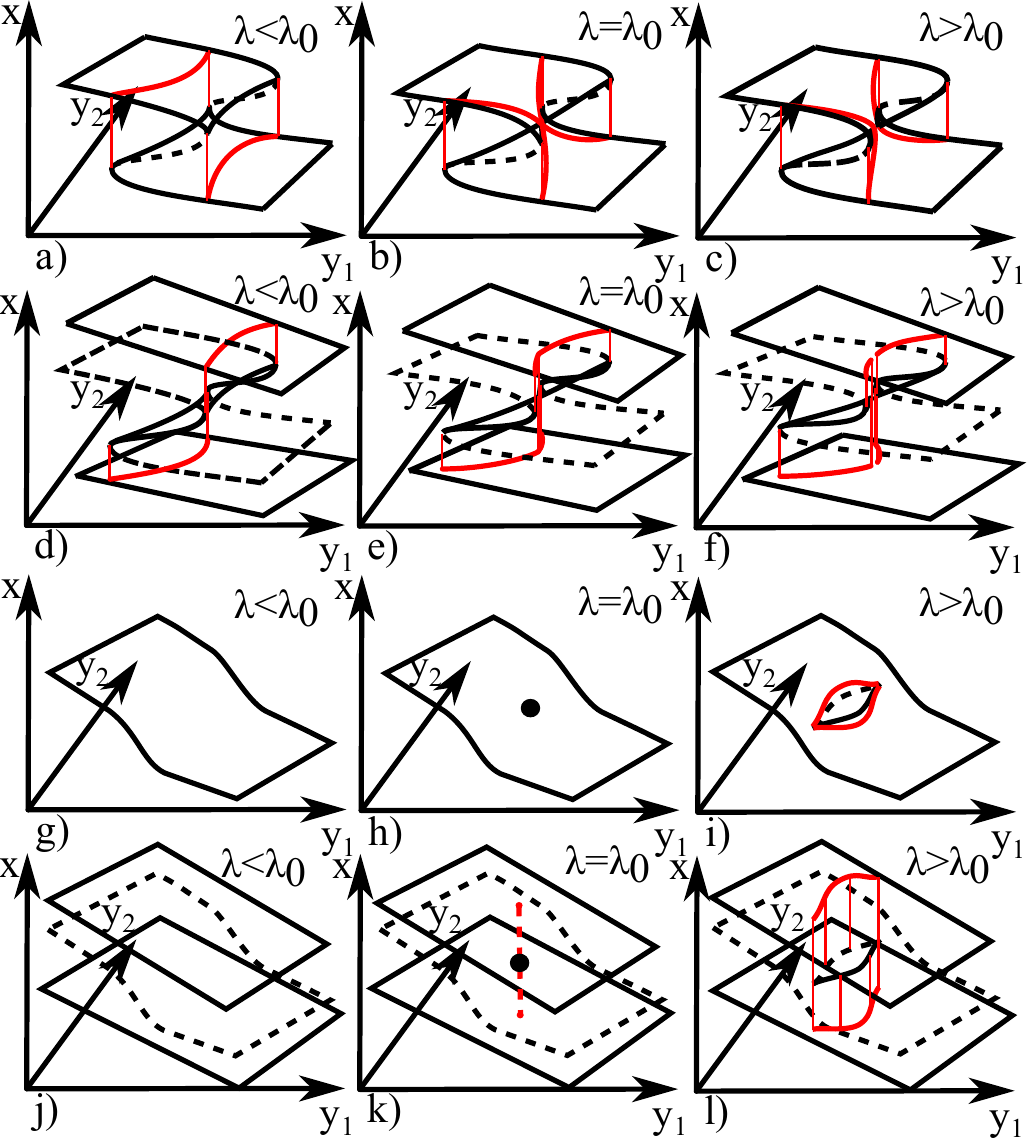}
\end{center}
\caption{
(Color online) Examples of codimension one cusp tangency bifurcation for one fast and two slow variables (see Table \ref{tab:localdegeneracies12}). Each row shows unfolding with a bifurcation parameter $\lambda$. Bifurcation occurs as $\lambda=\lambda_0$. Solid/dashed black lines show stable/unstable sheets of the critical set and red lines show the image of the fold under the umbral map
}
\label{fig:structureCM121D2}
\end{figure}

\begin{figure}
\begin{center}
\includegraphics[width=15cm]{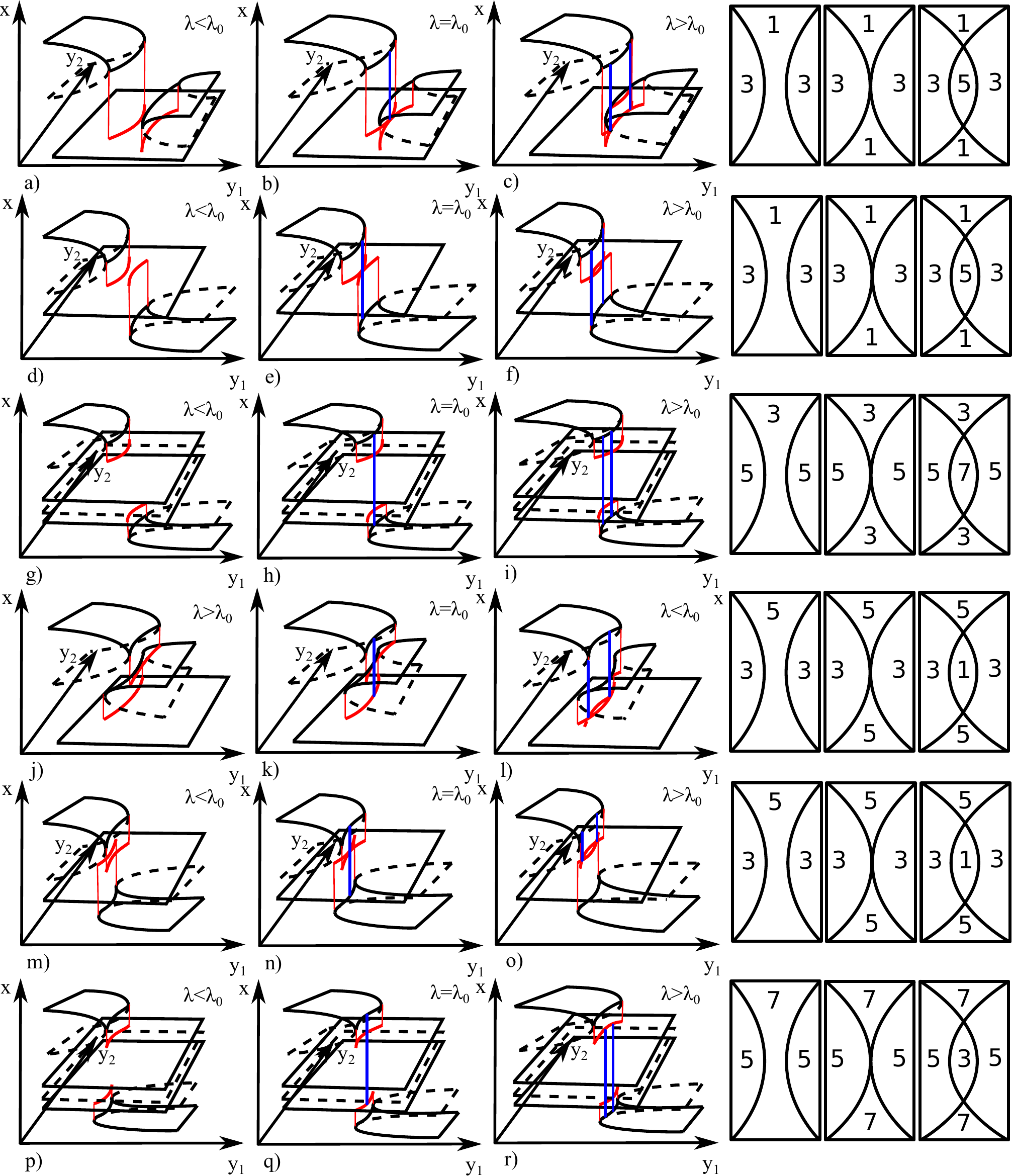}
\end{center}
\caption{
(Color online) Examples of codimension one opposed fold projection tangency bifurcation (see Table \ref{tab:globaldegeneracies12P2Folds}. Each row shows unfolding with a bifucation parameter $\lambda$. Bifurcation occurs as $\lambda=\lambda_0$. Solid/dashed black lines show the stable/unstable sheets of the critical set, red lines show the image of the fold under the umbral map and blue lines indicate tangency of projections of fold sheets. As a visual aid, the number of sheets of the critical set in a neighbourhood of the bifurcation is shown to the right, to be viewed as a projection onto the slow variables
\label{fig:structureCM121D4Opposed}
}
\end{figure}

\begin{figure}
\begin{center}
\includegraphics[width=15cm]{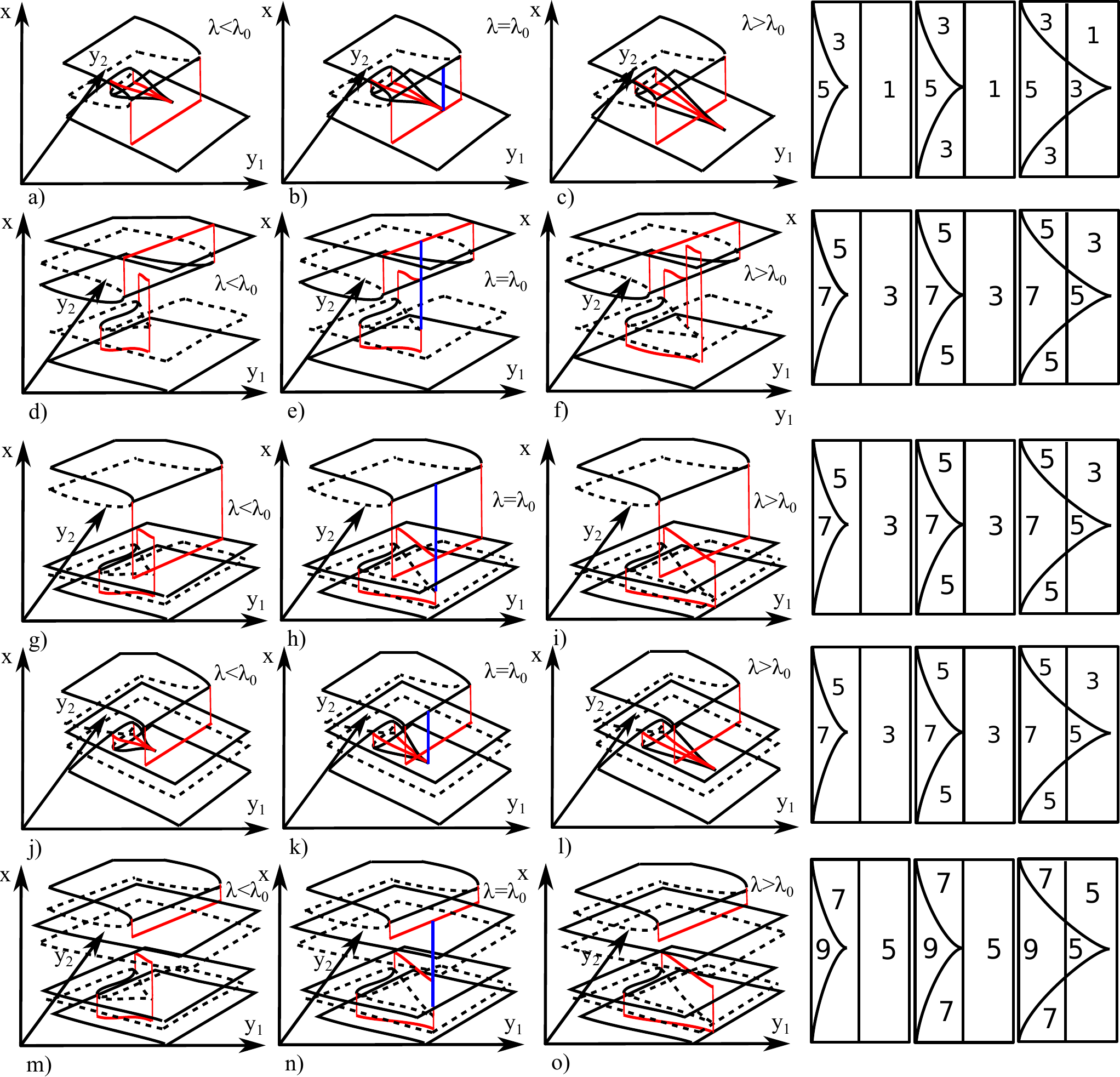}
\end{center}
\caption{
(Color online) Examples of codimension one aligned cusp-fold projection intersection bifurcation (see Table \ref{tab:globaldegeneracies12P2Cusps}). Each row shows unfolding with a bifucation parameter $\lambda$. Bifurcation occurs as $\lambda=\lambda_0$. Solid/dashed black lines show the stable/unstable critical set, red lines show umbrae, and blue lines indicate intersection of projections of cusps and fold onto the slow plane. As a visual aid, the number of sheets of the critical set in a neighbourhood of the bifurcation is shown to the right, to be viewed as a projection onto the slow variables
\label{fig:structureCM121D5Aligned}
}
\end{figure}

\begin{figure}
\begin{center}
\includegraphics[width=15cm]{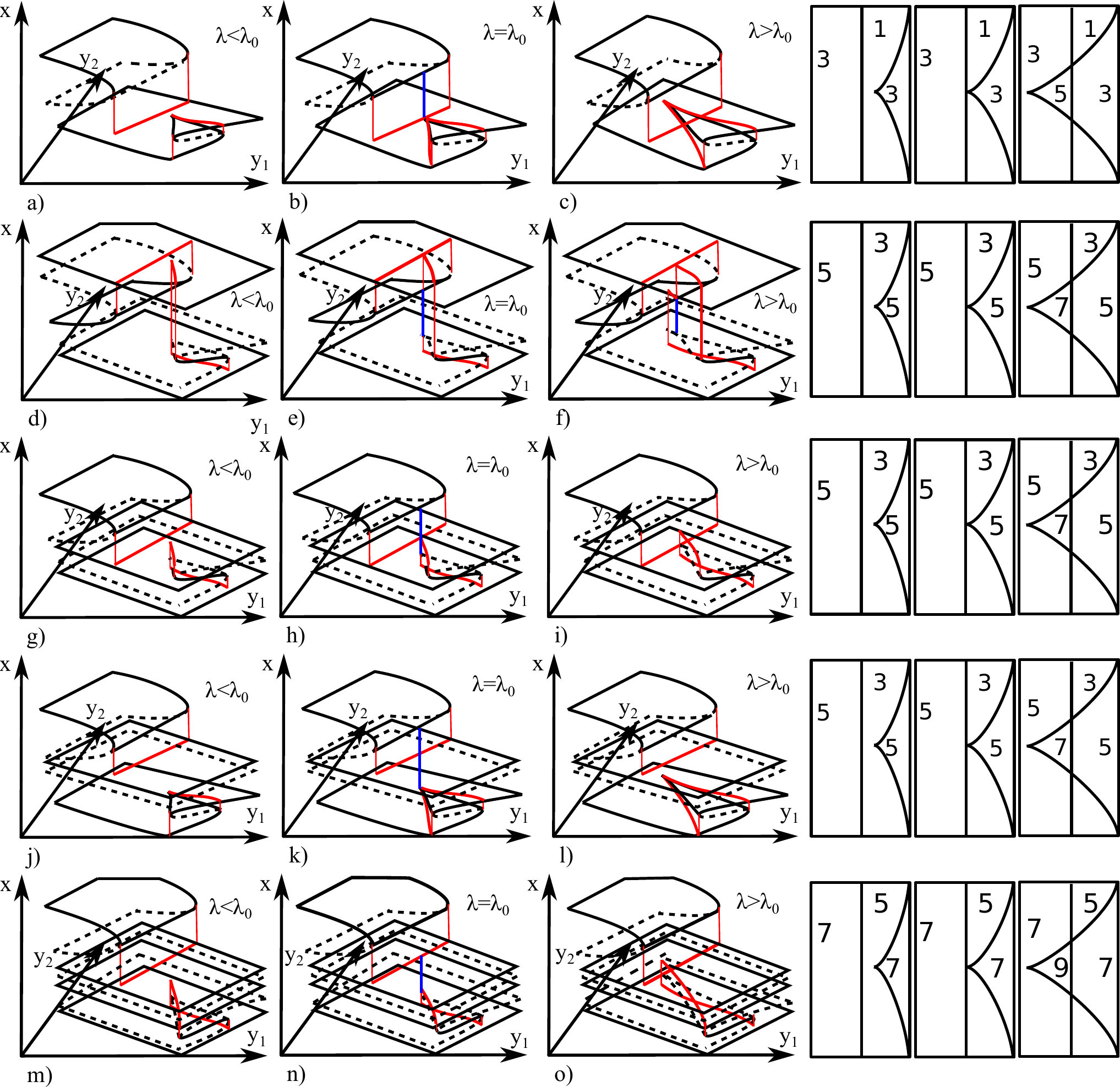}
\end{center}
\caption{
(Color online) Examples of codimension one opposed cusp-fold projection intersection bifurcation (see Table \ref{tab:globaldegeneracies12P2Cusps}). Each row shows unfolding with a bifucation parameter $\lambda$. Bifurcation occurs as $\lambda=\lambda_0$. Solid/dashed black lines show the stable/unstable critical set, red lines show umbrae, and blue lines indicate intersection of projections of cusps and fold onto the slow plane. As a visual aid, the number of sheets of the critical set in a neighbourhood of the bifurcation is shown to the right, to be viewed as a projection onto the slow variables
\label{fig:structureCM121D5Opposed}
}
\end{figure}

\begin{figure}
\begin{center}
\includegraphics[width=15cm]{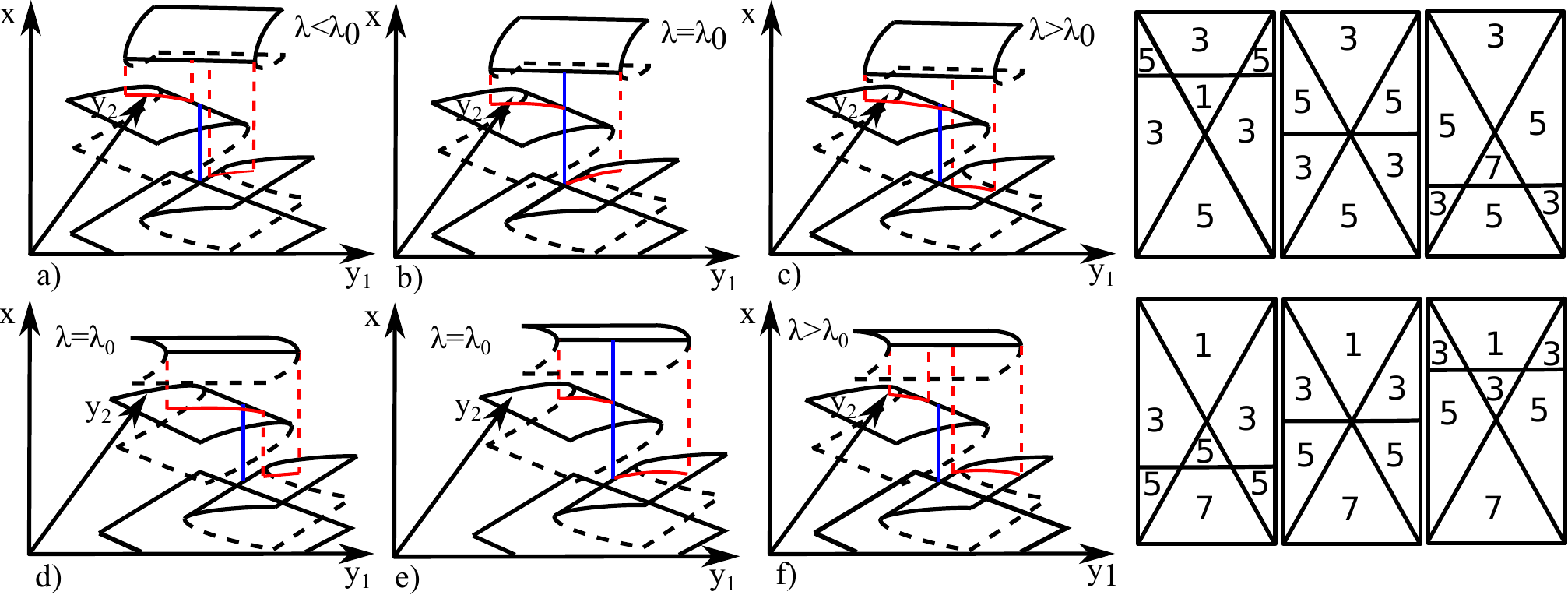}
\end{center}
\caption{
(Color online) Examples of covering a,b,c) and non-covering d,e,f) codimension one triple limit point bifurcation (see Table \ref{tab:globaldegeneracies12P3}). Each row shows unfolding with a bifucation parameter $\lambda$. Bifurcation occurs as $\lambda=\lambda_0$. Solid/dashed black curves show the stable/unstable critical set, red curves show umbrae, and blue curves show intersections of at least the lower two folds, but possibly all three folds. As a visual aid, the number of sheets of the critical set in a neighbourhood of the bifurcation is shown to the right, to be viewed as a projection onto the slow variables
\label{fig:structureCM121D6}
}
\end{figure}
\FloatBarrier

\begin{figure}
\begin{center}
\includegraphics[width=15cm]{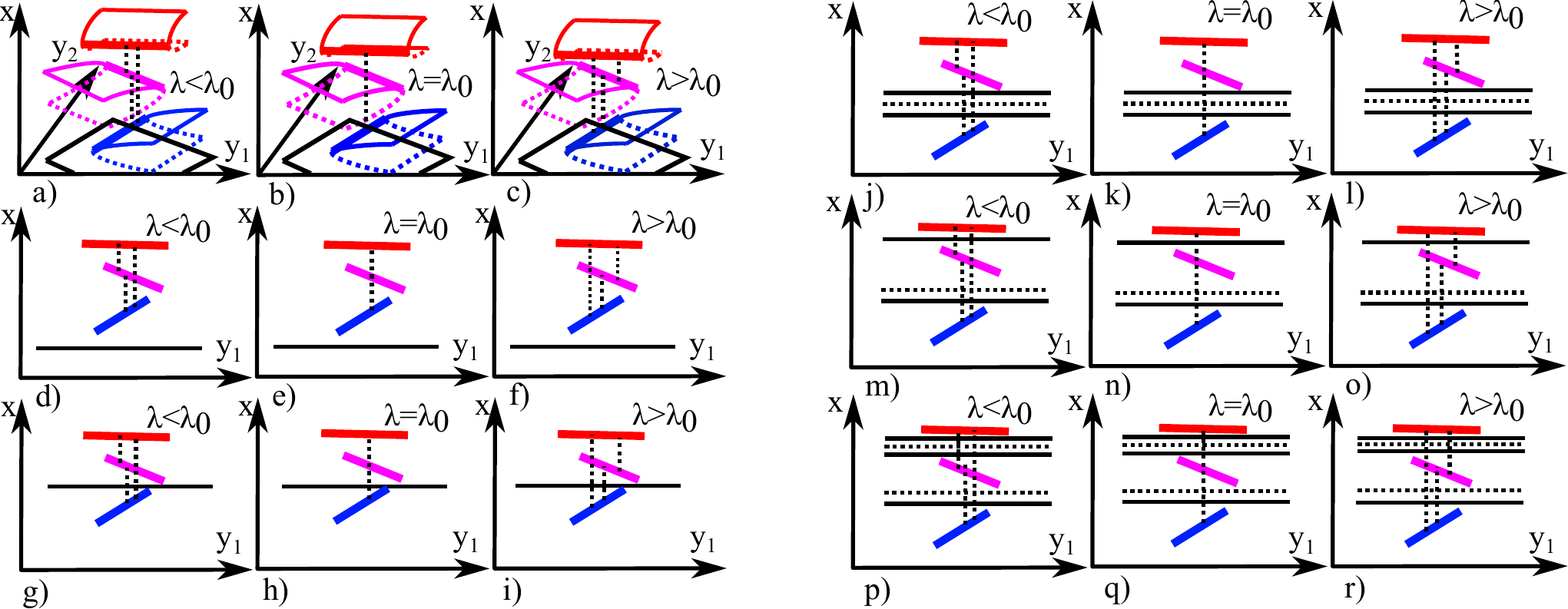}
\end{center}
\caption{
(Color online) Examples of codimension one triple limit point bifurcation (Table~\ref{tab:globaldegeneracies12P3}), to be viewed almost as a projection onto the $(x,y_1)$ plane. Each row shows unfolding with a bifucation parameter $\lambda$. Bifurcation occurs as $\lambda=\lambda_0$. Solid/dashed black horizontal curves show stable/unstable sheets of the critical set, and dashed vertical lines indicate coordinates where slow projections of folds intersect transversally. Red, blue and magenta lines are representations of three fold lines (the row a,b,c) illustrates this representaion for the d,e,f) row). Bifurcation occurs as the slow projections of all three folds intersect. Each case can be either covering or or non-covering, as shown in Fig.~\ref{fig:structureCM121D6}
\label{fig:structureCM121D6Sideview}
}
\end{figure}

\section{Examples of bifurcations of relaxation oscillations for one fast and one slow variable}

\label{sec:examplesection}

In this section we present the equations for the example fast-slow systems for $m=n=1$, showing bifurcations of relaxation oscillations due to the critical manifold in Figure~\ref{fig:slowfastbifexamples}, how they were constructed, and how the figures were produced.

\subsection{Bifurcation of relaxation oscillation due to hyperbolic fold tangency}

We seek fast and slow subsystems $g(x,y)$ and $h(x,y)$ such that (\ref{eq:mainsystem}) displays bifurcation of singular relaxation oscillations due to hyperbolic fold tangency bifurcation (Figure~\ref{fig:slowfastbifexamples} a,b,c)).

The fast subsystem $g(x,y)$ is written as a perturbed product of a hysteresis curve and a circle:
\begin{equation}
\begin{array}{rl}
g_{hyst}(x,y) & = x^3 - 2x + y \\
g_{circ}(x,y) & = (x - \lambda)^2 + (y - y_c)^2 - R^2 \\
g(x,y) & = -(g_{hyst}(x,y)g_{circ}(x,y)+\lambda x + q).
\end{array}
\end{equation}
$(x_c,y_c)=(0.81,-0.25)$ is the centre and $R=0.55$ the radius of the circle, $\lambda$ is the bifurcation parameter and $q=0.01$ is a genericity parameter. For some $\lambda\in[-0.02,0.02]$ tangency bifurcation off the critical set occurs.

The slow subsystem is taken to be
\begin{equation*}
h(x,y) = x - (-b(y-y_c)^2 + x_{max}), 
\end{equation*}
where $b=0.5$ and $x_{max} = x_c+R-0.1$. This choice of $h(x,y)$ makes the nullcline $h(x,y)=0$ intersect the critical set at where it is unstable, for $\lambda \in [-0.02,0.02]$.

In Figure~\ref{fig:slowfastbifexamples} a,b,c) we solve (\ref{eq:mainsystem})
with Matlab's stiff solver \mbox{{\tt ode23s}} for 2000 time units, starting from initial conditions $(x,y)=(2,-1)$ and with scale separation parameter $\epsilon=0.02$.

\subsection{Hysteresis bifurcation of relaxation oscillations}

We seek fast and slow subsystems $g(x,y)$ and $h(x,y)$ such that (\ref{eq:mainsystem}) displays bifurcation of singular relaxation oscillations due hysteresis bifurcation (Figure~\ref{fig:slowfastbifexamples} d,e,f)).

The function $g(x,y)$ is constructed such that the critical set $g(x,y)=0$ at bifurcation $\lambda=0$ is a fifth order polynomial in $x$ as a function of $y$ with quadratic extrema at points $x_1,x_2$ and a cubic double root at $x_3$, that is, $g_x$ satisfies:
\begin{equation*}
g_x(x,y) = a(x-x_1)(x-x_2)(x-x_3)^2-\lambda
\end{equation*}
We fix $x_1 = 1$,$x_3=1$ and leave $a$ and $x_2$ as free parameters. We take $g$ to be the primitive function of $g_x$ with constant term $q=-0.6$ and $y$-term $-y$ such that
\begin{equation*}
g(x,y) = \int a(x-x_1)(x-x_2)(x-x_3)^2dx -\lambda x -q - y,
\end{equation*}
and $g(0,0)=0$. For $\lambda=0$ we solve the linear pair of equations
\begin{equation*}
g(x_1,1)=g(x_3,1) = 0 \\
\end{equation*}
for $a$ and $x_2$, giving $x_2=-13/40$ and $a=640/49$. The system undergoes hysteresis bifurcation for $\lambda\in[-0.2,0.2]$. 

Finally, we reverse the sign of $x$, $x\mapsto-x$, such that
\begin{equation*}
g(x,y) = \int -a(x+x_1)(x+x_2)(x+x_3)^2dx +\lambda x - q - y.
\end{equation*}

The slow subsystem is set to be positive above the constant nullcline $x=x_{nc}=0.7$ and negative below such that
\begin{equation*}
h(x,y) = x - x_{nc}.
\end{equation*}

In Figure~\ref{fig:slowfastbifexamples} d,e,f) we solve (\ref{eq:mainsystem})
with Matlab's stiff solver \texttt{ode23s} for 1000 time units, starting from initial conditions $x=y=0$ and with scale separation parameter $\epsilon=0.05$.

\subsection{Aligned double limit point bifurcation of relaxation oscillations}

We seek fast and slow subsystems $g(x,y)$ and $h(x,y)$ such that (\ref{eq:mainsystem}) displays bifurcation of singular relaxation oscillations due to aligned double limit point bifurcation Figure~\ref{fig:slowfastbifexamples} g,h,i)).

$g(x,y)$ is constructed such that the critical set $g(x,y)=0$ is a fifth order polynomial in $x$ as a function of $y$ with extrema at points $x_1,x_2,x_3$ and $x_4$, that is, $g_x$ satisfies:
\begin{equation*}
g_x(x,y) = a(x-x_1)(x-x_2)(x-x_3)(x-x_4)
\end{equation*}
We fix $x_1 = -1$,$x_3=1/2,x_4=5/4$ and leave $a$ and $x_2$ as free parameters. We take $g$ to be the primitive function of $g_x$ with zero constant term (default of Matlab's \texttt{int} command) such that
\begin{equation*}
g(x,y) = \int a(x-x_1)(x-x_2)(x-x_3)(x-x_4)dx -y,
\end{equation*}
and $g(0,0)=0$. We solve the linear pair of equations
\begin{equation*}
g(x_1,1)=g(x_3,1) = 0 \\
\end{equation*}
for $a$ and $x_2$, giving $x_2=-13/40$ and $a=640/49$. Then we add a bifurcation parameter $\lambda$ breaking the degeneracy, giving
\begin{equation*}
g(x,y) = \int a(x-x_1)(x-x_2)(x-x_3)(x-x_4)dx -\lambda x - y.
\end{equation*}
In Figure~\ref{fig:slowfastbifexamples} g,h,i) $\lambda\in[-0.1,0.1]$. Finally, we reverse the sign of $x$, such that
\begin{equation*}
g(x,y) = \int -a(x+x_1)(x+x_2)(x+x_3)(x+x_4)dx +\lambda x - y.
\end{equation*}

The slow subsystem is set to be positive above the constant nullcline $x=x_{nc}=1.5$ and negative below such that
\begin{equation*}
h(x,y) = x - x_{nc}.
\end{equation*}

In Figure~\ref{fig:slowfastbifexamples} g,h,i) we solve (\ref{eq:mainsystem})
with Matlab's stiff solver \texttt{ode23s} for 1000 time units, starting from initial conditions $x=y=0$ and with scale separation parameter $\epsilon=0.01$.

\subsection{Opposed double limit point bifurcation of relaxation oscillations}

We seek fast and slow subsystems $g(x,y)$ and $h(x,y)$ such that (\ref{eq:mainsystem}) displays bifurcation of singular relaxation oscillations due to opposed double limit point bifurcation (Figure~\ref{fig:slowfastbifexamples} j,k,l)).

We construct the fast subsystem $g(x,y)$ the perturbed product of a hysteresis curve and a "bean" curve \footnote{see \texttt{http://www.2dcurves.com/higher/highergb.html}, accessed 4 August 2019} 
\begin{equation}
\begin{array}{rl}
g_{hyst}(x,y) & = 0.5x^3 - x + y \\
g_{bean,base}(\hat{x},\hat{y}) & = (\hat{x}^2 + \hat{y}^2)^3 - (\hat{x}^2 + (\hat{x}^2 + \hat{y}^2)^2\hat{y}^2)\\
g(x,y) & = -(g_{hyst}(x,y)g_{bean}(\hat{x},\hat{y})+\lambda x + q),
\end{array}
\end{equation}
where $\lambda$ is a bifurcation parameter $\lambda\in[-0.003,0.006]$, $q=0.01$ is a genericity parameter, and $(\hat{x},\hat{y})$ are scaled, rotated and translated coordinates $(x,y)$:
\begin{equation*}
(\hat{x},\hat{y}) = (Mx\cos{\theta} + My\sin{\theta} -x_c,-Mx\sin{\theta}+My\cos{\theta}-y_c),
\end{equation*}
where $M=1.5,\theta=13/40\pi$ and $(x_c,y_c)=(0.97,-0.55)$.

The slow subsystem
\begin{equation*}
h(x,y) = x - (ky + c),
\end{equation*}
with $k = (x_1-x_2)/(y_1-y_2)$, $c = x_1 - ky_1$, $x_1 = 0.7868$, $x_2=1.221$, $y_1 = -0.11$ and $y_2=-0.74$ is chosen to make the nullcline $h(x,y)=0$ pass through the unstable parts of the critical set and enable relaxation oscillation.

In Figure~\ref{fig:slowfastbifexamples} j,k,l) we integrate (\ref{eq:mainsystem}) with Matlab's stiff solver \texttt{ode23s} for $500$ time units, starting from initial conditions $(x_0,y_0) = (2.3,-1)$ and with scale separation parameter $\epsilon=0.001$.

\end{document}